\documentclass{amsart}

\newsymbol\pp 1275
\newsymbol\twoheadrightarrow 1310
\newsymbol\ltimes 226E
\newsymbol\rtimes 226F

\input xy
\xyoption{all}
\usepackage{graphicx}

\newcommand{\Spec}{\operatorname{Spec}}
\newcommand{\Proj}{\operatorname{Proj}}

\title{The combinatorics of quiver representations}
\author{Harm Derksen and Jerzy Weyman}
\thanks{The first author was supported by NSF, grant DMS 0349019
and the second author was supported by NSF, grant  DMS 0300064.}
\newcommand{\dd}{\underline{\operatorname{dim}}\,}
\newcommand{\rk}{\underline{\operatorname{rk}}\,}
\newcommand{\Hom}{\operatorname{Hom}}
\newcommand{\Mat}{\operatorname{Mat}}

\newcommand{\Ext}{\operatorname{Ext}}
\newcommand{\Rep}{\operatorname{Rep}}
\newcommand{\SI}{\operatorname{SI}}
\newcommand{\Gr}{\operatorname{Gr}}
\newcommand{\pperp}{\perp\!\!\!\perp}
\newcommand{\NN}{{\mathbb N}}
\newcommand{\ext}{\operatorname{ext}}

\newcommand{\CC}{{\mathbb C}}

\newcommand{\SL}{\operatorname{SL}}
\newcommand{\GL}{\operatorname{GL}}
\newcommand{\ZZ}{{\mathbb Z}}

\newcommand{\quo}{/\!/}
\newcommand{\supp}{\operatorname{Supp}}
\newcommand{\RR}{{\mathbb R}}
\newtheorem{theorem}{Theorem}[section]
\newtheorem{proposition}[theorem]{Proposition}
\newtheorem{lemma}[theorem]{Lemma}
\newtheorem{corollary}[theorem]{Corollary}

\theoremstyle{definition}
\newtheorem{definition}[theorem]{Definition}

\newtheorem{remark}[theorem]{Remark}
\newtheorem{example}[theorem]{Example}

\begin{document}
\maketitle
\tableofcontents
\section{Introduction}
\subsection{Main results}
Let $\alpha,\beta$ be dimension vectors for a quiver $Q$ without oriented cycles.
If $\langle\alpha,\beta\rangle_Q=0$, where $\langle\cdot,\cdot\rangle_Q$ is
the Euler form (or Ringel form), then we will define a number $\alpha\circ\beta$.
This number can be defined as the  dimension of
a space $\SI(Q,\beta)_{\langle\alpha,\cdot \rangle}$ of  semi-invariants
(see Definition~\ref{defCirc}).
It was shown in \cite{DSW} that $\alpha\circ\beta$ can also be
defined in terms of {\em Schubert calculus}. In the Schubert calculus approach,
$\alpha\circ\beta$ counts the number of $\alpha$-dimensional subrepresentations
of a general $(\alpha+\beta)$-dimensional representation.
In the special case of the triple flag quiver, the number $\alpha\circ\beta$ turns
out to a be Littlewood-Richardson coefficient $c_{\lambda,\mu}^\nu$
where the  partitions $\lambda(\alpha,\beta),\mu=\mu(\alpha,\beta),\nu=\nu(\alpha,\beta)$ depend on $\alpha$ and $\beta$. 
This allows us to prove
many new results about LR-coefficients, and to extend results
about Littlewood-Richardson coefficients to the more general setting of quiver representations.

Knutson and Tao (see~\cite{KT}) proved the saturation conjecture for LR-coefficients:
if $c_{N\lambda,N\mu}^{N\nu}>0$ for some positive
integer $N$, then $c_{\lambda,\mu}^\nu>0$ (see Theorem~\ref{theoSat}).
It was shown in \cite{DW} that the saturation theorem for LR-coefficients generalizes to quivers:
if $N\alpha\circ M\beta>0$ for some positive integers $M,N$, then $\alpha\circ\beta>0$.
Fulton conjectured that $c_{\lambda,\mu}^\nu=1$ implies
that $c_{N\lambda,N\mu}^{N\nu}=1$ for all positive integers $N$.
Knutson, Tao and Woodward
proved  this conjecture in \cite{KTW} (see Theorem~\ref{FultonConjecture}). Belkale gave
a geometric proof of this result (see~\cite{B}).
Belkale's proof generalizes to the
quiver setting; we show that if $\alpha\circ\beta=1$, then $N\alpha\circ M\beta=1$ for
all positive integers $N,M$. 
A proof, following Belkale, is in the Appendix.
The description of Knutson, Tao and Woodward of the walls
of the Klyachko cone in \cite{KTW} (see Theorem~\ref{theoIJK}) generalizes to the quiver setting as well.
In this paper, we extend this result by giving a combinatorial
description of the faces of the cone of {\em arbitrary} dimension (Theorem~\ref{bijection} and 
Corollary~\ref{walls}).
We also prove that, if a triple of partitions $(\lambda,\mu,\nu)$
lies on a wall of the Klyachko cone, then the LR-coefficient
$c_{\lambda,\mu}^\nu$ is a product of smaller LR-coefficients (Theorem~\ref{propprod}).

The main technical tool we introduce in this
paper is the notion of {\em Schur sequences}.
A Schur sequence is a sequence of Schur roots $\alpha_1,\dots,\alpha_s$
such that $\alpha_i\circ\alpha_j=1$ for all $i<j$.
Schur sequences are  a natural generalizations of exceptional sequences,
allowing imaginary Schur roots  to appear instead of only real Schur roots. 
Schur sequences occur naturally as the dimension vectors
appearing in the canonical decomposition of a dimension vector. 
In this paper we also study the dimension vectors of the factors
in a Jordan-H\"older filtration of a $\sigma$-semi-stable
representation. This leads to the notion of the $\sigma$-stable decomposition
of a dimension vector. Again, the dimensions in the $\sigma$-stable
decomposition form a Schur sequence. A crucial result is that every
Schur sequence can be refined to an exceptional sequence. 

Before proving our main results, we 
 review various notions such as perpendicular categories,
exceptional sequences and stability for quivers. 

The main results of the paper have interesting applications.
We 
 use exceptional sequences to ``embed'' the category
of  representation of a quiver $Q$ into the category
of representations of another quiver $Q'$. This is sometimes
possible even when $Q$ is not a subquiver $Q'$. Using this approach
we  construct a triple of partitions $(\lambda,\mu,\nu)$
lying on an extremal ray of the Klyachko cone, such that $c_{\lambda,\mu}^\nu>1$.
In~\cite{CDW}, Calin Chindris and the authors constructed a 
counterexample
to Okounkov's conjecture stating that LR-coefficients are log-concave functions
of the partitions, using this
embedding method. The Embedding Theorem was also used
by the first author to give examples of small Galois groups in Schubert type  problems
(see~\cite{Vakil}[2.10,5.13]).

The current  paper can be viewed as a culmination of the series of papers \cite{DW}, \cite{DW2}, \cite{DW3}, \cite{DSW}. It provides the tools to systematically use the ideas of these papers.

\subsection{Horn's conjecture and related problems}
A classical topic going back to Herman Weyl~\cite{Weyl} is to compare the
eigenvalues of two Hermitian $n\times n$ matrices $A,B$ with the eigenvalues of their sum $C:=A+B$.
For a Hermitian matrix with eigenvalues $\lambda_1\geq\lambda_2\geq \cdots \geq \lambda_n$
define $s(A)=(\lambda_1,\dots,\lambda_n)\in \RR^n$.  One would like
to understand the set
$$
{\mathcal K}_n:=\{(s(A),s(B),s(C))\in \RR^{3n}\mid A,B,C\mbox{ Hermitian, }C=A+B\}.
$$
It turns out that ${\mathcal K}_n$ is given by the {\it trace equation\/}
\begin{equation}\label{Trace}
\lambda_1+\cdots+\lambda_n+\mu_1+\cdots+\mu_n=\nu_1+\cdots+\nu_n
\end{equation}
the {\it weakly decreasing conditions}
\begin{equation}\label{weaklydecreasing}
\lambda_1\geq\cdots\geq \lambda_n\quad,\mu_1\geq\cdots\geq\mu_n,\quad\nu_1\geq\cdots\geq\nu_n.
\end{equation}
and finitely many inequalities of the form
\begin{equation}\label{ineqIJK}
\sum_{i\in I}\lambda_i+\sum_{j\in J}\mu_j\geq \sum_{k\in K}\nu_k
\end{equation}
where $I,J,K$ are subsets of $\{1,2,\dots,n\}$ of the same cardinality. 
We will denote the inequality (\ref{ineqIJK}) by $(\star_{I,J,K})$.
Horn made in 1962
a precise conjecture about triples $(I,J,K)$ for
which the corresponding inequalities $(\star_{I,J,K})$ define ${\mathcal K}_n$ (see~\cite{H}).
Horn's conjecture provides a recursive procedure to determine all those triples $(I,J,K)$.
This conjecture has been proved as a result of work by Klyachko, Totaro, Knutson and Tao
and others.
We will state here a closely related statement about the recursive nature
of the inequalities defining ${\mathcal K}_n$. For 
$$I=\{i_1,i_2,\dots,i_r\}\subseteq\{1,2,\dots,n\}$$
with $i_1<i_2,\dots,i_r$ we define
$$
\lambda(I)=(i_r-r+1,i_{r-1}-r+2,\dots,i_2-1,i_1).
$$
\begin{theorem}\label{theoHorn}
The set
${\mathcal K}_n\subseteq \RR^{3n}$
is given by the trace equation~(\ref{Trace}), the weakly decreasing
conditions~(\ref{weaklydecreasing}) and
all inequalities $(\star_{I,J,K})$ (see~(\ref{ineqIJK}))  with $0<r:=|I|=|J|=|K|<n$, such that
$$
(\lambda(I),\lambda(J),\lambda(K))\in {\mathcal K}_r.
$$
\end{theorem}
The theorem reflects the recursive nature of the cones ${\mathcal K}_n$. 
Once we have determined 
 the cones ${\mathcal K}_1,\dots,{\mathcal K}_{n-1}$, we can
determine a system of inequalities for the cone ${\mathcal K}_n$.

A crucial part of the solution of Horn's conjecture is its connection
to the representation theory of $\GL_n(\CC)$.
Irreducible representations $V_\lambda$ of $\GL_n(\CC)$ are parameterized by
non-increasing integer sequences $\lambda=(\lambda_1,\dots,\lambda_n)$.
The {\it Littlewood-Richardson coefficient\/} $c_{\lambda,\mu}^\nu$ is defined
as the multiplicity of $V_\nu$ inside the tensor product $V_\lambda\otimes V_\mu$, i.e.,
$$
c_{\lambda,\mu}^\nu:=\dim 
(V_\lambda\otimes V_\mu\otimes V_\nu^\star)^{\GL_n(\CC)}.
$$
Here $V_{\nu}^\star$ denotes the dual space of $V_\nu$ and
$(V_\lambda\otimes V_\mu\otimes V_\nu^\star)^{\GL_n(\CC)}$
denotes the $\GL_n(\CC)$-invariant tensors in $V_\lambda\otimes V_\mu\otimes V_\nu^\star$.
We define $c_{\lambda,\mu}^\nu=0$ if $\lambda,\mu,\nu$ are not weakly decreasing.
Let us define
$$
{\mathcal LR}_n=\{(\lambda,\mu,\nu)\in (\ZZ^{n})^3\mid c_{\lambda,\mu}^\nu\neq 0\}.
$$
The following results follow from Klyachko's paper \cite{Kl}.
\begin{theorem}\label{theoKl1}
Let $\RR_+$ be the set of nonnegative real numbers.
The cone
$\RR_+{\mathcal LR}_n\subseteq \RR^{3n}$ is equal to ${\mathcal K}_n$.
\end{theorem}
\begin{theorem}\label{theoKl2}
The set
${\mathcal K}_n\subseteq \RR^{3n}$
is given by the trace equation~(\ref{Trace}), the weakly decreasing conditions~(\ref{weaklydecreasing})
 and all inequalities $(\star_{I,J,K})$ with $0<r:=|I|=|J|=|K|<n$, such that
$$
(\lambda(I),\lambda(J),\lambda(K))\in {\mathcal LR}_r.
$$
\end{theorem}
Finally, the missing link for Theorem~\ref{theoHorn} is proved by Knutson and Tao in~\cite{KT}.
\begin{theorem}[Saturation Theorem]\label{theoSat}
The set ${\mathcal LR}_n\subseteq \ZZ^{3n}$ is saturated, i.e.,
$$
{\mathcal LR}_n=\RR_+{\mathcal LR}_n\cap \ZZ^{3n}.
$$
\end{theorem}
The Saturation Theorem can also be formulated as follows: if $\lambda,\mu,\nu\in \ZZ^{n}$
such that $c_{N\lambda,N\mu}^{N\nu}\neq 0$ for some positive integer
$N$, then $c_{\lambda,\mu}^\nu\neq 0$. 
By Theorems~\ref{theoKl1} and \ref{theoSat} we have ${\mathcal LR}_n={\mathcal K}_n\cap \ZZ^{3n}$.
Combining this with Theorem~\ref{theoKl2} implies Theorem~\ref{theoHorn}.
In~\cite{KT}, Knutson and Tao use their Honeycomb model to prove Theorem~\ref{theoSat}.
See also~\cite{Bu} for another version of the proof.
A geometric proof of Theorem~\ref{theoSat} was given by 
the authors in \cite{DW} and by
Belkale in~\cite{B2}.

By Theorem~\ref{theoKl2}, the set ${\mathcal K}_n$ is defined by (\ref{Trace}), (\ref{weaklydecreasing})
and all inequalities
$(\star_{I,J,K})$ for which $c_{\lambda(I),\lambda(J)}^{\lambda(K)}\neq 0$ are nonzero. 
C.~Woodward was first to note that some of these inequalities are redundant: they follow
from the other inequalities. P.~Belkale proved that all inequalities $(\star_{I,J,K})$
for which $c_{\lambda(I),\lambda(J)}^{\lambda(K)}>1$ are redundant. This
class includes the examples found by Woodward. As the following
theorem by Knutson, Tao and Woodward (\cite{KTW})
shows, none of the remaining inequalities can be omitted.
\begin{theorem}\label{theoIJK}
For $n\geq 3$, ${\mathcal K}_n$ is defined by the equation (\ref{Trace}), 
the inequalities (\ref{weaklydecreasing}) and all inequalities $(\star_{I,J,K})$
for which $c_{\lambda(I),\lambda(J)}^{\lambda(K)}=1$.
None of the inequalities the inequalities can be omitted.
\end{theorem}
For $n=2$, some of the inequalities (\ref{weaklydecreasing}) can be omitted.
The cone ${\mathcal K}_n$ has dimension $3n-1$. The inequality
$$
(\star_{I,J,K}):\sum_{i\in I}\lambda_i+\sum_{j\in J}\mu_j\geq \sum_{k\in K}\nu_k
$$
is necessary if and only if the hyperplane section
$$
\Big\{(\lambda,\mu,\nu)\in (\RR^n)^3\Big|
\sum_{i\in I}\lambda_i+\sum_{j\in J}\mu_j= \sum_{k\in K}\nu_k\Big\}\cap {\mathcal K}_n
$$
defines a {\it facet\/} (or {\it wall\/}) of the cone ${\mathcal K}_n$.

Along the way, Knutson, Tao and Woodward also proved (see~\cite{KTW}) the following Theorem, 
which was conjectured by W.~Fulton.
\begin{theorem}\label{FultonConjecture}
If $c_{\lambda,\mu}^\nu=1$ for some $\lambda,\mu,\nu$, then $c_{N\lambda,N\mu}^{N\nu}=1$
for all nonnegative integers $N$.
\end{theorem}
A geometric proof of Theorem~\ref{FultonConjecture} using Schubert calculus 
was given by P.~Belkale in \cite{B}.

\subsection{The quiver method}
A quiver is just a directed graph. If we attach vector spaces to the vertices and
linear maps to the arrows, we get a representation of that graph. 

Let $Q$ be a quiver without oriented cycles and
$K$ be an algebraically closed field. Suppose that
 $\alpha\in \NN^{Q_0}$, where $\NN=\{0,1,2,\dots\}$, $Q_0$
 is the set of vertices of the quiver and $\NN^{Q_0}$ 
 is the set of dimension vectors.
In \cite{DW} the authors studied the set 
$\Sigma (Q,\alpha)\subseteq \ZZ^{Q_0}$ of weights
occurring in the ring of semi-invariants 
on the space of $\alpha$-dimensional representations $\Rep(Q,\alpha)$
over the field $K$.
We showed that this set is given by one linear homogeneous equation 
and a finite set
of homogeneous linear inequalities. 
Thus the positive real span $\RR_{+}\Sigma(Q,\alpha)\subseteq \RR^{Q_0}$
 forms a rational polyhedral cone in $\RR^{Q_0}$. 
In the particular case of the triple flag quiver $T_{n,n,n}$ and a special dimension vector $\beta$ the cone
$\RR_+\Sigma (T_{n,n,n},\beta )$ turns out to be equal to the cone $\RR_+{\mathcal LR}_n$ defined above. 
Thus it is interesting to study the properties of the general cones $\RR_+\Sigma (Q,\alpha )$, in particular whether  the statements listed above concerning the cones $\RR_+{\mathcal LR}_n$ generalize to arbitrary quivers. If this is true, one could expect that such technique would allow to prove stronger results about the cones $\RR_+{\mathcal LR}_n$.

This program is carried out in the present paper for arbitrary quiver $Q$ without oriented cycles. The approach is based on studying the notions of semi-stable filtrations from Geometric Invariant Theory in terms of quiver representations. As a result of this study we generalize Theorems~\ref{theoSat}, \ref{theoIJK} and \ref{FultonConjecture}  above to arbitrary cones $\RR_+\Sigma(Q,\alpha )$. We also obtain combinatorial description of faces of all codimensions in the cones $\RR_+\Sigma(Q,\alpha )$. For the case of faces of codimension one this is equivalent to Theorem~\ref{theoIJK}. In the special case of extremal rays the results imply that the semi-invariants with weights lying on an extremal ray of $\RR_+\Sigma (Q,\alpha )$ form a subring isomorphic to the ring of semi-invariants for quivers with two vertices and multiple arrows. 
This indicates that one needs to study {\it all\/} quivers, not just the special class
of triple flag quivers,
so the quiver technique is the right tool for studying similar questions. 
The cones $\RR_+{\mathcal LR}_n$ are also related to a problem of existence of short exact
sequences of abelian $p$-groups (see~\cite{F2,Klein}).  Recent results of C.~Chindris (\cite{CH1,CH2}) show that quivers can be successfully applied to get similar results for longer exact sequences.

We use the developed techniques to prove some new results on the faces of cones $\RR_+{\mathcal LR}_n$, in particular a product formula for Littlewood-Richardson coefficients. The general approach explains why one could expect such formula.

\subsection{Organization}
The paper is organized as follows. In Section~\ref{sec2} we give the basic notation and review the needed results on quivers and their semi-invariants. In particular we review the content of papers \cite{CB},  \cite{DW}, \cite{S2} and \cite{S1}, in particular the Saturation Theorem for the cones $\RR_+\Sigma(Q,\alpha )$, the exceptional sequences and orthogonal categories. 

There are some reformulations and extensions, notably the Embedding Theorem~\ref{theoembedding}. 
Schofield's technique of orthogonal categories allows to embed the category of representations
of a smaller quiver (not necessarily a subquiver) into the category of representations
of the original quiver. Theorem~\ref{theoembedding} shows that this embedding respects
semi-invariants. This gives us a method for proving results about the cones
$\RR_+\Sigma(Q,\alpha)$ by induction. This technique  requires one to work with arbitrary quivers,
not just triple flag quivers.

We also formulate the Generalized Fulton's Conjecture (Theorem~\ref{thmBelkale}) whose proof 
(essentially due to P.~Belkale, see~\cite{B}) is given in the Appendix.

In Section~\ref{sec3} we study the notions of semi-stability and stability of quiver representations. We relate the geometric notion of $(\sigma :\tau)$-stability and the algebraic notion of $\sigma$-stability, using the approach of \cite{R}. We introduce the notions of Harder-Narasimhan, Jordan-H\"older filtrations, and their combination - the HNJH filtrations. The key statement is Lemma~\ref{lemlocallyclosed} which shows that the subsets of the representation space $\Rep(Q,\alpha )$ where the  dimensions of the factors of Harder-Narasimhan and HNJH filtrations are constant are locally closed. Then, filtrations are used to define
the $\sigma$-stable decompositions and $(\sigma :\tau )$-stable decompositions of  dimension vectors, and to prove their basic properties. In particular, Theorem~\ref{invar} relates the semi-invariants in weights $m\sigma$ for a dimension vector $\alpha$
to those for the dimension vectors that are factors in $\sigma$-stable decomposition of $\alpha$.

In Section~\ref{sec4} we introduce the key notions of Schur sequences and Schur quiver sequences, inspired by the notion of exceptional sequences. These notions provide the tools for the descriptions of the faces of the cones $\RR_+\Sigma (Q,\alpha )$. Then we prove the 
Refinement Theorem (Theorem~\ref{theorefinement}) which says that every Schur sequence can be refined to an exceptional sequence. This theorem makes it easy to understand Schur sequences in terms of exceptional sequences.

In Section~\ref{sec5} we prove the basic Theorem~\ref{bijection} giving a bijection between the faces of dimension $n-r$ in $\Sigma (Q,\alpha )$, and Schur quiver sequences of $r$ dimension vectors summing to $\alpha$. This is the main result of the paper. Then we draw consequences for faces of codimension $1$ and for extremal rays.

In Section~\ref{sec6} we study the dual problem of how the $\sigma$-stable decomposition of $\alpha$ varies when $\alpha$ varies and $\sigma$ is fixed. Theorem~\ref{sigmastable} gives a general combinatorial criterion of when $\alpha$ is $\sigma$-stable.
We also extend the notion of $\sigma$-stable decomposition to quivers with oriented cycles.

In Section~\ref{sec7} we apply our results to triple flag quivers and the cones $\RR_+{\mathcal LR}_n$. We recover Theorem~\ref{theoIJK} of Knutson, Tao and Woodward  on faces of codimension one.

We also investigate the faces of $\RR_+{\mathcal LR}_n$ of arbitrary codimension. In particular we show that for $n\le 7$ all weight multiplicities along extremal rays of the cones $\RR_+{\mathcal LR}_n$ are equal to $1$, and for $n=8$ we give an example of an extremal ray with weight multiplicities bigger than $1$.

Finally we prove the product formula for Littlewood-Richardson coefficients (Theorem~\ref{propprod}). It shows that if a weight $\sigma$ corresponding to a triple of partitions $(\lambda ,\mu,\nu )$ lies on a face of $\Sigma (T_{n,n,n},\beta )$ of positive codimension, then the Littlewood-Richardson coefficient $c_{\lambda ,\mu}^{\nu}$ decomposes to a product of smaller Littlewood-Richardson coefficients.\medskip

\noindent{\bf Acknowledgment.} Both authors would like to thank Prakash Belkale for the permission to include his proof of the Fulton Conjecture, Calin Chindris for helpful discussions, and Jiarui Fei for proofreading.

\section{Preliminaries}\label{sec2}
\subsection{Basic notions for quivers}
A quiver $Q$ is a quadruple $Q= (Q_0 ,Q_1,h,t)$ where 
 $Q_0$ is a finite set of vertices,
 $Q_1$ is a finite set of arrows and $h,t:Q_1\to Q_0$ are maps.
For each arrow $a\in Q_1$, its
head is $ha:=h(a)\in Q_0$ and its tail is $ta:=t(a)\in Q_1$. 

We fix an algebraically closed field $K$. 
A representation $V$ of $Q$ is a family of
 finite dimensional
$K$-vector spaces 
$$\lbrace V(x) \mid x\in Q_0\rbrace$$ 
together with a family of
$K$-linear maps
$$
\lbrace V(a):V(ta)\rightarrow V(ha)\mid a\in Q_1\rbrace.
$$ 
The dimension vector of a representation
$V$ is the function $\dd V:Q_0\rightarrow \NN$ defined by
$$(\dd V)(x):= \dim V(x),\quad x\in Q_0.$$

The dimension vectors $\NN^{Q_0}$ 
are contained in the set
$\Gamma:=\ZZ^{Q_0}$ of integer-valued functions
on $Q_0$. 
A morphism $\phi :V\rightarrow W$ between two representations is the collection of linear maps
$$
\{\phi(x): V(x)\rightarrow W(x)\mid x\in Q_0\}$$ 
such that for each $a\in Q_1$ we have
$$W(a)\phi(ta)=\phi (ha)V(a).$$
We denote the (finite dimensional) linear space of morphisms from $V$ to $W$ by $\Hom_Q(V,W)$.

A {\it nontrivial} path $p$ in the quiver $Q$ of length $n\geq 1$
is a sequence $p=a_na_{n-1}\cdots a_1$ of arrows, such that
$ta_{i+1}=ha_i$ for $i=1,2,\dots,n-1$. We define the {\it head\/} and
the {\it tail\/} of the path $p$ as $hp:=ha_n$ and $tp:=ta_1$
respectively. Besides the nontrivial paths one usually also defines
a trivial path $\epsilon_x$ of length 0 for every vertex $x\in Q_0$
whose head and tail are equal to $h\epsilon_x=t\epsilon_x=x$.
An oriented cycle is a nontrivial path $p$ satisfying $hp=tp$.
{\it Throughout this paper, we will assume that $Q$ has no oriented cycles,
unless stated otherwise.}

The category $\Rep(Q)=\Rep_K(Q)$ 
of representations of $Q$ over $K$ is hereditary, i.e., 
the subobject of a projective
object is projective. This means that every representation has projective 
dimension $\le 1$, i.e., $\Ext^i_Q(V,W)=0$ for all representations $V,W$ and all $i>1$.
\begin{lemma}[See~\cite{Ri}]\label{lemkercoker}
The spaces $\Hom_Q (V, W)$ and $\Ext_Q (V, W):=\Ext^1_Q(V,W)$ are the kernel
and cokernel of the following linear map
\begin{equation}
d^V_W : \bigoplus_{x\in Q_0} \Hom (V(x), W(x))\longrightarrow\bigoplus_{a\in Q_1} 
\Hom (V(ta), W(ha))\label{eq1}
\end{equation}
where  $d^V_W$ is given  by
$$
\{\phi(x)\mid x\in Q_0\}\mapsto
\{W(a)\phi(ta)-\phi(ha)V(a)\mid 
a\in Q_1\}.
$$
\end{lemma}
Let $\alpha ,\beta$ be two elements of $\Gamma$. 
We define the {\it Euler inner product\/} (or {\it Ringel form\/}) by
\begin{equation}
\langle\alpha ,\beta \rangle_Q = \sum_{x\in Q_0} \alpha (x)\beta (x)-
\sum_{a\in Q_1} \alpha(ta)\beta(ha).\label{eq2}
\end{equation}
It follows from Lemma~\ref{lemkercoker} and (\ref{eq2}) that 
\begin{equation}\label{eq41}
\langle\dd V,\dd W\rangle_Q= \dim\Hom_Q (V,W)-\dim\Ext_Q (V,W).
\end{equation}
We will omit the subscript $Q$ and just write $\langle\cdot,\cdot\rangle$
instead of $\langle\cdot,\cdot\rangle_Q$ if there is no chance of confusion.

\subsection{Semi-invariants for quiver representations}
Suppose that $V$
is a $\beta$-dimensional representation. Choose
a basis of each of the vector spaces $V(x)$, $x\in Q_0$.
The matrix of $V(a)$ with respect to the bases
in $V(ta)$ and $V(ha)$ is an element
of $\Mat_{\beta(ha),\beta(ta)}(K)$,
where $\Mat_{p,q}(K)$ denotes the $p\times q$ matrices with
entries in $K$. This way we can associate to $V$
an element of the representation space
$$
\Rep(Q,\beta):=\prod_{a\in Q_1}\Mat_{\beta(ha),\beta(ta)}(K)
$$ 
 The group
$$\GL(Q,\beta):=\prod_{x\in Q_0} \GL_{\beta (x)}(K)$$
acts on $\Rep(Q,\beta)$ as follows
$$
\{A(x)\mid x\in Q_0\}\cdot \{V(a)\mid a\in Q_1\}:=
\{A(ha)V(a)A(ta)^{-1}\mid a\in Q_1\},
$$
for $\{A(x)\mid x\in Q_0\}\in \GL(Q,\beta)$ and
$\{V(a)\mid a\in Q_1\}\in \Rep(Q,\beta)$. 
The action of $\GL(Q,\beta)$
 on $\Rep(Q,\beta)$ corresponds to base changes in
 each of the vector
spaces $V(x)\cong K^{\beta(x)}$, $x\in Q_0$. The orbits
of $\GL(Q,\beta)$ in $\Rep(Q,\beta)$ correspond to the isomorphism
classes of $\beta$-dimensional representations of $Q$.

Define $\SL(Q,\beta)\subseteq \GL(Q,\beta)$ by
$$\SL(Q,\beta)=\prod_{x\in Q_0} \SL_{\beta (x)}(K).$$
We are interested in the rings of semi-invariants
$$
\SI(Q,\beta) = K[\Rep (Q,\beta)]^{\SL(Q,\beta)}.
$$
The ring $\SI(Q,\beta)$ has a weight space decomposition
$$
\SI(Q,\beta)=\bigoplus_\sigma \SI(Q,\beta)_\sigma
$$
where $\sigma$ runs through the multiplicative characters of $\GL(Q,\beta)$ and 
$$
\SI(Q,\beta)_\sigma = 
\lbrace f\in K[\Rep(Q,\beta)]\mid g(f)=\sigma (g)f\ \forall g\in 
\GL(Q,\beta)\rbrace.
$$
Any {\it character\/} or {\it weight\/} of $\GL(Q,\beta)$ has the form
\begin{equation}\label{eqchar}
\{A(x)\mid x\in Q_0\}\mapsto \prod_{x\in Q_0} (\det A(x))^{\sigma(x)}
\end{equation}
with $\sigma(x)\in \ZZ$ for all $x\in Q_0$. This way, the multiplicative
character (\ref{eqchar}) of $\GL(Q,\beta)$ can be identified with $\sigma\in \Gamma:=\ZZ^{Q_0}$.

If $\alpha\in \Gamma$
 then we define
$$
\sigma(\alpha)=\sum_{x\in Q_0} \sigma(x)\alpha(x).
$$
We will identify the set of weights  with 
$\Gamma^\star= \Hom_{\ZZ}(\Gamma,\ZZ)\cong \ZZ^{Q_0}$.
Note that $\Gamma^\star$ and $\Gamma$ are canonically isomorphic, but we still
would like to distinguish between  them.

Let us choose the dimension
vectors $\alpha,\beta\in \NN^{Q_0}$ such that
$\langle \alpha ,\beta \rangle=0$. 
If $V\in \Rep(Q,\alpha)$ and $W\in \Rep(Q,\beta)$, then
 the matrix of $d^V_W$ in (\ref{eq1}) is a square
matrix. Following \cite{S2} we  define the semi-invariant
$$
c(V,W) := \det d^V_W
$$
of the action of $\GL(Q,\alpha)\times \GL(Q,\beta)$ on $\Rep(Q,\alpha)\times \Rep(Q,\beta)$.
For a fixed $V$ the restriction of $c$ to $\lbrace V\rbrace\times \Rep(Q,\beta)$
defines a semi-invariant $c^V$ in $\SI(Q, \beta)$. 
Schofield proved (\cite[Lemma 1.4]{S2}) that the weight of $c^V$ equals 
$\langle\alpha , \cdot\rangle$. Note that $\langle\alpha,\cdot\rangle$
can be  viewed as an element in $\Gamma^\star$.
Similarly, for a fixed
$W$ the restriction of $c$ to $\Rep(Q,\alpha)\times\lbrace W\rbrace$
defines a semi-invariant $c_W$ in $\SI(Q, \alpha)$ of weight 
$-\langle\cdot, \beta\rangle$
(\cite[Lemma 1.4]{S2}).
\begin{lemma}[Lemma 1 of \cite{DW}]\label{lemexseq}
Suppose that 
$$
0\to V_1\to V\to V_2\to 0
$$
is an exact sequence of representations of $Q$ and
$\langle \dd V_1,\beta\rangle=\langle\dd V_2,\beta\rangle=0$,
then as a function on $\Rep(Q,\beta)$, $c^V$ is, up to a scalar, equal
to $c^{V_1}\cdot c^{V_2}$.
\end{lemma}

\begin{theorem}[Theorem 1 of~\cite{DW}]\label{theoDW}
The ring $\SI(Q,\beta)$ is spanned by semi-invariants of the form $c^V$
for which $\langle \dd V,\beta\rangle=0$.
It is also spanned by semi-invariants of the form $c_W$ for which
$\langle \beta,\dd W\rangle=0$.
\end{theorem}
For a more general statement for quivers with oriented cycles,
see~\cite{DW3,SvdB}.
\begin{remark}\label{rem1}
If $\langle \dd V,\dd W\rangle=0$ then we have  
$c(V,W)=c^V(W)=c_W(V)=0$ 
if and only if
$\Hom_Q (V, W)\ne 0$ which is equivalent to $\Ext_Q(V, W)\ne 0$
by Lemma~\ref{lemkercoker}.
\end{remark}
It was also shown in \cite{DW} that 
$$\dim \SI(Q,\beta)_{\langle \alpha,\cdot\rangle}=
\dim \SI(Q,\alpha)_{-\langle \cdot,\beta\rangle}.$$ 
\begin{definition}\label{defCirc}
For dimension vectors $\alpha,\beta$ with $\langle \alpha,\beta\rangle=0$
we define
$$
(\alpha \circ \beta)_Q:=\dim \SI(Q,\beta)_{\langle \alpha,\cdot\rangle}=
\dim \SI(Q,\alpha)_{-\langle \cdot,\beta\rangle}.
$$
Again, we will drop the subscript $Q$ most of the time and write $\alpha\circ \beta$
instead of $(\alpha\circ\beta)_Q$.
\end{definition}
\subsection{Representations in general position}
A representation $V$ is called indecomposable if it is
not isomorphic to the direct sum of two nonzero representations.
The set of dimension vectors $\alpha$ for which there exists
an $\alpha$-dimensional indecomposable representation can be
identified with the set of positive roots for the Kac-Moody algebra
associated with the graph $Q$ (where we forget the orientation).
This was proved in \cite{K1}. We will call a dimension vector $\alpha$
a {\it root\/} if there exists an indecomposable representation of dimension $\alpha$.
Kac proved that $\langle \alpha,\alpha\rangle\leq 1$ for every
root $\alpha$. If $\alpha$ is a root, then we will call $\alpha$ {\it real\/}
if $\langle\alpha,\alpha\rangle=1$ and {\it imaginary\/} if
$\langle\alpha,\alpha\rangle\leq 0$. We call $\alpha$ {\it isotropic\/} if
$\langle\alpha,\alpha\rangle=0$.

A representation $V$ is called a Schur representation (or a {\it brick\/}) if
$\Hom_Q(V,V)\cong K$. Note that every Schur representation must
be indecomposable. If $\Rep(Q,\alpha)$ contains a Schur
representation, then $\alpha$ is called a {\it Schur root}.

A representation $V$ is called in general position of dimension $\alpha$ if
$V\in \Rep(Q,\alpha)$ lies in a sufficiently small  Zariski open
subset (``sufficient'' here depends on the context). 
Suppose that $\alpha$ is a Schur root.
Since $V\mapsto \dim \Hom_Q(V,V)$ depends upper semi-continuously
on $V\in \Rep(Q,\alpha)$, its minimal value 1 is attained
on some open dense subset $U\subseteq \Rep(Q,\alpha)$.
This shows that a general representation of $\Rep(Q,\alpha)$ is indecomposable.
Conversely, if a general representation of dimension $\alpha$ is
indecomposable, then $\alpha$ must be a Schur root (see~\cite[Proposition
1]{K2}). 

We define
$$\hom_Q(\alpha,\beta)=\min\{\dim \Hom_Q(V,W)\mid V\in \Rep(Q,\alpha),W\in
\Rep(Q,\beta)\},
$$
where $\min$ denotes the minimum.
The function $(V,W)\mapsto \dim \Hom_Q(V,W)$ is upper semi-continuous, so
$\dim \Hom_Q(V,W)=\hom_Q(\alpha,\beta)$
if $(V,W)\in \Rep(Q,\alpha)\times \Rep(Q,\beta)$ is in general position
(see~\cite{S1}).
Similarly, we define 
$$\ext_Q(\alpha,\beta)=\min\{\dim \Ext_Q(V,W)\mid V\in \Rep(Q,\alpha),W\in
\Rep(Q,\beta)\}.
$$
We will drop the subscript and
write $\hom(\alpha,\beta)$ and $\ext(\alpha,\beta)$ if there is no confusion.
From (\ref{eq2}) follows that
\begin{equation}\label{eq552}
\langle \alpha,\beta\rangle=\hom(\alpha,\beta)-\ext(\alpha,\beta).
\end{equation}

\begin{definition}\label{defperp}
If
$\hom_Q(\alpha,\beta)=\ext_Q(\alpha,\beta)=0$, then we write $\alpha\perp \beta$ and
we will say that $\alpha$ is left perpendicular to $\beta$. 
\end{definition}
By Remark~\ref{rem1}
we have $\alpha\perp\beta$ if and only if $\alpha\circ\beta\neq 0$.
Following Schofield, we write $\alpha\hookrightarrow\beta$ if a general representation of
dimension $\beta$ contains a subrepresentation of dimension $\alpha$.
We write $\alpha\twoheadrightarrow\beta$ if a general representation
of dimension $\alpha$ has a factor of dimension $\beta$.
The proof of the following Theorem can be found in  \cite{S1} for a base field of characteristic 0.
For a proof that also works in positive characteristic, see~\cite{CB2}.

\begin{theorem}[Theorem 3.3 of~\cite{S1}]\label{lemScho}
We have
$$
\alpha\hookrightarrow \alpha+\beta\quad\Leftrightarrow\quad
\ext_Q(\alpha,\beta)=0\quad(\Leftrightarrow \alpha+\beta\twoheadrightarrow
\beta).
$$
\end{theorem}
\begin{definition}
For a dimension vector $\beta$, we define 
$$\Sigma(Q,\beta)=\{\sigma\in \Gamma^\star\mid \SI(Q,\beta)_\sigma\neq 0\}.
$$
\end{definition}
\begin{theorem}[see~\cite{DW}]\label{theosubrep}
We have
$$\Sigma(Q,\beta)=\{\sigma\in \Gamma^\star\mid \sigma(\beta)=0\mbox{ and 
$\sigma(\gamma)\leq 0$ for all $\gamma\hookrightarrow \beta$}\}.
$$
\end{theorem}
For some $\gamma\hookrightarrow\beta$, the inequality
$\sigma(\gamma)\leq 0$ can be omitted because it follows
from the other inequalities. Later, we will describe a minimal
list of inequalities for $\Sigma(Q,\beta)$.
\begin{theorem}[see~\cite{DSW}]\label{theoDSW}
Suppose that $\alpha,\beta$ are dimension vectors satisfying
$\hom(\alpha,\beta)=\ext(\alpha,\beta)=0$. Then a general
representation of dimension $\alpha+\beta$ has exactly $\alpha\circ\beta$
subrepresentations of dimension $\alpha$.
\end{theorem}
\begin{lemma}\label{lemfiberlargest}
Under the assumptions of Theorem~\ref{theosubrep}, if $V\in \Rep(Q,\alpha+\beta)$
is arbitrary such that $V$ has exactly $r$ subrepresentations, where $r$ is finite,
then $r\leq \alpha\circ\beta$.
\end{lemma}
\begin{proof}
Schofield constructs a variety $Z:=R(Q,\alpha\subset \alpha+\beta)$ (see~\cite{S1})
and a projective morphism $p:Z\to \Rep(Q,\alpha+\beta)$ such that the
fiber $p^{-1}(V)$ of $V\in \Rep(Q,\alpha+\beta)$ can be identified with the set of
all subrepresentations of $V$. Let $U\subseteq \Rep(Q,\alpha+\beta)$
be the set of all $V\in \Rep(Q,\alpha+\beta)$ such that the fiber $p^{-1}(V)$
is finite. Because $p$ is projective it follows by the semi-continuity of the
dimension of a fiber that $U$ is open. Let us restrict $p$ to
$p^{-1}(U)\to U$.
Now $p:p^{-1}(U)\to U$
is a projective, quasi-finite map, hence it is finite. 
It follows that all fibers have the same cardinality if
counted with multiplicity. It was shown in~\cite{CB2} that a general fiber
of $p$ is reduced (this is not immediately clear in positive characteristic).
Therefore, a general fiber is, set-theoretically, the largest among all fibers
$p^{-1}(V)$, $V\in U$.

\end{proof}

\subsection{The canonical decomposition}
Following Kac, we make the following definition.
\begin{definition}[Section~4 of \cite{K2}]We call
$$
\alpha=\alpha_1\oplus \alpha_2\oplus \cdots\oplus \alpha_s
$$
the canonical decomposition of $\alpha$ if a general representation
of dimension $\alpha$ decomposes into indecomposable representations
of dimensions $\alpha_1,\alpha_2,\dots,\alpha_s$.
\end{definition}
For more details on the canonical decomposition, see \cite{K2,DW2}.

\begin{theorem}[Proposition~3 of \cite{K2}]\label{theoKac}
The expression
$$
\alpha=\alpha_1\oplus \alpha_2\oplus \cdots\oplus \alpha_s
$$
is the canonical decomposition if and only if
$\alpha_1,\dots,\alpha_s$ are Schur roots, and $\ext(\alpha_i,\alpha_j)=0$
for all $i\neq j$.
\end{theorem}
\begin{theorem}[Theorem 3.8 of \cite{S1}]\label{theoremmultiple}
Suppose that $\alpha$ has the canonical decomposition
$$
\alpha_1\oplus\alpha_2\oplus\cdots\oplus\alpha_s
$$
and $p\in \NN$.
Then the canonical decomposition of $p\alpha$ is
$$
[p\alpha_1]\oplus [p\alpha_2]\oplus\cdots\oplus [p\alpha_s],
$$
where 
$$
[p\beta]=\left\{\begin{array}{ll}
\beta^{\oplus p}:=\underbrace{\beta\oplus\cdots\oplus \beta}_p &
\mbox{if $\beta$ is a real or isotropic Schur root;}\\
p\beta & \mbox{if $\beta$ is a non-isotropic imaginary Schur root.}
\end{array}\right.
$$
\end{theorem}
\begin{lemma}[Lemma 5.2 of \cite{S3}] \label{lemrearrange}
Suppose that
$$\alpha=\alpha_1^{\oplus r_1}\oplus \alpha_2^{\oplus r_2}\oplus \cdots \oplus \alpha_s^{\oplus r_s}$$
is the canonical decomposition of $\alpha$,
where $\alpha_1,\alpha_2,\dots,\alpha_s$ are distinct dimension vectors
and $r_1,\dots,r_s$ are positive integers. Then we may assume, after
rearranging $\alpha_1,\dots,\alpha_s$, that $\hom_Q(\alpha_i,\alpha_j)=0$
for all $i<j$.
\end{lemma}
In~\cite{DW2} an efficient algorithm was given to compute the canonical decomposition
of a given dimension vector. A similar recursive procedure was given in in \cite{S3}.
Lemma~\ref{lemrearrange} follows immediately from the correctness of the algorithm in~\cite{DW2},
because the output of the algorithm has the desired property.

For a representation $V\in \Rep(Q,\alpha)$, we have
\begin{equation}\label{eq31}
\langle\alpha,\alpha\rangle=\sum_{x\in Q_0}\alpha(x)^2-
\sum_{a\in Q_1}\alpha(ta)\alpha(ha)=
\dim \GL(Q,\alpha)-\dim \Rep(Q,\alpha).
\end{equation}
On the other hand,
\begin{equation}\label{eq32}
\dim \GL(Q,\alpha)=\dim \GL(Q,\alpha)_V+\dim \GL(Q,\alpha)\cdot V
\end{equation}
where $\GL(Q,\alpha)_V$ is the stabilizer of $V$ and
$\GL(Q,\alpha)\cdot V$ is the orbit of $V$. 
Because $\GL(Q,\alpha)_V$ is equal to the invertible elements
of $\Hom_Q(V,V)$, it follows from (\ref{eq41}) that
\begin{equation}\label{eq33}
\dim \GL(Q,\alpha)_V=\dim \Hom_Q(V,V)=\langle \alpha,\alpha\rangle+\dim
\Ext_Q(V,V).
\end{equation}
Adding (\ref{eq31}), (\ref{eq32}) and (\ref{eq33}) yields
\begin{equation}\label{eq34}
\dim \Rep(Q,\alpha)-\dim \GL(Q,\alpha)\cdot V=\dim \Ext(V,V)
\end{equation}
(see also \cite[Lemma 4]{K2} and \cite{Ri}). In other words,
$\dim \Ext_Q(V,V)$ is the codimension of the orbit of $V$ in $\Rep(Q,\alpha)$.
\begin{definition}
If $\Rep(Q,\alpha)$ contains a dense $\GL(Q,\alpha)$-orbit, then $\alpha$
is called a {\it prehomogeneous\/} dimension vector.
\end{definition}
From (\ref{eq34}) follows that $\alpha$ is prehomogeneous if and only
if there exists a representation $V\in \Rep(Q,\alpha)$ satisfying
$\Ext_Q(V,V)=0$.

Let us prove the following known fact.
\begin{lemma}\label{lemrealSchur}
If $\alpha$ is a real Schur root, then $\alpha$ is prehomogeneous.
\end{lemma}
\begin{proof}
Suppose that $\alpha$ is a real Schur root and $V\in \Rep(Q,\alpha)$
is a Schur representation. From $\langle\alpha,\alpha\rangle=1$
and $\dim \Hom_Q(V,V)=1$ follows that $\dim \Ext_Q(V,V)=0$
by (\ref{eq2}) and (\ref{eq41}).
\end{proof}

\begin{theorem} [See Proposition~4 in \cite{K2}]
Suppose that 
$$
\alpha=\alpha_1\oplus \alpha_2\oplus \cdots\oplus \alpha_s
$$
is the canonical decomposition of $\alpha$. Then $\alpha$
is prehomogeneous if and only of $\alpha_1,\dots,\alpha_s$
are real Schur roots.
\end{theorem}
%

\subsection{The combinatorics of dimension vectors}
We point out that many of the notions we just introduced can be defined
combinatorially. For example the quantity $\ext(\alpha,\beta)$ can be,
in principle, computed using a recursive procedure using
the following result.
\begin{theorem}[Theorem 5.4 of \cite{S1}]\label{Schofield5.4}
We have
$$
\ext(\alpha,\beta)=
\max\{-\langle\alpha',\beta'\rangle\mid \alpha'\hookrightarrow \alpha,\beta\twoheadrightarrow
\beta'\}=
$$
$$=
\max\{-\langle\alpha,\beta'\rangle\mid \beta\twoheadrightarrow
\beta'\}=
\max\{-\langle\alpha',\beta\rangle\mid \alpha'\hookrightarrow \alpha\}.
$$
\end{theorem}
By Theorem~\ref{lemScho}, the conditions $\alpha'\hookrightarrow \alpha$
and $\beta\twoheadrightarrow \beta'$ can be verified by
computing $\ext$-numbers for smaller dimension vectors.

Using (\ref{eq552}) we can also compute $\hom(\alpha,\beta)$ recursively.
\begin{corollary}
The numbers $\hom(\alpha,\beta),\ext(\alpha,\beta)$ do not depend
on the base field $K$.
\end{corollary}
\begin{proposition}(see~\cite{DSW})
The numbers $\alpha\circ\beta$ do not depend on the base field $K$.
\end{proposition}
The numbers $\alpha\circ\beta$ can be computed either in terms
of Schur functors, or equivalently, in terms of Schubert calculus.
This way, $\alpha\circ\beta$ can be expressed as a (perhaps large) sum
of products of Littlewood-Richardson coefficients. See~\cite{DSW}
for more details.

The following definition will be important later.
\begin{definition}
Suppose that $\alpha,\beta\in\NN^{Q_0}$.
We say that two dimension vectors $\alpha,\beta$ are {\it strongly}
perpendicular if
$$
\alpha\circ \beta=1
$$
We will denote this by $\alpha\pperp\beta$.
\end{definition}
We have 
$$
\alpha\pperp\beta\Rightarrow\alpha\perp\beta\Rightarrow\langle\alpha,\beta\rangle=0,
$$
and none of the implications can be reversed.
The following result will be crucial for this paper.
\begin{theorem}[Generalized Fulton Conjecture, Belkale, see the Appendix]\label{thmBelkale}
If $\alpha\circ\beta=1$, then 
$$
p\alpha\circ q\beta=1
$$
 for all $p,q\in \NN$.
\end{theorem}
\begin{remark}
Theorem~\ref{thmBelkale} can be thought of as a generalization of
Fulton's Conjecture (Theorem~\ref{FultonConjecture}). For partitions $\lambda,\mu,\nu$
one can construct a quiver $Q$ and dimension vectors $\alpha,\beta$ such 
that $\alpha\circ\beta=c_{\lambda,\mu}^\nu$,
and $\alpha\circ(n\beta)=c_{n\lambda,n\mu}^{n\nu}$.
We will explain this in more detail in Section~\ref{sec7}.
\end{remark}

A dimension vector $\alpha$ is a Schur root if and only if
there exist no nonzero dimension vectors $\beta,\gamma$
with $\alpha=\beta+\gamma$ and $\ext(\beta,\gamma)=\ext(\gamma,\beta)=0$.
Therefore, the set of Schur roots does not depend on the base field.

\subsection{Perpendicular categories}\label{secperp}
\begin{definition}
A representation $V$ is called {\it exceptional\/} if $\Hom_Q(V,V)\cong K$
and $\Ext_Q(V,V)=0$.
\end{definition}
If $V\in \Rep(Q,\alpha)$ is exceptional, then $V$ is a Schur representation
and $\alpha$ is a Schur root. Moreover, 
$$\langle \alpha,\alpha\rangle=\dim \Hom_Q(V,V)-\dim \Ext_Q(V,V)=1-0=1,
$$
so $\alpha$ is a real Schur root. Conversely, if $\alpha$ is a real Schur root,
then there exists a Schur representation $V\in \Rep(Q,\alpha)$.
From
$$
1=\langle\alpha,\alpha\rangle=\dim\Hom_Q(V,V)-\dim\Ext_Q(V,V)=1-\dim\Ext_Q(V,V)
$$
follows that $\Ext_Q(V,V)=0$. This means that the orbit $\GL(Q,\alpha)\cdot V$
is open and dense in $\Rep(Q,\alpha)$. Therefore, a general representation
of dimension $\alpha$ is isomorphic to $V$. This shows that there
is a natural bijection between real Schur roots and exceptional representations.
\begin{definition}
Suppose that $V$ is an exceptional representation. The {\it right\/}
perpendicular category $V^{\perp}$ is
the full subcategory of $\Rep_K(Q)$ whose objects are all representations
$W$ such that $\Hom_Q(V,W)=\Ext_Q(V,W)=0$.
Similarly, we define the {\it left\/} perpendicular category ${}^\perp V$. 
as the full subcategory of $\Rep_K(Q)$ whose objects are all representations
$W$ for which $\Hom_Q(W,V)=\Ext_Q(W,V)=0$.
\end{definition}
Note that if $V\perp W$ then $\dd V\perp \dd W$ (see Definition~\ref{defperp}).
Conversely, if $\alpha,\beta$ are dimension vectors with $\alpha\perp \beta$
then there exist $V\in \Rep(Q,\alpha)$ and $W\in \Rep(Q,\beta)$ with
$V\perp W$.

The subcategory $V^{\perp}$ (respectively ${}^\perp V$) are closed
under taking kernels, cokernels, direct sums, images and extensions.
\begin{theorem}[Theorem 2.2 of \cite{S2}]
Suppose that $V$ is a sincere representation, i.e., $V(x)\neq 0$ for all $x\in
Q_0$. Then the categories $V^\perp$ and ${}^\perp V$ are equivalent.
\end{theorem}
The equivalence in the Theorem is given by the Auslander-Reiten transform.
If $W$ is an object of ${}^\perp V$, then $W$ does not contain any projective
summands. Then one can define the Auslander-Reiten translate $\tau(W)$ of $W$.
From Auslander-Reiten duality (see properties (5), (6), (7)
on pages 75--76 in~\cite{Ri2}) follows that $\tau(W)$ lies in the
right perpendicular category. The Auslander-Reiten transform induces
an equivalence of categories.

\begin{theorem}[Theorem 2.3 of \cite{S2}]\label{rightperp}
Suppose that $V$ is an exceptional representation of
a quiver $Q$ without oriented cycles and with $n=\#Q_0$ vertices. 
Then $V^{\perp}$
(resp. ${}^\perp V$) is equivalent to $\Rep_K(Q')$ where $Q'$
is a quiver without oriented cycles such that $\# Q_0'=n-1$.
\end{theorem}
Suppose we are in the situation of Theorem~\ref{rightperp}.
The category $V^{\perp}\cong \Rep_{K}(Q')$ has exactly $n-1$
simple objects, say $E_1,E_2,\dots,E_{n-1}$. Now $Q'$ is
the graph with vertices $1,2,\dots,n-1$ and $r_{i,j}:=\dim \Ext_{Q'}(E_i,E_j)$
arrows from $i$ to $j$ for all $i,j$.
We have
$$
\Hom_{Q'}(E_i,E_i)\cong \Hom_{Q}(E_i,E_i)\cong K
$$
for all $i$. This shows that $E_1,\dots,E_{n-1}$ are Schur representations.
The category $V^{\perp}$ is closed under extensions, so every
nontrivial extension of $E_i$ with itself in the category $\Rep(Q)$
would yield a nontrivial extension of $E_i$ with itself in the category $\Rep(Q')$.
Since $\Ext_{Q'}(E_i,E_i)=0$, we have $\Ext_Q(E_i,E_i)=0$.
Therefore, $E_1,\dots,E_{n-1}$ are exceptional representations for $Q$.
Let $W$ be an object of $V^{\perp}\cong \Rep_K(Q')$.
Suppose that, as a representation of $Q'$, its dimension vector
is $\alpha'=(\alpha_1',\dots,\alpha_{n-1}')$. Then $W$ can be build up
from extensions using $\alpha_i'$ copies of $E_i$ for $i=1,2,\dots,n-1$.
This shows that the dimension vector of $W$, as a representation of $Q$
is equal to
$$
\alpha:=\sum_{i=1}^{n-1}\alpha_i'\varepsilon_i
$$
where $\varepsilon_i=\dd_Q E_i$, the dimension vector of $E_i$ seen
as a representation of $Q$. Let us define
$$
I:\NN^{Q_0'}\cong \NN^{n-1}\to \NN^{Q_0}
$$
by
$$
I(\beta_1,\dots,\beta_{n-1})=\sum_{i=1}^{n-1}\beta_i\varepsilon_i.
$$
So if $W$ is a representation of $Q'$ of dimension $\beta$,
then $W$, viewed as a representation of $Q$, has dimension $I(\beta)$.

If $W$ and $Z$ are representations of $Q'$, then 
$$
\Hom_Q(W,Z)\cong\Hom_{Q'}(W,Z)
$$
because $\Rep_K(Q')$ is a full subcategory of $\Rep_K(Q)$. Since
$V^{\perp}$ is closed under extensions, we also have
$$
\Ext_Q(W,Z)\cong \Ext_{Q'}(W,Z).
$$
From this follows that
\begin{equation}\label{homsame}
\hom_{Q'}(\beta,\gamma)=\hom_{Q}(I(\beta),I(\gamma))
\end{equation}
and
\begin{equation}\label{extsame}
\ext_{Q'}(\beta,\gamma)=\ext_{Q}(I(\beta),I(\gamma)).
\end{equation}
Now we also get
\begin{multline}\label{bilsame}
\langle \beta,\gamma\rangle_{Q'}=
\hom_{Q'}(\beta,\gamma)-\ext_{Q'}(\beta,\gamma)
=\\
=\hom_Q(I(\beta),I(\gamma))-\ext_Q(I(\beta),I(\gamma))=
\langle I(\beta),I(\gamma)\rangle_{Q}.
\end{multline}
\begin{lemma}\label{betaSchur}
Suppose that $\beta\in \NN^{Q_0'}$. Then $\beta$ is a Schur
root (for $Q'$) if and only if $I(\beta)$ is a Schur root (for $Q$).
\end{lemma}
\begin{proof}
If $W$ is a Schur representation of dimension $\beta$ for
the quiver $Q'$, then $W$ is also a Schur representation
of dimension $I(\beta)$ as a representation of $Q$. 

Conversely, suppose that $I(\beta)$ is a Schur root. Then a
general representation of dimension $I(\beta)$ is a Schur representation.
Since there exists a representation of dimension $I(\beta)$ in $V^\perp$,
we have that a general representation of dimension $I(\beta)$ lies
in $V^\perp$ (because $V\mapsto \dim\Hom_Q(V,W)$ and $V\mapsto \dim \Ext_Q(V,W)$
are upper semi-continuous). A general representation of dimension $I(\beta)$,
can be seen as a $\beta$-dimensional Schur representation for $Q'$.
\end{proof}
\begin{theorem}\label{Isame}
If $\beta,\gamma\in \NN^{Q_0'}$ and $\beta\perp \gamma$, then
$$
(I(\beta)\circ I(\gamma))_Q=(\beta\circ \gamma)_{Q'}.
$$
\end{theorem}
\begin{proof}
Choose $W\in \Rep(Q',\beta+\gamma)$ in general position.
 So $W$ has $(\beta\circ\gamma)_{Q'}$
subrepresentations of dimension $\beta$. These subrepresentations
correspond to $I(\beta)$-dimensional subrepresentations of $W$,
seen as representations of $Q$.
Suppose that $Z$ is an
$I(\beta)$-dimensional subrepresentation of $W$.
Since $Z$ is a subrepresentation of $W$, 
 $\Hom_Q(V,W)=0$ implies $\Hom_Q(V,Z)=0$. Since
$\langle\beta,\gamma\rangle=0$ we get $\Ext_Q(V,Z)=0$ as well.
This implies that $Z$ lies in $V^{\perp}$, so $Z$ may
be viewed as a representation of $Q'$.
As a representation of $Q$, $W$ has exactly $(\beta\circ\gamma)_{Q'}$
subrepresentations.
By Lemma~\ref{lemfiberlargest} we obtain
$$
(I(\beta)\circ I(\gamma))_Q\geq(\beta\circ \gamma)_{Q'}.
$$

Choose $W\in \Rep(Q,I(\beta)+I(\gamma))$ in general position.
Then $W$ has exactly $(I(\beta)\circ I(\gamma))_Q$  subrepresentations
of dimension $I(\beta)$.
We have $\Hom_Q(V,W)=\Ext_Q(V,W)=0$ (by semicontinuity) so $W$ lies in $V^\perp$.
We may view $W$ as a representation of $Q'$ of dimension $\beta+\gamma$.
Again, the $I(\beta)$-dimensional subrepresentations of $W$,
seen as a representation of $Q$ are exactly the $\beta$-dimensional
subrepresentations of $W$, seen as a representation of $Q'$.
So as a representation of $Q'$, $W$ has exactly
$$
(I(\beta)\circ I(\gamma))_Q
$$
subrepresentations of dimension $\beta$.
Again, Lemma~\ref{lemfiberlargest}  implies
that 
$$
(I(\beta)\circ I(\gamma))_Q\leq (\beta\circ\gamma)_{Q'}.
$$
We conclude that
$$
(I(\beta)\circ I(\gamma))_Q=(\beta\circ \gamma)_{Q'}
$$
\end{proof}
\subsection{Exceptional Sequences}\label{secexseq}
We will introduce exceptional sequences and their basic properties. For more details,
see~\cite{CB,Ri3}.
\begin{definition}
An {\it exceptional sequence} is a sequence $E_1,\dots,E_r$ of exceptional representations
such that $E_i\perp E_j$ for $i<j$. 
\end{definition}
Define $\varepsilon_i:=\dd E_i$ for all $i$. The matrix $(\langle \varepsilon_i,\varepsilon_j\rangle)_{i,j}$
is lower triangular with 1's on the diagonal, and is therefore invertible.
It follows that $\varepsilon_1,\dots,\varepsilon_r$ are linearly independent,
hence $r\leq n:=\#Q_0$.
\begin{definition}
An exceptional sequence $E_1,\dots,E_r$ is called {\it maximal\/} or {\it complete\/} if $r=n$.
\end{definition}
\begin{lemma}[Lemma~1 of \cite{CB}]\label{CBlem1}\label{exext}
If $E_1,E_2,\dots,E_i,E_j,E_{j+1},\dots,E_n$ is an exceptional sequence ($i<j$)
then there exist $E_{i+1},\dots,E_{j-1}$ such that
$E_1,E_2,\dots,E_n$ is an exceptional sequence.
\end{lemma}
In particular, every exceptional sequence $E_1,\dots,E_r$ can be extended
to a maximal exceptional sequence $E_1,\dots,E_n$. To see this, consider
the full subcategory of all representations $V$ such that
$$
E_i \perp V\mbox{ for $i=1,2,\dots,r$}.
$$
Let us denote this category by
$$
E_1^\perp\cap E_2^\perp\cap \cdots \cap E_r^\perp
$$
or simply $E^\perp$ where $E=(E_1,\dots,E_r)$.
Using Theorem~\ref{rightperp} and induction on $r$ 
we see that this category $E^\perp$ is equivalent to the category
of representations $\Rep(Q')$ of a quiver $Q'$ with $n-r$ vertices
and with no oriented cycles. Let $E_{r+1},\dots,E_{n}$ be the 
simple representations (pairwise non-isomorphic) in $E^\perp$
corresponding to the $n-r$ vertices of $Q'$. We can order $E_{r+1},\dots,E_n$
in such way that $E_{j}\perp E_k$ for all $j,k$ with $r+1\leq j<k\leq n$, because $Q'$
has no oriented cycles.
 Then $E_1,\dots,E_n$
is a (maximal) exceptional sequence.
Lemma~\ref{CBlem1} is proved in a similar fashion~(see \cite{CB}).
\begin{definition}
If $E_1,E_2,\dots,E_r$ is an exceptional sequence, then we define ${\mathcal C}(E_1,\dots,E_r)$
as the full subcategory of $\Rep(Q)$ which contains $E_1,\dots,E_r$ and is closed under
extensions, kernels of epimorphisms, and cokernels of monomorphisms.
\end{definition}
\begin{lemma}[Lemma~4 of \cite{CB}]\label{experp}
If $E_1,\dots,E_n$ is a maximal exceptional sequence, then 
${\mathcal C}(E_1,\dots,E_r)$ is equivalent to the category
$${}^\perp E_{r+1}\cap \cdots\cap {}^\perp E_n={}^\perp F
$$
where $F=(E_{r+1},\dots,E_n)$.
\end{lemma}
\begin{lemma}\label{lemsimple}
Suppose that $E_1,E_2,\dots,E_r$ is an exceptional sequence, and
$\Hom_Q(E_i,E_j)=0$ for all $i\neq j$. Then $E_1,\dots,E_r$ are
exactly all simple objects in ${\mathcal C}(E_1,E_2,\dots,E_r)$.
\end{lemma}
\begin{proof}
Let ${\mathcal D}(E_1,E_2,\dots,E_r)$ be the the smallest full subcategory
of $\Rep(Q)$ which contains $E_1,\dots,E_r$ and which is closed
under extensions. The objects of ${\mathcal D}(E_1,\dots,E_r)$ are all
representations which allow a filtration such that each factor
is isomorphic to one of the representations $E_1,\dots,E_r$.
We claim that ${\mathcal D}(E_1,\dots,E_r)={\mathcal C}(E_1,\dots,E_r)$.
To show this, it suffices to show that the category
${\mathcal D}(E_1,\dots,E_r)$ is closed under taking kernels of
epimorphisms and taking cokernels of monomorphisms.
We will show that ${\mathcal D}(E_1,\dots,E_r)$ is closed
under taking cokernels of monomorphisms. Dualizing the arguments
one can then show that ${\mathcal D}(E_1,\dots,E_r)$ is
also closed under taking kernels of epimorphisms.
Suppose that $\phi:V\to W$ is a monomorphism and $V,W$
are objects of ${\mathcal D}(E_1,\dots,E_r)$.
We have filtrations
$$
V=F^0(V)\supset F^{1}(V)\supset \cdots \supset F^{s}(V)=\{0\}
$$
$$
W=F^0(W)\supset F^{1}(W)\supset \cdots \supset F^t(W)=\{0\}
$$
such that all quotients 
$F^i(V)/F^{i+1}(V)$, $F^i(W)/F^{i+1}(W)$
are isomorphic
to one of the representations $E_1,\dots,E_r$.

The case $s=0$ is trivial. Suppose that $s=1$. Then we have that
$V=F^0(V)$
is isomorphic to one of the representations $E_1,\dots,E_r$.
We prove the statement by induction on $t$. If $t=1$, then 
$\phi$ must be an isomorphism because $V$ and $W$ are simple.
So $V/W=0$ and we are done.
Suppose that $t>0$. Let $\psi$ be the composition 
$V\to W\to W/F^{1}(W)$.
Suppose that $\psi=0$. Then $V$ is a subrepresentation of 
$F^1(W)$ 
By induction 
$F^1(W)/V$
is an object of ${\mathcal D}(E_1,\dots,E_r)$,
and 
$W/F^1(W)$
is also an object of ${\mathcal D}(E_1,\dots,E_r)$.
From the exact sequence
$$
0\to F^1(W)/V\to W/V\to W/F^1(W)\to 0
$$
follows that $W/V$ is an object of ${\mathcal D}(E_1,\dots,E_r)$.

Suppose that $\psi\neq 0$. Both 
$V$ and $W/F^1(W)$
are isomorphic to one of the $E_1,\dots,E_r$.
Because $V$ and $W/F^1(W)$ are simple, 
 they are  isomorphic to
each other and  to $E_i$ for some $i$.
Since $E_i$ is exceptional and $\psi$ is nonzero, we must
have that $\psi$ is an isomorphism. It follows that 
$V+F^1(W)=W$
and 
$V\cap F^1(W)=0$.
But then
$$
F^1(W)=F^1(W)/(V\cap F^1(W))\cong (F^1(W)+V)/V=W/V
$$
 This shows that $W/V$ is an object
of ${\mathcal D}(E_1,\dots,E_r)$.
We have proved the case $s=1$.

Suppose now that $s>1$. We will prove the theorem
by induction on $s$. By the above, we now that 
$W/F^{s-1}(V)$ and $V/F^{s-1}(V)$
are objects of ${\mathcal D}(E_1,\dots,E_r)$.
From induction and the exact sequence
$$
0\to V/F^{s-1}(V)\to W/F^{s-1}(V)\to W/V\to 0
$$
we conclude that $W/V$ also is an object of ${\mathcal D}(E_1,\dots,E_r)$.
We have proved that ${\mathcal D}(E_1,\dots,E_r)={\mathcal C}(E_1,\dots,E_r)$.

Since every object in ${\mathcal D}(E_1,\dots,E_r)={\mathcal C}(E_1,\dots,E_r)$ has a 
subobject isomorphic
to one of the representations $E_1,\dots,E_r$, the only possible
simple objects of ${\mathcal C}(E_1,\dots,E_r)$ are $E_1,\dots,E_r$.
 It is also easy to see that each $E_i$ is simple. Indeed,
 if $W$ is a proper subrepresentation of $E_i$, then $W$ has a proper subrepresentation isomorphic
 to $E_k$ for some $k\neq i$. But then $E_k$ is a proper subrepresentation of $E_i$
 which contradicts the assumption that $\Hom_Q(E_k,E_i)=0$.

\end{proof}
As we have noted before, there is a bijection between real Schur roots and exceptional representations.
Suppose that $E_1,\dots,E_r$ are exceptional representations and let $\varepsilon_i:=\dd E_i$.
We  have seen that the orbit of $E_i$ in $\Rep(Q,\varepsilon_i)$ is dense.
From this follows that
$$
\dim \Hom_Q(E_i,E_j)=\hom_Q(\varepsilon_i,\varepsilon_j)
$$
and
$$
\dim \Ext_Q(E_i,E_j)=\ext_Q(\varepsilon_i,\varepsilon_j)
$$
for all $i,j$. This allows us to give a more combinatorial definition
of an exceptional sequence.
\begin{definition}
A sequence of dimension vectors $\varepsilon_1,\dots,\varepsilon_r$ is
called an exceptional sequence if $\varepsilon_1,\dots,\varepsilon_r$
are real Schur roots, and $\varepsilon_i\perp \varepsilon_j$ for all $i<j$.
\end{definition}
So if $E_1,\dots,E_r$ is an exceptional sequence of quiver representations,
then $\varepsilon_1,\dots,\varepsilon_r$ is an exceptional sequence of dimension
vectors, where $\varepsilon_i=\dd E_i$ for all $i$. Conversely,
suppose that $\varepsilon_1,\dots,\varepsilon_r$ is an exceptional sequence
of dimension vector. Since $\varepsilon_i$ is a real Schur root,
the exists a unique dense orbit in $\Rep(Q,\varepsilon_i)$. Let $E_i$ be
a representation that lies in that orbit ($E_i$ is unique op to isomorphism).
Then $E_1,\dots,E_r$ is an exceptional sequence of representations.
\begin{theorem}[Theorem 4.1 of  \cite{S1}]\label{theohomorext}
Suppose that $\alpha,\beta$ are Schur roots such that $\ext(\alpha,\beta)=0$. 
Then $\hom(\beta,\alpha)=0$
or $\ext(\beta,\alpha)=0$. Moreover, if both $\alpha$ and $\beta$ are imaginary,
then $\hom(\beta,\alpha)=0$.
\end{theorem}

\begin{theorem}[Embedding Theorem]\label{theoembedding}
Suppose that $\varepsilon_1,\dots,\varepsilon_r$ is an exceptional
sequence for the quiver $Q$ (without orientations). 
Suppose that $\langle \varepsilon_i,\varepsilon_j\rangle\leq 0$ for all $i>j$.
We define a quiver $Q'$ with vertices $Q'_0=\{1,2,\dots,r\}$.
We draw $-\langle \varepsilon_i,\varepsilon_j\rangle$ arrows
from $i$ to $j$ for all $i>j$.
We define
$$
I:\NN^{Q'_0}\cong \NN^r\to \NN^{Q_0}
$$
by
$$
I(\beta_1,\dots,\beta_r)=\sum_{i=1}^r \beta_i\varepsilon_i
$$
for all $\beta=(\beta_1,\dots,\beta_r)\in \NN^{Q'_0}\cong \NN^r$.
\begin{enumerate}
\renewcommand{\theenumi}{\alph{enumi}}
\item For all $\beta,\gamma\in \NN^{Q'_0}$ we have
$$
\hom_Q(I(\beta),I(\gamma))=\hom_{Q'}(\beta,\gamma),
$$
$$
\ext_Q(I(\beta),I(\gamma))=\ext_{Q'}(\beta,\gamma),
$$
$$
\langle I(\beta),I(\gamma)\rangle_Q=\langle\beta,\gamma\rangle_{Q'}.
$$
\item The dimension vector $\beta$ is a Schur root (for $Q'$) if and only if $I(\beta)$ is a Schur root (for $Q$).
\item For all $\beta,\gamma\in \NN^{Q'_0}$ with $\beta\perp \gamma$ we have
$$
(I(\beta)\circ I(\gamma))_Q=(\beta\circ\gamma)_{Q'}
$$
\end{enumerate}
\end{theorem}
\begin{proof}
Let $E_i$ be a representation
corresponding to the dense orbit in $\Rep(Q,\varepsilon_i)$.
Then $E_1,E_2,\dots,E_r$ is an exceptional sequence. 
For $i>j$ we have 
$$\langle \varepsilon_i,\varepsilon_j\rangle=\hom(\varepsilon_i,\varepsilon_j)-\ext(\varepsilon_i,\varepsilon_j)\leq 0
$$
Since $\hom(\varepsilon_i,\varepsilon_j)=0$ or $\ext(\varepsilon_i,\varepsilon_j)=0$ by 
Theorem~\ref{theohomorext},
we must have $\hom(\varepsilon_i,\varepsilon_j)=0$. It follows that $\Hom_Q(E_i,E_j)=0$
for $i>j$. From Lemma~\ref{lemsimple} we see that $E_1,E_2,\dots,E_r$ are exactly all
simple objects in ${\mathcal C}(E_1,\dots,E_r)$. We can extend $E_1,\dots,E_r$ to
a maximal exceptional sequence $E_1,E_2,\dots,E_n$ using  Lemma~\ref{exext}.
Now, we have equality
$${\mathcal C}(E_1,\dots,E_r)={}^\perp E_{r+1}\cap \cdots \cap {}^\perp E_n
$$
by Lemma~\ref{experp}.

Suppose that $r=n-1$. Then we are exactly in the situation of Theorem~\ref{rightperp}.
We note that the definition
of the quiver $Q'$ and the map $I$ coincide with the definitions in the discussion
after Theorem~\ref{rightperp}. In particular, we get (see~(\ref{homsame}), (\ref{extsame}),
 (\ref{bilsame}))
$$
\hom_{Q'}(\beta,\gamma)=\hom_Q(I(\beta),I(\gamma))
$$
$$
\ext_{Q'}(\beta,\gamma)=\ext_Q(I(\beta),I(\gamma)).
$$
and
$$
\langle \beta,\gamma\rangle_{Q'}=\langle I(\beta),I(\gamma)\rangle_{Q}.
$$
Lemma~\ref{betaSchur} implies that $\beta\in \NN^{Q'_0}$ is a Schur root if and only
if $I(\beta)$ is a Schur root.
Finally, Theorem~\ref{Isame} implies that
$$
(\beta\circ\gamma)_{Q'}=(I(\beta)\circ I(\gamma))_Q.
$$

If $r<n$, the theorem can be proved in a similar way, using induction on $n-r$.
\end{proof}
Theorem~\ref{theoembedding} allows us to ``embed'' the combinators of
a quiver $Q'$ into the combinators of another quiver $Q$, using exceptional sequences.
This is a very useful tool. Of course, we would like to be able to 
construct such exceptional sequences. Crawley-Boevey defined an action
of the braid group on the set of all maximal exceptional sequences
and he proved that this action is transitive (see~\cite{CB}).

So let us assume that $Q$ is a quiver without orientations and $n$ vertices.
Assume that the vertices are labeled by $1,2,\dots,n$. Since there are no oriented
cycles, we can arrange the vertices in such way that for every arrow $i\to j$ in the
quiver we have $i>j$. If $S_i$ is the simple representation corresponding
to the vertex $i$, then $S_1,S_2,\dots,S_n$ is a maximal exceptional sequence.

The braid group ${\mathcal B}_n$ is the group given by generators
$s_1,s_2,\dots,s_{n-1}$ and relations
$$
s_is_j=s_js_i,\quad\mbox{if $|i-j|\geq 2$ and}
$$
$$
s_is_{i+1}s_i=s_{i+1}s_is_{i+1}
$$
for $i=1,2,\dots,n-2$. Let ${\mathcal E}_Q$ be the set of all maximal exceptional
sequences for $Q$. We will define the action of ${\mathcal B}_n$ on ${\mathcal E}_Q$ as follows.
Suppose that $E_1,E_2,\dots,E_n$ is a maximal exceptional sequence.
Let $\varepsilon_1,\varepsilon_2,\dots,\varepsilon_n$ be their dimension vectors.
Then we define
$$
s_i(E_1,E_2,\dots,E_n)=(E_1,\dots,E_{i-1},E_{i+1}',E_i,E_{i+2},\dots,E_n)
$$
where $E_{i+1}'$ is the unique simple object in ${}^\perp E_i$ within
the subcategory ${\mathcal C}(E_i,E_{i+1})$. Let $\varepsilon_{i+1}'$ be the dimension
vector of $E_{i+1}'$. 
From Lemma~\ref{experp} follows that ${\mathcal C}(E_i,E_{i+1})={\mathcal C}(E_{i+1}',E_i)$.
Therefore, we have 
$$
\ZZ\varepsilon_i+\ZZ\varepsilon_{i+1}=\ZZ\varepsilon_i+\ZZ\varepsilon_{i+1}'.
$$
It follows that $\varepsilon_{i+1}'=\varepsilon_{i+1}+k\varepsilon_{i}$
or $\varepsilon_{i+1}'=-\varepsilon_{i+1}+k\varepsilon_i$ for some $k\in \ZZ$.
We know that $\langle \varepsilon_{i+1}',\varepsilon_i\rangle=0$. There
can only be two possibilities, namely 
$$
\varepsilon_{i+1}'\in \{
\varepsilon_{i+1}-\langle \varepsilon_{i+1},\varepsilon_i\rangle \varepsilon_{i},
-(\varepsilon_{i+1}-\langle \varepsilon_{i+1},\varepsilon_i\rangle \varepsilon_{i})\}.
$$
This uniquely determines $\varepsilon_{i+1}'$ because only one
of two vectors on the right can lie in $\NN^{Q_0}$.

We also define
$$
s_{i}^{-1}(E_1,E_2,\dots,E_n)=(E_1,\dots,E_{i-1},E_{i+1},E_i',E_{i+2},\dots,E_n)
$$
where
$E_i'$ is the unique simple object in $E_{i+1}^{\perp}$ within the subcategory
${\mathcal C}(E_i,E_{i+1})$. Let $\varepsilon_i'=\dd E_i'$.
One can prove that
$$
\varepsilon_i'\in \{\varepsilon_i-\langle \varepsilon_i,\varepsilon_{i+1}\rangle\varepsilon_{i+1},
-(\varepsilon_i-\langle \varepsilon_i,\varepsilon_{i+1}\rangle\varepsilon_{i+1})\}.
$$
One can show that this gives a well-defined action of the braid group ${\mathcal B}_n$
on the set ${\mathcal E}_Q$ of maximal exceptional sequences.
\begin{theorem}[\cite{CB}]
The action of ${\mathcal B}_n$ on ${\mathcal E}_Q$ is transitive.
\end{theorem}
\subsection{Stability and GIT-quotients}\label{secKing}
King studied in \cite{K} moduli spaces for representations of quivers (and more generally, finite dimensional
algebras). Here we will study the moduli space of a quiver $Q$ without oriented cycles.
Suppose that $\sigma\in \Gamma^\star$ is a weight. A representation $V\in \Rep(Q,\alpha)$
is called $\sigma$-semi-stable if there exists a positive integer $m$ and a semi-invariant
 $f\in \SI(Q,\alpha)_{m\sigma}$ such that $f(V)\neq 0$.
King gave
the following nice characterization of $\sigma$-stable and $\sigma$-semi-stable
representations.
\begin{theorem}[Proposition 3.1 of~\cite{K}]\ \label{KingCrit}
\begin{enumerate}
\item
A representation $V$ is $\sigma$-semi-stable if $\sigma(\dd V)=0$
and $\sigma(\dd W)\leq 0$ for all subrepresentations $W\subseteq V$.
\item
 representation $V$ is $\sigma$-stable if $\sigma(\dd V)=0$
and $\sigma(\dd W)< 0$ for all nonzero, proper subrepresentations $W\subset V$.
\end{enumerate}
\end{theorem}
In King's setup the
inequalities go the other way. His 
 convention for writing weights are opposite to ours.
\begin{example}
Suppose that $\sigma=0$.
Then a representation is $\sigma$-stable if and only if
it is simple. 
Every representation is $\sigma$-semi-stable.
\end{example}
In general, $\sigma$-stable representations are indecomposable. Indeed,
if $V=V_1\oplus V_2$ for some nontrivial representations $V_1$ and $V_2$,
then $0=\sigma(\dd V)=\sigma(\dd V_1)+\sigma(\dd V_2)$,
so $\sigma(\dd V_1)\leq 0$ or $\sigma(\dd V_2)\leq 0$.
If $d$ is a positive integer, then a representation $V$ is $\sigma$-stable
if and only if it is $(d\sigma)$-stable. 
The set of $\sigma$-semi-stable representations in $\Rep(Q,\alpha)$ is denoted by
$\Rep(Q,\alpha)_{\sigma}^{\rm ss}$. We will assume that $\Rep(Q,\alpha)_{\sigma}^{\rm ss}$
is nonempty. In other words $\SI(Q,\alpha)_{n\sigma}\neq 0$ for some positive integer $n$.

Recall that $\sigma$ defines a character $\GL(Q,\alpha)\to K^\star$. Let $\GL(Q,\alpha)_\sigma$ be the
kernel of this character.
\begin{lemma}
Suppose that $\sigma$ is not divisible by the characteristic of $K$.
The invariant ring of $K[\Rep(Q,\alpha)]$ with respect to $\GL(Q,\alpha)_\sigma$ is equal to 
$$
K[\Rep(Q,\alpha)]^{\GL(Q,\alpha)_\sigma}=\bigoplus_{n\geq 0}\SI(Q,\alpha)_{n\sigma}.
$$
\end{lemma}
\begin{proof}
The inequality $\supseteq$ is easy to see. We will prove $\subseteq$. Suppose that $f\in K[\Rep(Q,\alpha)]^{\GL(Q,\alpha)_\sigma}$. Then in particular 
$$
f\in \SI(Q,\alpha)=\bigoplus_{\tau\in \Gamma^\star} \SI(Q,\alpha)_\tau.
$$
Write $f=\sum_{\tau}f_\tau$ with $f_\tau\in \SI(Q,\alpha)_\tau$ for all $\tau\in \Gamma^\star$.
Suppose that $f_\eta\neq 0$. If $\eta$ is not an integer multiple of $\sigma$,
then there exists an $A\in \GL(Q,\alpha)_\sigma$ such that $\eta(A)\neq 1$.
But then 
$$
\sum_\tau f_\tau=f=A\cdot f=\sum_\tau  A\cdot f_\tau=\sum_\tau \tau(A) f_\tau
$$
so $f_\eta=\eta(A)f_\eta$ and $f_\eta=0$. So $f_\eta=0$ unless $\eta$ is an integer
multiple of $\sigma$. This shows that $f\in \bigoplus_{n\in \ZZ}\SI(Q,\alpha)_{n\sigma}$.
We have assumed that 
there exists a nonzero $t\in \SI(Q,\alpha)_{\sigma}\neq 0$ for some positive $n$.
This implies that $\SI(Q,\alpha)_{-m\sigma}=0$ for all $m>0$, for
if $u\in \SI(Q,\alpha)_{-m\sigma}$ then $u^nt^m\in K[\Rep(Q,\alpha)]^{\GL(Q,\alpha)}=K$
and $u$ must be constant. But the only constant function in
 $ \SI(Q,\alpha)_{-m\sigma}$ is 0.
 \end{proof}
 Let  $X=\Spec(R)$ be the affine variety corresponding to the ring 
$$
R= K[\Rep(Q,\alpha)]^{\GL(Q,\alpha)_\sigma}=\bigoplus_{n\in \NN}\SI(Q,\alpha)_{n\sigma}.
$$
(Here $\Spec$ denotes the spectrum of {\it maximal\/} ideals.)
 The inclusion $R\subseteq K[\Rep(Q,\alpha)]$ corresponds to a a  morphism
 $$
 \psi:\Rep(Q,\alpha)\twoheadrightarrow X
 $$
 This morphism is an affine quotient with respect to a reductive group.
 Such a quotient map is known to be surjective. Note that $\psi^{-1}(0)$
 is exactly the complement of $\Rep(Q,\alpha)_\sigma^{\rm ss}$
 in $\Rep(Q,\alpha)$. Here $0\in X$ is the point corresponding
 to the homogeneous maximal ideal ${\mathfrak m}=\bigoplus_{n>0}\SI(Q,\alpha)_{n\sigma}$.
Define
$$
Y=\Proj\Big(\bigoplus_{n\geq 0}\SI(Q,\sigma)_{n\sigma}\Big).
$$
There is a natural surjective map $X\setminus\{0\}\to Y$. So we have a diagram
$$
\xymatrix{
\Rep(Q,\alpha)_\sigma^{\rm ss}\ar[r]\ar[rd] & X\setminus\{0\}\ar[d]\\
 & Y}
$$

The map $\pi:\Rep(Q,\alpha)_{\sigma}^{\rm ss}\to Y$ is a GIT-quotient
of $\Rep(Q,\alpha)$. Note that this quotient depends on $\sigma$.
Again, the quotient $\pi$ is surjective. 

\section{Stability}\label{sec3}
\subsection{Harder-Narasimhan and Jordan-H\"older filtrations}
The notion of stability was compared (see~\cite[Section 3]{R}) 
to the stability defined through a slope of two
additive functions. In the quiver setting, the set
of additive functions on the category can be identified with $\Gamma^\star$.
Suppose that $\sigma,\tau\in \Gamma^\star$ and $\tau(\alpha)>0$ for all
nonzero dimension vectors $\alpha$.
The {\it slope\/} of a representation $V$ is defined by
$$\mu(V):=\frac{\sigma(\dd V)}{\tau(\dd V)}.$$
\begin{definition}\label{defi1}
A representation $V$ is called $(\sigma:\tau)$-semi-stable if
 $$
 \mu(W)\leq \mu(V)
 $$
 for every proper subrepresentation $W\subseteq V$.
 It is called $(\sigma:\tau)$-stable
 if
 $$
 \mu(W)<\mu(V)
 $$
 for every proper subrepresentation $W\subset V$.
\end{definition}
Rudakov showed (see~\cite[Lemma 3.2]{R}) that the order $\prec$ defined by
$$
V\prec W\Leftrightarrow \mu(V)<\mu(W)
$$
and equivalence relation defined by
$$
V\asymp W\Leftrightarrow \mu(V)=\mu(W)
$$
defines a stability order in the sense of \cite[Section~1]{R}. 
\begin{remark}\label{sigmasigmatau}
The notions of $\sigma$-stability and $(\sigma:\tau)$-stability
are closely related. 

Suppose that $V$ is a representation of $Q$ and choose $\tau$ such that
$\tau(\alpha)>0$ for every dimension vector $\alpha$. For example,
define $\tau(\alpha):=\sum_{x\in Q_0}\alpha(x)$.
If $V$ is $\sigma$-(semi-)stable, then $V$ is also
$(\sigma:\tau)$-(semi-)stable.

Also, if $V$ is $(\sigma:\tau)$-(semi-)stable, and
$\sigma(\dd V)/\tau(\dd V)=a/b$ where $a,b\in \ZZ$ and $b$ is positive, then
$V$ is $(b\sigma-a\tau)$-(semi-)stable.
\end{remark}
 Theorems 2 and 3 from \cite{R} give immediately the following results.
\begin{proposition}
With the assumptions of Definition~\ref{defi1} we have:
\begin{enumerate}
\renewcommand{\theenumi}{\arabic{enumi}}
\item (Harder-Narasimhan filtration) Every object $V$ has a filtration
$$
V= F^0_H(V)\supset F^1_H(V)\supset\ldots\supset F^m_H(V)\supset F^{m+1}_H(V)=0
$$
such that
\begin{enumerate}
 \item each factor $G^i_H(V)=F^i_H(V)/ F^{i+1}_H(V)$ is $(\sigma:\tau)$-semi-stable; 
\item 
$$
G^0_H(V)\prec G^1_H(V)\prec\ldots\prec G^m_H(V).
$$
\end{enumerate}
The filtration with properties (a), (b) is unique.
\item  (Jordan-H\"older filtration) Every $(\sigma:\tau)$-semi-stable object $V$ has a filtration
$$
V= F^0_J(V)\supset F^1_J(V)\supset\ldots\supset F^m_J(V)\supset F^{m+1}_J(V)=0
$$
such that 
\begin{enumerate}
\item each factor $G^i_J(V)=F^i_J(V)/ F^{i+1}_J(V)$ is $(\sigma:\tau)$-stable; 
\item
$$
G^0_J(V)\asymp G^1_J(V)\asymp\ldots\asymp G^m_J(V).
$$
\end{enumerate}
The set of factors $\lbrace G^i_J(V)\rbrace$ is uniquely determined by the properties (a), (b).
\end{enumerate}
\end{proposition}
In part (2), the Jordan-H\"older filtration itself may not be unique.
\begin{proposition}\label{propFm}
Suppose that a  representation $V$ has the Harder-Narasimhan filtration
$$
V=F^0(V)\supset F^1(V)\supset \dots\supset F^m(V)\supset F^{m+1}(V)=0.
$$
If $W$ is a nonzero subrepresentation of $V$ then we have
$\mu(W)\leq \mu(F^m(V))$.
If moreover $\mu(W)=\mu(F^m(V))$, then $W\subseteq F^m(V)$.
\end{proposition}
\begin{proof}
This follows immediately from Proposition 1.9 and Proposition 1.13 in \cite{R}.
\end{proof}
\begin{remark}
Of course, by taking the Harder-Narasimhan filtration and using a Jordan-H\"older
filtration on each quotient, one obtains a filtration of an object $V$
such that every quotient is $(\sigma:\tau)$-{\it stable}. Such a filtration 
we will call
a Harder-Narasimhan-Jordan-H\"older filtration or HNJH-filtration for short.
\end{remark}
\begin{lemma}\label{FilterClosed}
Let $Q$ be a quiver and
$\alpha_0=\alpha,\alpha_1,\alpha_2,\dots,\alpha_m,\alpha_{m+1}=0$
are dimension vectors.
The set of all $V\in \Rep(Q,\alpha)$ of representations which allow
a filtration 
$$
V=F^0(V)\supset F^1(V)\supset\cdots\supset F^{m+1}(V)=0
$$
with $\dd F^i(V)=\alpha_i$ is Zariski closed.
\end{lemma}
\begin{proof}
The proof goes exactly as in \cite[Section 3]{S1}.
\end{proof}
The following Lemma follows immediately from the definition of
$\sigma$-semi-stability.
\begin{lemma}\label{semisub}
Suppose that $V$ is a $\sigma$-semi-stable representation, and $W$
is a subrepresentation with $\sigma(\dd W)=0$, then 
$W$ is $\sigma$-semi-stable as well.
\end{lemma} 
\begin{lemma}\label{lemlocallyclosed}
Let $Q$ be a quiver and
let $\alpha_0=\alpha,\alpha_1,\alpha_2,\dots,\alpha_m,\alpha_{m+1}=0$
be dimension vectors.
\begin{enumerate}
\renewcommand{\theenumi}{\alph{enumi}}
\item 
Define $U=U_{HN}(\alpha_0,\alpha_1,\dots,\alpha_{m+1})$
as the set of all representations
$V\in \Rep(Q,\alpha)$  which 
have a Harder-Narasimhan filtration
$$
V=F^0(V)\supset F^1(V)\supset\cdots\supset F^{m+1}(V)=0
$$
with $\dd F^i(V)=\alpha_i$ for all $i$. Then $U$ is locally closed.
\item
Define $U=U_{HNJH}(\alpha_0,\alpha_1,\dots,\alpha_{m+1})$
as the set of all representations $V\in \Rep(Q,\alpha)$  which 
have a HNJH-filtration
$$
V=F^0(V)\supset F^1(V)\supset\cdots\supset F^{m+1}(V)=0
$$
with $\dd F^i(V)=\alpha_i$ for all $i$.
Then $U$ is locally closed.
\end{enumerate}
\end{lemma}
\begin{proof}
(a) Suppose that $U$ is not locally closed.
This means that $\overline{U}\setminus U$ is not Zariski closed.
Now $\overline{U}\setminus U$ is contained
in the union of all sets of the form $U_{HN}(\beta_0,\dots,\beta_{r+1})$.
Since there are only finitely many such sets,  we can find
$\beta_0,\beta_1,\dots,\beta_{r+1}$ 
such that  the Zariski closure of $U'\cap (\overline{U}\setminus U)$ is
not contained in $\overline{U}\setminus U$ for 
$U'=U_{HN}(\beta_0,\dots,\beta_{r+1})$.
So we have $\overline{U'}\cap U\neq \emptyset$
and $U'\cap \overline{U}\neq \emptyset$.

Let $V\in U\cap \overline{U'}$. Since $V\in U$, it has a Harder-Narasimhan filtration
$$
V= F^0_H(V)\supset F^1_H(V)\supset\ldots\supset F^m_H(V)\supset F^{m+1}_H(V)=0
$$
with $\dd F^i_H(V)=\alpha_i$ for all $i$.
Because $V\in \overline{U'}$,  Lemma~\ref{FilterClosed} implies that it also has a filtration
$$
V= \widehat{F}^0(V)\supset \widehat{F}^1(V)\supset\ldots\supset 
\widehat{F}^r(V)\supset \widehat{F}^{r+1}(V)=0
$$
with $\dd \widehat{F}^i(V)=\beta_i$ for all $i$.

Let $V'\in \overline{U}\cap U'$. Because $V\in U'$, it has a Harder-Narasimhan filtration
$$
V'= \widehat{F}^0_H(V')\supset \widehat{F}^1_H(V')\supset\ldots\supset 
\widehat{F}^r_H(V')\supset \widehat{F}^{r+1}_H(V')=0
$$
with $\dd \widehat{F}^i_H(V')=\beta_i$.
Because $V'\in \overline{U}$,  Lemma~\ref{FilterClosed} implies that it has
a filtration
$$
V'= F^0(V')\supset F^1(V')\supset\ldots\supset F^m(V')\supset F^{m+1}(V')=0
$$
with $\dd F^i(V')=\alpha_i$ for all $i$.

By induction on $i$ we will prove that
$\widehat{F}^{r-i}(V)=F_{H}^{m-i}(V)$ and
$F^{m-i}(V')=\widehat{F}_H^{r-i}(V')$.

We start the induction with the case $i=0$. 
Since $\widehat{F}^r(V)$ is a proper subrepresentation of $V$,
Proposition~\ref{propFm} implies that
$\mu(\beta_r)\leq\mu(\alpha_m)$.
Similarly, because $F^{m}(V')$ is a subrepresentation of $V'$,
$$
\mu(\alpha_m)\leq
\mu(\beta_r)
$$
by Proposition~\ref{propFm}.
So we get $\mu(\alpha_m)=\mu(\beta_r)$.
Proposition~\ref{propFm}
shows that $\widehat{F}^{r}(V)\subseteq F_H^m(V)$
and  $\widehat{F}_H^r(V')\supseteq F^m(V')$.
We obtain $\beta_r=\alpha_m$, $\widehat{F}^r(V)=F^m_N(V)$ and
$\widehat{F}^r_H(V')=F^m(V')$.

Suppose we already have shown that 
$\widehat{F}^{r-j}(V)=F_H^{m-j}(V)$
and  $\widehat{F}_H^{r-j}(V')=F^{m-j}(V')$
for all $j<i$.
Now $\widehat{F}^{r-i}(V)/\widehat{F}^{r-i+1}(V)$
is a subrepresentation of $V/F_H^{m-i+1}(V)$
and $F^{m-i}(V')/F^{m-i+1}(V')$ is a subrepresentation
of $V'/\widehat{F}_H^{r-i+1}(V')$. Similar reasoning
as before shows that
$$
\mu(\beta_{r-i}-\beta_{r-i+1})=
\mu(\alpha_{m-i}-\alpha_{m-i+1}).
$$
It follows that
$$\widehat{F}^{r-i}(V)/\widehat{F}^{r-1+1}(V)\subseteq
F^{m-i}_H(V)/F^{m-i+1}_H(V)$$
and
$$F^{m-i}(V')/F^{m-i+1}(V')\subseteq
\widehat{F}_H^{r-i}(V')/\widehat{F}_H^{r-i+1}(V').$$
We get $\beta_{r-i}-\beta_{r-i+1}=\alpha_{m-i}-\alpha_{m-i+1}$
and equality in the previous inclusions.
We conclude that $\beta_{r-i}=\alpha_{m-i}$,
$F^{m-i}_H(V)=\widehat{F}^{r-i}(V)$ and
$F^{m-i}(V')=\widehat{F}_H^{r-i}(V')$.

We have proved by induction that $\beta_{r-i}=\alpha_{m-i}$ for all $i$.
In particular we get $r=m$. We have $U=U'$. But
then the Zariski closure of
$$U'\cap (\overline{U}\setminus U)=U\cap (\overline{U}\setminus
U)=\emptyset$$ is contained in $\overline{U}\setminus U$, contrary to 
our assumptions.
 
(b) Suppose that $U$ is not locally closed. Then $\overline{U}\setminus U$ is not Zariski
closed.
Now $\overline{U}\setminus U$ is contained
in the union of all $U_{HNJH}(\beta_0,\dots,\beta_{r+1})$. 
Since this is only a finite union, there exist
$\beta_0,\beta_1,\dots,\beta_{r+1}$ 
such that the Zariski closure of $U'\cap(\overline U\setminus U)$
is not contained in $\overline{U}\setminus U$, where
 $U'=U_{HNJH}(\beta_0,\dots,\beta_{r+1})$. So we have
 $U'\cap \overline{U}\neq \emptyset$ and $U\cap \overline{U'}\neq \emptyset$.
From (a) follows that 
elements of $U$ and of $U'$ have
the same dimensions in the Harder-Narasimhan filtrations,
say
$$
U\subseteq U_{HN}(\gamma_0,\gamma_1,\dots,\gamma_{t+1})\supseteq
U'.
$$
for some dimension vectors $\gamma_0=\alpha,\gamma_1,\dots,\gamma_{t+1}=0$.
It follows that
$$
\{\alpha_0,\dots,\alpha_{m+1}\}\supseteq \{\gamma_0,\dots,\gamma_{t+1}\}
$$
and
$$
\{\beta_0,\dots,\beta_{r+1}\}\supseteq \{\gamma_0,\dots,\gamma_{t+1}\}.
$$
Let $V'\in U'\cap \overline{U}$. Now $V'$
has also a filtration $\{F^{i}(V')\}$ with dimensions
 $\alpha_0,\alpha_1,\dots,\alpha_{m+1}$ by Lemma~\ref{FilterClosed}.
This filtration is a refinement of the Harder-Narasimhan filtration of $V'$.
The quotients of this filtrations are $(\sigma:\tau)$-semi-stable
(see Remark~\ref{sigmasigmatau} and Lemma~\ref{semisub}).
This means that the filtration $\{F^i(V')\}$
can be refined to a HNJH-filtration.
We deduce that $m\leq r$.
Considering a representation $V\in U\cap \overline{U'}$ implies
the reverse inequality. Therefore $m=r$. 
If $V'\in U'\cap \overline{U}$, then the filtration $\{F^i(V)\}$ with
dimensions $\alpha_0,\alpha_1,\dots,\alpha_{m+1}$
Any HNJH-filtration of $V'$ has $r=m$ stable factors.
This implies that the filtration $\{F^i(V)\}$ cannot be further refined
to a HNJH-filtration. This means that it already is a HNJH-filtration.
This shows that $U'\cap\overline{U}\subseteq U$. 
This contradicts the fact that $U'\cap (\overline{U}\setminus U)\neq \emptyset$.
\end{proof}

\begin{proposition}\label{prop2}
Let $Q$ be a quiver without oriented cycles,
$\sigma,\tau$ ($\tau(\alpha)>0$ for
all nonzero dimension vectors) be weights and $\alpha$ be a dimension
vector. 
\begin{enumerate}
\renewcommand{\theenumi}{\alph{enumi}}
\item There exists a nonempty open set $U\subseteq \Rep_K (Q,\alpha )$ such that for 
$V\in U$ the dimensions of the factors of the Harder-Narasimhan 
filtration with respect to $(\sigma:\tau)$ of $V$ are constant.
\item There exists a nonempty open set $U\subseteq \Rep_K (Q,\alpha )$ such that for 
$V\in U$ the dimensions of the factors of the HNJH-filtration
with respect to $(\sigma:\tau)$ of $V$ are constant.
\item If a general representation of dimension $\alpha$ is $\sigma$-semi-stable, then
there exists an nonempty open set $U\subseteq \Rep_K (Q,\alpha )$ such that for 
$V\in U$ the dimensions of the factors of a Jordan-H\"older-filtration
of $V$ are constant.
\end{enumerate}
\end{proposition}
\begin{proof}\ 

(a) There exist dimension vectors
$\alpha_0=\alpha,\alpha_1,\dots,\alpha_m,\alpha_{m+1}=0$ such that
$$
U=U_{HN}(\alpha_0,\dots,\alpha_{m+1})
$$
lies dense in $\Rep(Q,\beta)$. Since $U$ is locally closed by
Lemma~\ref{lemlocallyclosed},
we conclude that $U$ is open and dense in $\Rep(Q,\beta)$.

(b) goes similarly.

(c) Follows directly from (b) and Remark~\ref{sigmasigmatau}.
\end{proof}
\begin{lemma}\label{lemmorphism}\ 

\begin{enumerate}
\renewcommand{\theenumi}{\alph{enumi}}
\item Suppose $\sigma\in \Gamma^\star$ is a weight and $V,W$  are $\sigma$-stable
representations of the quiver. 
If $\phi:V\to W$ is a morphism then either $\phi=0$ or $\phi$ is
an isomorphism.
\item Suppose that $\sigma,\tau\in \Gamma^\star$ and $\tau(\alpha)>0$ for all 
nonzero dimension vectors $\alpha$. Assume that $V,W$ are $(\sigma:\tau)$-stable
representations with $\mu(V)\geq \mu(W)$, where $\mu=\sigma/\tau$. If $\phi:V\to W$ is a morphism
then either $\phi=0$ or $\phi$ is an isomorphism.
\end{enumerate}
\end{lemma}
\begin{proof}
This follows immediately from \cite[Theorem 1(d)]{R}.
\end{proof}

\subsection{The $\sigma$-stable decomposition}
Let $Q$ be a quiver without oriented cycles.
\begin{definition}
Suppose that
$\alpha$ is a dimension vector and $\sigma$ is a weight such that
$\sigma(\alpha)=0$.
A dimension vector $\alpha$ is called $\sigma$-(semi-)stable if a 
general representation
of dimension $\alpha$ is $\sigma$-(semi-)stable. 
The expression
$$
\alpha=\alpha_1\pp\alpha_2\pp\ldots \pp\alpha_s
\label{eqStable2}
$$
is called
the {\it $\sigma$-stable decomposition\/}
of a $\sigma$-semi-stable dimension vector $\alpha$ if a
general representation $V$ of dimension $\alpha$ has a
 Jordan-H\"older filtration
with factors of dimension $\alpha_1,\alpha_2,\dots,\alpha_s$
(in some order).
\end{definition}
\begin{proposition}
A dimension vector $\alpha$ is a Schur root if and only if $\alpha$ is $\sigma$-stable
for some weight $\sigma$.
\end{proposition}
\begin{proof}
If $\alpha$ is $\sigma$-stable, then a general  $\alpha$-dimensional representation
is $\sigma$-stable, hence indecomposable. Conversely, if $\alpha$ is a Schur root,
then $\alpha$ is $\sigma$-stable where 
$$\sigma=\langle\alpha,\cdot\rangle-\langle\cdot,\alpha\rangle.$$
This follows from \cite[Theorem 6.1]{S1}.
\end{proof}
\begin{definition}
Suppose that
$\alpha$ is a dimension vector. Let $\sigma,\tau$ be weights.
such that $\tau$ is positive on nonzero dimension vectors.
A dimension vector $\alpha$ is called $(\sigma:\tau)$-stable if a 
general representation
of dimension $\alpha$ is $(\sigma:\tau)$-stable. 
The expression
$$
\alpha=\alpha_1\pp\alpha_2\pp\ldots \pp\alpha_s
\label{eqStable}
$$
is called
the {\it $(\sigma:\tau)$-stable decomposition\/}
of $\alpha$ if a
general representation $V$ of dimension $\alpha$ has a
 HNJH-filtration
with factors of dimension $\alpha_1,\alpha_2,\dots,\alpha_s$
(in some order).
\end{definition}
\begin{lemma}\label{lemsubsum}
Let $\sigma,\tau$ be weights such that $\tau$ is positive
on nonzero dimension vectors.
Suppose that $\alpha\in \NN^{Q_0}$ has $(\sigma:\tau)$-stable decomposition
$$
\alpha_1\pp \alpha_2\pp \cdots \pp \alpha_s.
$$
For any indices $i_1<\cdots<i_r$ the $(\sigma:\tau)$-stable decomposition
of $\beta=\alpha_{i_1}+\alpha_{i_2}+\cdots+\alpha_{i_r}$
is
$$
\alpha_{i_1}\pp \alpha_{i_2}\pp \cdots \pp \alpha_{i_r}.
$$
\end{lemma}
\begin{proof}
After rearranging $\alpha_1,\dots,\alpha_s$ we may assume that
$$
U_\alpha:=U_{HNJH}(\alpha_1+\cdots+\alpha_s,\alpha_1+\cdots+\alpha_{s-1},\dots,\alpha_1,0)
$$
lies dense in $\Rep(Q,\alpha)$. By Lemma~\ref{lemlocallyclosed}, 
$U_{\alpha}$ is
locally closed. Therefore, $U_{\alpha}$ is dense and {\it open}.

Let $j_1<j_2<\cdots<j_{s-r}$ such that
$$
\{1,2,\dots,s\}=\{i_1,\dots,i_r\}\cup \{j_1,\dots,j_{s-r}\}.
$$
Then we have
$$
\alpha-\beta=\alpha_{j_1}+\cdots+\alpha_{j_{s-r}}
$$
Suppose that the $(\sigma:\tau)$-stable decomposition of $\beta$ is
$$
\beta_1\pp \beta_2\pp \cdots \pp \beta_t.
$$

Choose $V_{i}\in \Rep(Q,\alpha_i)$  $(\sigma:\tau)$-stable for $i=1,2,\dots,s$.
By Lemma~\ref{lemlocallyclosed}, there exists an open dense set $U_\beta\subseteq \Rep(Q,\beta)$
and dimension vectors $\gamma_0=\alpha,\gamma_1,\dots,\gamma_{l+1}=0$ such that
$$
V'\oplus V_{j_1}\oplus \cdots \oplus V_{j_{s-r}}\in U_{\alpha}':=
U_{HNJH}(\gamma_0,\dots,\gamma_{l+1})
$$
for all $V'\in U_{\beta}$. Note that
$$
\gamma_0-\gamma_1,\gamma_1-\gamma_2,\dots,\gamma_l-\gamma_{l+1}
$$
is a rearrangement of
$$
\beta_1,\dots,\beta_t,\alpha_{j_1},\dots,\alpha_{j_{s-r}}.
$$
We have
$$
V_1\oplus V_2\oplus \cdots\oplus V_s\in 
\overline{U_{\beta}\oplus V_{j_1}\oplus \cdots \oplus V_{j_{s-r}}}\subseteq
\overline{U_{\alpha}'}.
$$
This shows that $\overline{U_{\alpha}'}\cap U_{\alpha}\neq \emptyset$.
Since $U_{\alpha}$ is open, we have $U_{\alpha}'\cap U_{\alpha}\neq \emptyset$.
Therefore $U_{\alpha}'=U_{\alpha}$. Hence
$$
\alpha_1,\alpha_2,\dots,\alpha_r
$$
is a rearrangement of
$$
\alpha_{j_1},\dots,\alpha_{j_{r-s}},\beta_1,\dots,\beta_t
$$
and
$$\beta_1,\dots,\beta_t
$$
is a rearrangement of
$$
\alpha_{i_1},\dots,\alpha_{i_{r}}.
$$
\end{proof}
\begin{corollary}\label{corsubsum}
Let $\sigma$ be a weight and $\alpha$ be a $\sigma$-semi-stable
dimension vector such that $\sigma(\alpha)=0$.
Suppose that $\alpha\in \NN^{Q_0}$ has $\sigma$-stable decomposition
$$
\alpha_1\pp \alpha_2\pp \cdots \pp \alpha_s.
$$
For any indices $i_1<\cdots<i_r$ the $\sigma$-stable decomposition
of $\beta=\alpha_{i_1}+\alpha_{i_2}+\cdots+\alpha_{i_r}$
is
$$
\alpha_{i_1}\pp \alpha_{i_2}\pp \cdots \pp \alpha_{i_r}.
$$
\end{corollary}
\begin{proof}
This follows from Lemma~\ref{lemsubsum} and Remark~\ref{sigmasigmatau}.
\end{proof}

\begin{proposition}\label{multiplestable}
Suppose that $\alpha$ is a $\sigma$-stable dimension vector and $p\in \NN$. Then
the $\sigma$-stable decomposition of $p\alpha$ is
$$
\{p\alpha\}:=\left\{
\begin{array}{ll}
p\cdot \alpha:=\underbrace{\alpha\pp\cdots\pp\alpha}_p &
\mbox{if $\alpha$ is real or isotropic}\\
p\alpha &\mbox{if $\alpha$ is imaginary and nonisotropic.}
\end{array}\right.
$$
\end{proposition}
\begin{proof}
Suppose that the $\sigma$-stable decomposition of $p\alpha$
is
$$
\gamma_1\pp\gamma_2\pp\cdots \pp \gamma_r
$$
Let $W\in \Rep(Q,\alpha)$ be a $\sigma$-stable representation.
Then $V=W^p=W\oplus \cdots \oplus W$ is $\sigma$-semi-stable
and it has a filtration
with factors of dimensions $\gamma_1,\dots,\gamma_r$ (see Lemma~\ref{FilterClosed}).
This filtration can be refined to a Jordan-H\"older filtration of $V$
whose factors are all isomorphic to the $\alpha$-dimensional representation $W$.
This shows that $\gamma_1,\dots,\gamma_r$ are multiples of $\alpha$.
If $\alpha$ is real or isotropic, then the canonical decomposition
of $p\alpha$ is $\alpha^{\oplus p}$. In particular, a general representation $Z$
of dimension $p\alpha$ has a filtration $\{F^i(Z)\}$ with factors of dimension $\alpha$.
This filtration can be refined to a J\"ordan H\"older filtration.
But since the factors of the J\"ordan H\"older filtration have dimensions
that are {\it multiples} of $\alpha$, the filtration $\{F^i(Z)\}$ is already
a Jordan-H\"older filtration.
Therefore $r=p$ and $\gamma_1,\dots,\gamma_r=\alpha$.

Suppose $\alpha$ is imaginary and not isotropic. Then $\gamma_i\hookrightarrow p\alpha$
for some $i$. Assume that $\gamma_i=q\alpha$. Then $q\alpha\hookrightarrow p\alpha$,
so $\ext_Q(q\alpha,(p-q)\alpha)=0$. It follows that 
$$q(p-q)\langle \alpha,\alpha\rangle=
\langle q\alpha,(p-q)\alpha\rangle\geq 0.$$
Since $\langle\alpha,\alpha\rangle<0$, we must have $p=q$.
It follows that $i=1$ and $\gamma_1=p\alpha$.
\end{proof}
\begin{proposition}\label{multiples}
If the $\sigma$-stable decomposition of $\alpha$ is
$$
\alpha=\alpha_1\pp\alpha_2\pp\cdots\pp\alpha_s
$$
then the $\sigma$-stable decomposition of $p\alpha$ is
\begin{equation}
p\alpha=\{p\alpha_1\}\pp\{p\alpha_2\}\pp\cdots\pp\{p\alpha_s\}\label{eqpalph}
\end{equation}
for every positive integer $p$.
\end{proposition}
\begin{proof}
If $\beta\hookrightarrow\alpha$ then $\ext_Q(\beta,\alpha-\beta)=0$.
It follows that $\ext_Q(p\alpha,p(\alpha-\beta))=0$ and $m\beta\hookrightarrow m\alpha$.
Suppose that $V\in \Rep(Q,p\alpha)$. Then $V$ has a filtration
with factors $V_1,V_2,\dots,V_s$ of dimensions $p\alpha_1,\dots,p\alpha_s$.
From Proposition~\ref{multiplestable} and Lemma~\ref{FilterClosed} follows
that if $\alpha_i$ is isotropic or real, then
each $V_i$ has a filtration such that all quotients
have dimension $\alpha_i$. This shows that $V$ has a filtration $\{F^i(V)\}$
such that the dimensions of the factors are exactly all dimensions
appearing on the right-hand side of (\ref{eqpalph}).
Assume now that
$V$ is a general representation of dimension $p\alpha$. Then all factors of the filtration $\{F^i(V)\}$
are $\sigma$-semi-stable. There exists a Jordan-H\"older filtration
$\{(F')^i(V)\}$ which is a refinement of the filtration $\{F^i(V)\}$.
If $W$ is any representation of dimension $p\alpha$, then
$W$ has a filtration $\{(F')^i(W)\}$ such that $\dd (F')^i(W)=\dd (F')^i(V)$ by Lemma~\ref{FilterClosed}.
Take now
$W=W_1\oplus W_2\oplus \cdots\oplus W_s$ with $W_1,\dots,W_s$
in general position of dimension $p\alpha_1,\dots,p\alpha_s$.
The filtration $\{(F')^i(W)\}$ can be refined to a Jordan-H\"older filtration
of $W$. But we know that the dimensions of the factors of the Jordan-H\"older filtration
of $W$ are exactly the dimensions appearing on the right-hand side of (\ref{eqpalph}) by Proposition~\ref{multiplestable}.
It follows that the filtrations $\{F^i(W)\}$ and $\{(F')^i(W)\}$ are the same,
so the filtrations $\{F^i(V)\}$ and $\{(F')^i(V)\}$ are the same. 
This shows that the factors of the Jordan-H\"older filtration of $V$
are exactly given by the dimensions on the right-hand side of (\ref{eqpalph}).
\end{proof}
\begin{proposition}\label{propstable}
Let $Q$ be a quiver without oriented cycles, $\sigma$ be a weight
and $\alpha$ be a $\sigma$-semi-stable dimension vector.
Let
$$
\alpha =c_1\cdot\alpha_1\pp c_2\cdot \alpha_2\pp \cdots \pp c_s\cdot\alpha_s
$$
be the $\sigma$-stable decomposition where the  $\alpha_i$ are distinct. Then we have
\begin{enumerate}
\renewcommand{\theenumi}{\alph{enumi}}
\item  All $\alpha_i$ are Schur roots;
\item if $c_i>1$ then $\alpha_i$ must be a real or an isotropic root;
\item $\hom_Q(\alpha_i,\alpha_j)=0$ if $i\neq j$;
\item after rearranging one may assume that $\alpha_i\pperp \alpha_j$
for all $i<j$.
\end{enumerate}
\end{proposition}
\begin{proof}\ 

(a) If  general representations of dimension $\alpha_i$ are decomposable,
then a general representation of dimension $\alpha_i$ is not $\sigma$-stable.
Since a general representation of dimension $\alpha_i$ has a non-trivial
Jordan-H\"older filtration, {\it every\/} representation of
dimension $\alpha_i$ has a nontrivial Jordan-H\"older filtration by
Lemma~\ref{FilterClosed}.
It follows that there are no $\sigma$-stable representations
in dimension $\alpha_i$.

(b) If $c_i\geq 2$, then by Corollary~\ref{corsubsum}, the $\sigma$-stable
decomposition of $2\alpha_i$ is $\alpha_i\pp\alpha_i$.
In particular we have $\alpha_i\hookrightarrow 2\alpha_i$, so 
$\ext_Q(\alpha_i,\alpha_i)=0$ and $\langle \alpha_i,\alpha_i\rangle \geq 0$.

(c) Let $V_i$ and $V_j$ be $\sigma$-stable representations of
dimensions $\alpha_i$ and $\alpha_j$ respectively. From Lemma~\ref{lemmorphism}
follows that any nonzero homomorphism between $V_i$ and $V_j$
must be an isomorphism. Since $\alpha_i\neq\alpha_j$, we have
$\Hom_Q(V_i,V_j)=0$ and therefore $\hom_Q(\alpha_i,\alpha_j)=0$.

(d) Suppose that there exists an $r$-cycle: 
$$
\ext_Q(\alpha_{i_1},\alpha_{i_2})\neq 0,\ \ext_Q(\alpha_{i_2},\alpha_{i_3})\neq 0,\dots,
\ \ext_Q(\alpha_{i_{r-1}},\alpha_{i_r})\neq 0,\ \ext_Q(\alpha_{i_r},\alpha_{i_1})\neq 0.
$$
We assume that $r\geq 2$ is minimal such that an $r$-cycle exists.
After rearranging $\alpha_1,\alpha_2,\dots,\alpha_r$ we may assume that
$$
\ext_Q(\alpha_{1},\alpha_{2})\neq 0,\ \ext_Q(\alpha_{2},\alpha_{3})\neq 0,\dots,
\ \ext_Q(\alpha_{r-1},\alpha_r)\neq 0,\ \ext_Q(\alpha_r,\alpha_1)\neq 0.
$$
Also, by the minimality of $r$, 
 $\ext_Q(\alpha_i,\alpha_j)=0$ if $1\leq i,j\leq r$, unless
 $j=i$ or $j\equiv i+1\ (\bmod\ r)$.
The $\sigma$-stable decomposition of
 $\beta=\alpha_1+\alpha_2+\cdots+\alpha_r$ is
$$
\alpha_1\pp\alpha_2\pp\cdots\pp\alpha_r
$$
by Corollary~\ref{corsubsum}.
We have $\alpha_i\hookrightarrow \beta$ for some $i$, and after
cyclic relabeling we may assume that $\alpha_1\hookrightarrow \beta$,
so $\ext_Q(\alpha_1,\beta-\alpha_1)=0$.
If $3\leq i\leq r$, then we have $\ext_Q(\alpha_i,\alpha_2)=0$ by minimality of
$r$.
It follows that $\ext_Q(\alpha_3+\alpha_4+\cdots+\alpha_r,\alpha_2)=0$
or equivalently, $\beta-\alpha_1\twoheadrightarrow \alpha_2$.
Consider an exact sequence
\begin{equation}\label{exact}
0\to V'\to V\to V''\to 0
\end{equation}
with $V',V,V''$ of dimension $\beta-\alpha_1-\alpha_2,\beta-\alpha_1,\alpha_2$
respectively, and $V$ is in general position.
Let $W$ be a general representation of dimension $\alpha_1$ and
apply the functor $\Hom_Q(W,\cdot)$ to (\ref{exact}) to obtain a long
exact sequence
$$
\cdots \to \Ext_Q(W,V')\to \Ext_Q(W,V)\to\Ext_Q(W,V'')\to 0.
$$
Since $\Ext_Q(W,V)=0$, we have $\Ext_Q(W,V'')=0$. It follows that
$\ext_Q(\alpha_1,\alpha_2)=0$.
Contradiction.

 Therefore, there are no $r$-cycles.
So it is possible to
rearrange $\alpha_1,\dots,\alpha_s$ such that $\ext_Q(\alpha_i,\alpha_j)=0$
for all $i<j$. Together with $\hom_Q(\alpha_i,\alpha_j)=0$ by (c)
we get $\alpha_i\perp \alpha_j$ for $i<j$.

The number $p=\alpha_i\circ\alpha_j$ is finite and nonzero.
We would like to show that $p=1$. Suppose that $p\geq 2$.
Let $V$ be a general representation of dimension $\alpha_i+\alpha_j$.
Then $V$ 
has $p\geq 2$ subrepresentations of dimension $\alpha_i$ (see Theorem~\ref{theoDSW}).
Let $V_1$ and $V_2$ be two distinct subrepresentations of $V$
of dimension $\alpha_i$. The $\sigma$-stable decomposition of $\alpha_i+\alpha_j$
is $\alpha_i\pp\alpha_j$ by 
Corollary~\ref{corsubsum}.
This implies that $V_1,V_2,V/V_1,V/V_2$ are $\sigma$-stable.
Let $\varphi:V_1\to V/V_2$ be the composition
$V_1\to V\to V/V_2$. Since $V_1$ and $V_2$ are distinct of the same dimension,
$\varphi$ is not equal the 0 map. Since $V_1$, $V/V_2$ are $\sigma$-stable
$\varphi$ must be an isomorphism. This implies that $\alpha_i=\alpha_j$.
Contradiction.
So we conclude that $p=\alpha_i\circ\alpha_j=1$ and $\alpha_i\pperp\alpha_j$.
\end{proof}
\begin{corollary}\label{corstable}
Let $Q$ be a quiver without oriented cycles,
$\sigma,\tau$ two weights with $\tau$ nonnegative
on nonzero dimension vectors, and let $\alpha$ be
a dimension vector. Let
$$
\alpha =c_1\cdot\alpha_1\pp c_2\cdot \alpha_2\pp \cdots \pp c_s\cdot\alpha_s
$$
be the $(\sigma:\tau)$-stable decomposition where the  $\alpha_i$ are distinct.
\begin{enumerate}
\renewcommand{\theenumi}{\alph{enumi}}
\item  All $\alpha_i$ are Schur roots;
\item if $c_i>1$ then $\alpha_i$ must be a real root or an isotropic root;
\item after rearranging one can assume that $\alpha_i\pperp \alpha_j$
for all $i<j$.
\end{enumerate}
\end{corollary}
\begin{proof}
(a) and (b) follow from Proposition~\ref{propstable} (a), (b).
Let $\mu=\sigma/\tau$. We may assume that $\mu(\alpha_i)\geq \mu(\alpha_j)$
for all $i<j$. We have $\hom_Q(\alpha_i,\alpha_j)=0$ for $i<j$ by Lemma~\ref{lemmorphism}(b).
Suppose that $\mu(\alpha_i)>\mu(\alpha_j)$. The $(\sigma:\tau)$-stable
decomposition of $\alpha_i+\alpha_j$ is $\alpha_i\pp\alpha_j$ (see
Lemma~\ref{lemsubsum}). 
We have either $\alpha_i\hookrightarrow \alpha_i+\alpha_j$
or $\alpha_j\hookrightarrow \alpha_i+\alpha_j$. Because
$\mu(\alpha_i)>\mu(\alpha_j)$, $\alpha_i$ (and not
$\alpha_j$) is the dimension of the submodule in
the Harder-Narasimhan filtration
of $\alpha_i+\alpha_j$, so $\alpha_i\hookrightarrow \alpha_i+\alpha_j$.
This shows that $\ext_Q(\alpha_i,\alpha_j)=0$, and therefore $\alpha_i\perp \alpha_j$.
Since the Harder-Narasimhan filtration is unique, a general representation
of $\alpha_i+\alpha_j$ has a unique $\alpha_i$-dimensional subrepresentation.
We conclude that $\alpha_i\pperp \alpha_j$
by Theorem~\ref{theoDSW}.

We can arrange $\alpha_1,\alpha_2,\dots,\alpha_s$ such that
$\mu(\alpha_1)\geq\mu(\alpha_2)\geq \cdots \geq \mu(\alpha_s)$. 
Using Proposition~\ref{propstable}(d)  and Remark~\ref{sigmasigmatau} 
we see that we can even rearrange
$\alpha_1,\dots,\alpha_s$ such that $\alpha_i\pperp \alpha_j$ for all $i<j$.
\end{proof}

\begin{theorem}\label{invar}
Suppose that $\sigma$ is an indivisible weight.
If 
$$\alpha=c_1\cdot \alpha_1\pp c_2\cdot \alpha_2\pp\cdots \pp c_r\cdot\alpha_r$$
is the $\sigma$-stable decomposition of $\alpha$, then
there exists an isomorphism
$$
\SI(Q,\alpha)_{m\sigma}\cong S^{c_1}(\SI(Q,\alpha_1)_{m\sigma})\otimes 
S^{c_2}(\SI(Q,\alpha_2)_{m\sigma})
\otimes\cdots\otimes S^{c_r}(\SI(Q,\alpha_r)_{m\sigma})
$$
\end{theorem}
\begin{proof}
First we assume that the base field $K$ has characteristic 0.
Let
$$
X:=\Rep(Q,\alpha_1)^{c_1}\oplus \Rep(Q,\alpha_2)^{c_2}\oplus \cdots\oplus 
\Rep(Q,\alpha_r)^{c_r}.
$$
We have a natural embedding
$$
\varphi:X\hookrightarrow \Rep(Q,\alpha).
$$
Let $\GL(Q,\alpha)_{\sigma}$ be the kernel
of $\sigma$, where $\sigma$ is interpreted as
a multiplicative character $\GL(Q,\alpha)\to K^\star$.
Let $G$ be the stabilizer of $X$ within $\GL(Q,\alpha)_{\sigma}$.
This group $G$ is isomorphic to the
intersection of 
$$
(\Sigma_{c_1}\ltimes \GL(Q,\alpha_1)^{c_1})\times  
(\Sigma_{c_2}\ltimes \GL(Q,\alpha_2)^{c_2})\times\cdots
\times
(\Sigma_{c_r}\ltimes \GL(Q,\alpha_r)^{c_r})$$ 
and $\GL(Q,\alpha)_\sigma$. Here $\Sigma_c$ is the symmetric group on 
$c$ elements.
Let $\pi_X:X\to X\quo G$ and
$\pi:\Rep(Q,\alpha)\to \Rep(Q,\alpha)\quo \GL(Q,\alpha)_\sigma$
be the categorical quotients.
The embedding $\varphi$ induces a morphism between categorical quotients
$$
\psi:X\quo G\to \Rep(Q,\alpha)\quo \GL(Q,\alpha)_\sigma.
$$
We will show that $\psi$ is an isomorphism.

First we will show that $\psi$ is dominant. 
Let $V\in \Rep(Q,\alpha)$ be a general
representation and suppose $V$ is $\sigma$-semi-stable, i.e.,
$\pi(V)\neq 0$.
Let
$$
0=F^0_J(V)\subset F^1_J(V)\subset\cdots\subset F^s_J(V)=V
$$
be a Jordan-H\"older filtration with $\sigma$-stable quotients
$G^i(V)=F^i(V)/F^{i-1}(V)$.
Now $W:=\bigoplus_i G^i(V)\in X$ lies in the $\GL(Q,\alpha)_\sigma$
closure of $V$. In particular $\psi(\pi_X(W))=\pi(V)$. This proves that $\psi$ is
surjective.

Note that $W\in X$ is $G$-semi-stable if and only if all the
summands in $\Rep(Q,\alpha_i)$ are $\sigma$-semi-stable.
In particular, if $W\in X$ is $G$-semi-stable, then
$W$ $\sigma$-semi-stable, so $W$ is $\GL(Q,\alpha)_\sigma$-stable.
This shows
that $\psi^{-1}(0)=\{0\}$, so $\psi$ is a finite map.

For a general representation $V\in \Rep(Q,\alpha)$, the
quotients of a Jordan-H\"older filtration corresponding to $\sigma$
are unique up to permutation. This shows that a generic
fiber of $\psi$ consists of only one point. So $\psi$ is birational.

Because $\psi$ is birational and finite, and $\Rep(Q,\alpha)\quo \GL(Q,\alpha)_{\sigma}$
is normal, $\psi$ must be an isomorphism.
Now the graded coordinate ring
of $\Rep(Q,\alpha)\quo \GL(Q,\alpha)_{\sigma}$ is
$\bigoplus_{m\geq 0} \SI(Q,\alpha)_{m\sigma}$
and $X\quo G$ has graded coordinate ring
$$
\bigoplus_{m\geq 0}
S^{c_1}(\SI(Q,\alpha_1)_{m\sigma})\otimes 
S^{c_2}(\SI(Q,\alpha_2)_{m\sigma})
\otimes\cdots\otimes S^{c_r}(\SI(Q,\alpha_r)_{m\sigma}).
$$
In particular we get
$$
\SI(Q,\alpha)_{m\sigma}\cong
S^{c_1}(\SI(Q,\alpha_1)_{m\sigma})\otimes 
S^{c_2}(\SI(Q,\alpha_2)_{m\sigma})
\otimes\cdots\otimes S^{c_r}(\SI(Q,\alpha_r)_{m\sigma})
$$
for all $m\geq 0$.

Suppose that $K$ has arbitrary characteristic. 
We note that $\dim \SI(Q,\beta)_{\tau}$
is independent of the base field $K$ (see~\cite{DSW})
for any dimension vector $\beta$ and any weight $\tau$.
The vector spaces 
$\SI(Q,\alpha)_{m\sigma}$
and
$$
S^{c_1}(\SI(Q,\alpha_1)_{m\sigma})\otimes 
S^{c_2}(\SI(Q,\alpha_2)_{m\sigma})
\otimes\cdots\otimes S^{c_r}(\SI(Q,\alpha_r)_{m\sigma})
$$
have the same dimension if $K$ has characteristic 0. But
then they have the same dimension even in the case where $K$ has
positive characteristic. We have an isomorphism
$$
\SI(Q,\alpha)_{m\sigma}\cong
S^{c_1}(\SI(Q,\alpha_1)_{m\sigma})\otimes 
S^{c_2}(\SI(Q,\alpha_2)_{m\sigma})
\otimes\cdots\otimes S^{c_r}(\SI(Q,\alpha_r)_{m\sigma})
$$
although this isomorphism may not be canonical. If we
replace the symmetric powers $S^c$ by the divided powers $D^c$, then
the argument in characteristic $0$ still shows that we
have a canonical injective map
$$
\SI(Q,\alpha)_{m\sigma}\hookrightarrow
D^{c_1}(\SI(Q,\alpha_1)_{m\sigma})\otimes 
D^{c_2}(\SI(Q,\alpha_2)_{m\sigma})
\otimes\cdots\otimes D^{c_r}(\SI(Q,\alpha_r)_{m\sigma})
$$
Comparing dimensions shows that this map is an isomorphism.
\end{proof}

\section{Schur sequences}\label{sec4}
Proposition~\ref{propstable} and Corollary~\ref{corstable}
motivates us to define  Schur sequences.
The notion of Schur sequences is closely related to the notion
of exceptional sequences (see Section~\ref{secexseq}).
A Schur sequence is similar to an exceptional sequence, but also
imaginary Schur roots are allowed.
\begin{definition}
A sequence $\underline{\gamma}=(\gamma_1,\gamma_2,\dots,\gamma_s)$
is called a Schur sequence if
\begin{enumerate}
\renewcommand{\theenumi}{\arabic{enumi}}
\item $\gamma_i$ is a Schur root for every $i$;
\item $\gamma_i\pperp \gamma_j$ for all $i<j$.
\end{enumerate}
\end{definition}
\begin{lemma}\label{lemdense}
If $\Rep(Q,\alpha)$ has a dense $\GL(Q,\alpha)$-orbit or
$\Rep(Q,\beta)$ has a dense $\GL(Q,\beta)$-orbit, then
$\alpha\circ \beta\leq 1$.
\end{lemma}
\begin{proof}
If $\Rep(Q,\beta)$ has a dense $\GL(Q,\beta)$-orbit, then
there are no rational $\GL(Q,\beta)$-invariants in $K[\Rep(Q,\beta)]$.
In particular, any quotient of two semi-invariants of the same
weight must be constant. This shows that
$$
\alpha\circ\beta=\dim \SI(Q,\beta)_{\langle \alpha,\cdot\rangle}\leq 1.
$$
If $\Rep(Q,\alpha)$ has a dense $\GL(Q,\alpha)$-orbit then the proof is similar.
\end{proof}
\begin{remark}
An exceptional sequence $\underline{\varepsilon}=(\varepsilon_1,\varepsilon_2,\dots,\varepsilon_s)$
is a Schur sequence. The space $\Rep(Q,\varepsilon_j)$ has a dense
$\GL(\varepsilon_j)$-orbit, so by Lemma~\ref{lemdense}
$$
\varepsilon_i\circ \varepsilon_j=
\dim \SI(Q,\varepsilon_j)_{\langle \varepsilon_i,\cdot\rangle}\leq 1.
$$
This shows that $\varepsilon_i\pperp \varepsilon_j$ for all $i<j$.
\end{remark}

\begin{lemma}\label{lemex}
Suppose that $\alpha\perp \gamma$ and $\beta\perp\gamma$. 
\begin{enumerate}
\renewcommand{\theenumi}{\alph{enumi}}
\item If $(\alpha+\beta)\circ \gamma=1$,
then $\alpha\circ\gamma=1$ and $\beta\circ\gamma=1$;
\item if $\ext_Q(\alpha,\beta)=0$, $\alpha\circ\gamma=1$ and $\beta\circ\gamma=1$
then $(\alpha+\beta)\circ \gamma=1$.
\end{enumerate}
\end{lemma}
\begin{proof}
(a) Since $\alpha\perp\gamma$ and $\beta\perp\gamma$
we have $\alpha\circ\gamma\geq 1$ and $\beta\circ\gamma\geq 1$.

Choose $f\in \SI(Q,\gamma)_{\langle \alpha,\cdot\rangle}$.
Then we have
$$
f\SI(Q,\gamma)_{\langle \beta,\cdot\rangle}\subseteq 
\SI(Q,\gamma)_{\langle \alpha+\beta,\cdot\rangle}.
$$
This shows that 
$$1\leq \beta\circ\gamma=\dim\SI(Q,\gamma)_{\langle \beta,\cdot\rangle}=
\dim f\SI(Q,\gamma)_{\langle \beta,\cdot\rangle}\leq
\dim SI(Q,\gamma)_{\langle\alpha+\beta,\cdot\rangle}=
(\alpha+\beta)\circ \gamma\leq 1.$$
Hence we get $\beta\circ\gamma=1$.
Similarly we obtain $\alpha\circ\gamma=1$.

(b) Any $(\alpha+\beta)$-dimensional representation $V$ has
an $\alpha$-dimensional subrepresentation $V'$. If $V''=V/V'$ then 
$c^V=c^{V'}\cdot c^{V''}$ up to a scalar
as functions on $\Rep(Q,\gamma)$ by Lemma~\ref{lemexseq}.
Since $\SI(Q,\gamma)_{\langle\alpha+\beta,\cdot\rangle}$ is spanned
by semi-invariants of the form $c^V$ (see Theorem~\ref{theoDW}),
we have
$$
\SI(Q,\gamma)_{\langle\alpha+\beta,\cdot\rangle}=
\SI(Q,\gamma)_{\langle\alpha,\cdot\rangle}
\SI(Q,\gamma)_{\langle\beta,\cdot\rangle}.$$
It follows that
$$
1\leq (\alpha+\beta)\circ \gamma\leq
(\alpha\circ\gamma)(\beta\circ\gamma)=1\cdot 1=1,
 $$
so $(\alpha+\beta)\circ\gamma=1$. 
\end{proof}
\begin{corollary}\label{corantirefine}
If $\gamma_1,\gamma_2,\dots,\gamma_s$ is a Schur sequence,
and $p\gamma_i+q\gamma_{i+1}$ is a Schur root,
then 
$$
\gamma_1,\gamma_2,\dots,\gamma_{i-1},p\gamma_i+q\gamma_{i+1},\gamma_{i+2},\dots,
\gamma_s
$$
is a Schur sequence.
\end{corollary}
\begin{proof}
This follows from Lemma~\ref{lemex}(b) and Theorem~\ref{thmBelkale}.
\end{proof}
\begin{remark}
Suppose that
 $\alpha=\alpha_1^{\oplus c_1}\oplus \alpha_2^{\oplus c_2}\oplus \cdots \oplus
\alpha_s^{\oplus c_s}$ is the canonical decomposition
with all $\alpha_i$ distinct. By Theorem~\ref{theoKac}
$\ext_Q(\alpha_i,\alpha_j)=0$ for all $i\neq j$. After reordering we may assume that 
$\hom_Q(\alpha_i,\alpha_j)=0$ for
all $i<j$ by Lemma~\ref{lemrearrange}. 
So $\alpha_i\perp \alpha_j$ for all $i<j$.
We claim that in fact $\alpha_1,\alpha_2,\dots,\alpha_s$
is a Schur sequence. This follows from the algorithm in \cite{DW2} for
finding the canonical decomposition and Corollary~\ref{corantirefine}.
\end{remark}
\begin{remark}
Suppose that $\alpha=c_1\cdot \alpha_1\pp c_2\cdot \alpha_2\pp \cdots \pp c_r\cdot\alpha_r$
is a $(\sigma:\tau)$-stable decomposition. By Corollary~\ref{corstable} (c)
we may assume that $\alpha_i\pperp\alpha_j$ for all $i<j$.
 By Corollary~\ref{corstable}(a), $\alpha_1,\dots,\alpha_r$ is a Schur sequence.
  \end{remark}

\begin{definition}
A Schur sequence  $\underline{\gamma}=(\gamma_1,\gamma_2,\dots,\gamma_s)$
is called a quiver Schur sequence if $\langle \gamma_j,\gamma_i\rangle\leq 0$
for all $i<j$.
\end{definition}
\begin{remark}
Suppose that $\alpha=c_1\cdot \alpha_1\pp c_2\cdot \alpha_2\pp \cdots \pp c_r\cdot\alpha_r$
is a $\sigma$-stable decomposition. By Proposition~\ref{propstable} (d)
we may assume that $\ext_Q(\alpha_i,\alpha_j)=0$ for all $i<j$.
 By Proposition~\ref{propstable} (a),(c),   $\alpha_1,\alpha_2,\dots,\alpha_r$ is
a quiver Schur sequence.
\end{remark}

\begin{definition}
Let ${\underline \gamma} = (\gamma_1,\ldots,\gamma_r)$, 
${\underline \beta}= (\beta_1,\ldots,\beta_s)$ be two
sequences of dimension vectors. 
We say that $\underline \beta$ is a refinement of $\underline \gamma$ if 
there exists a sequence $0=b_0 <b_1<\ldots<b_{r-1}<s=b_r$ such that for each 
$j=1,\ldots ,r$ the dimension
vector $\gamma_j$ is a positive linear combination of 
$\beta_{b_{j-1}+1},\ldots,\beta_{b_j}$.
\end{definition}
\begin{theorem}[Refinement Theorem]\label{theorefinement}
Let $\underline{\gamma}=(\gamma_1,\dots,\gamma_r)$ be a Schur sequence. 
Then there exists an exceptional
sequence $\underline{\varepsilon}=(\varepsilon_1,\dots,\varepsilon_s)$ 
such that $\underline{\varepsilon}$ is a refinement of $\underline{\gamma}$.
\end{theorem}
\begin{proof} 
For a dimension vector $\alpha$ we define $\tau(\alpha)=\sum_{x\in Q_0} \alpha(x)$.
We will prove the theorem by induction on
the number of vertices $n$ in the quiver $Q$,
and by induction on $\tau(\gamma_1)$.
 If $n=1$ there is nothing to prove.
If $\tau(\gamma_1)=0$ then there is nothing to prove either since this
is impossible.
 
Let us assume that in a Schur sequence the first dimension vector $\gamma_1$ 
is a 
real Schur root. Let $V$ be the unique indecomposable representation
corresponding to the dense orbit in $\Rep(Q,\gamma_1)$.
 Then the dimension
vectors $\gamma_2 ,\ldots ,\gamma_r$ are Schur roots in the right 
orthogonal category 
$V^\perp$. By Theorem~\ref{rightperp}
the category $V^\perp$ is equivalent to the category of representations of a 
quiver $Q^\prime$ with no oriented cycles and $n-1$ vertices. 
Let $I:\NN^{n-1}\to \NN^{Q_0}$ as in Section~\ref{secperp}.
We can write $\gamma_i=I(\delta_i)$. Then $\delta_2,\dots,\delta_n$
is a Schur sequence for $Q'$. We can refine $\delta_2,\dots,\delta_n$
to an exceptional sequence $\underline{\varepsilon}=(\varepsilon_1,\dots,\varepsilon_t)$
by induction on $n$.

Then the sequence 
${\underline \varepsilon}=(\gamma_1 ,I(\varepsilon_2),\ldots,I(\varepsilon_n))$ is 
clearly an exceptional sequence for $Q$ which refines $\underline \gamma$.

Next, assume that $\gamma_1$ is an imaginary root.
 Since $\gamma_1\circ \gamma_j=1$ for all $j\geq 2$,
it follows by induction from Lemma~\ref{lemex}~(b) that 
$\gamma_1\circ \delta=1$  where
 $\delta=\gamma_2+\cdots+\gamma_r$. 
Let $\gamma_1^{\perp}$ be the set of all dimension vectors $\alpha$
with $\gamma_1\perp \alpha$. 
By Theorem~2 from \cite{DW},
$\RR_+\gamma_1^\perp$ is a rational polyhedral cone in $\RR^{n-1}$. 
Suppose that $\delta$ 
is in the interior of the cone.
For each $\alpha\in
\gamma_1^\perp$ there exists $\beta\in \gamma_1^\perp$ such that 
$\alpha+\beta =
p\delta$ for some positive integer $p$.
From Lemma~\ref{lemex}~(a) follows that $\gamma_1\circ\alpha=1$.
This shows that for all $\sigma\in \Gamma^\star$, $\dim \SI(Q,\gamma_1)_{\sigma}\leq 1$.
Indeed, if $\SI(Q,\gamma_1)_\sigma\neq 0$ then
$\sigma=-\langle\cdot,\alpha\rangle$ for some $\alpha\in \gamma_1^\perp$
and $\dim \SI(Q,\gamma_1)_\sigma=\gamma_1\circ \alpha=1$.
This implies that all $\GL(\gamma_1)$-invariant rational functions on $\Rep(Q,\gamma_1)$
are constant. From this follows that $\GL(\gamma_1)$ has a dense
orbit in $\Rep(Q,\gamma_1)$.
Therefore,  $\gamma_1$ must be a 
{\it real\/}
Schur root. Contradiction, so it follows that $\delta$ is not
in the interior of $\RR_+\gamma_1^\perp$.

Let $\sigma = -\langle \cdot,\delta\rangle$
and let
us study the
$\sigma$-stable decomposition of $\gamma_1$. 
By Theorem~\ref{theoDW} there exists a $\beta\hookrightarrow\gamma_1$
such that $\sigma(\beta)=0$. In particular, the $\sigma$-stable decomposition
of $\gamma_1$ is nontrivial. 
Suppose that 
$$
\gamma_1=c_1\cdot \beta_1\pp c_2\cdot \beta_2\pp\cdots\pp c_l\cdot\beta_l
$$
is the $\sigma$-stable decomposition of $\gamma_1$.
We may assume that $\beta_i\pperp \beta_j$ for $i<j$.
From $\gamma_1\pperp \delta=1$ and Lemma~\ref{lemex} follows
that $\beta_i\pperp \gamma_j$ for all $j\geq 2$ and all $i$.
Therefore
\begin{equation}\label{eqS}
\underline{\gamma'}=(\beta_1,\beta_2,\dots,\beta_l,\gamma_2,\dots,\gamma_r)
\end{equation}
is a Schur sequence using Lemma~\ref{lemex}~(a). Notice that $\beta_1$ is
smaller than $\gamma_1$.
Now $\tau(\beta_1)<\tau(\gamma_1)$
so by induction there exists an exceptional sequence which is
a refinement of $\underline{\gamma'}$.
\end{proof}
\begin{corollary}
Let ${\underline \gamma}=(\gamma_1,\ldots,\gamma_r)$ be a Schur sequence. 
Then the vectors $\gamma_1,\ldots ,\gamma_r$ are linearly independent.
\end{corollary}
\begin{proof}
First assume that $\underline{\gamma}$ is an exceptional sequence.
We have
$$
\langle \gamma_i,\gamma_j\rangle=\left\{
\begin{array}{ll}
1 & \mbox{if $i=j$,}\\
0 & \mbox{if $i<j$.}
\end{array}\right.
$$
The matrix
$$
\big(\langle \gamma_i,\gamma_j\rangle\big)_{i,j=1}^n
$$
is invertible, so $\gamma_1,\dots,\gamma_r$ are linearly independent.

If $\underline{\gamma}$ is not an exceptional sequence
then it has a refinement which is an exceptional sequence.
So again, it follows that $\gamma_1,\dots,\gamma_r$ are linearly independent.

\end{proof} 

\section{The faces of the cone $\RR_+\Sigma(Q,\alpha)$}\label{sec5}
As before, $\Sigma(Q,\alpha)$ denotes the set of all weights $\sigma\in \Gamma^\star$
such that $\alpha$ is $\sigma$-semi-stable. If $\alpha$ is a sincere dimension vector,
(i.e., $\alpha(x)>0$ for all $x\in Q_0$,) then there is
a bijection between $\alpha^\perp=\{\beta\in \NN^{Q_0}\mid \alpha\perp \beta\}$
and $\Sigma(Q,\alpha)$ by
$$
\beta\in \alpha^\perp\leftrightarrow \sigma:=-\langle\cdot,\beta\rangle\in \Sigma(Q,\alpha).
$$
Similarly there is also a bijection between $\Sigma(Q,\alpha)$ and ${}^\perp\alpha$.

Let $\RR_+$ be the set of nonnegative
real numbers. We consider the cone $\RR_+\Sigma(Q,\alpha)\subseteq 
\RR_+\otimes_\ZZ \Gamma^\star\cong (\RR^n)^\star$. This cone is of particular interest.
For example, in a special case, this cone is equal to the  Klyachko cone. 
In this section, we will unravel the geometry of this cone.

The refinement Theorem (Theorem~\ref{theorefinement}) allows us to obtain a beautiful 
description of the faces of  $\RR_+\Sigma(Q,\alpha)$.
Let us denote by ${\mathcal W}_r(Q, \alpha)$  the set of all
sets $\{\gamma_1,\dots,\gamma_r\}$ such that
 ${\underline \gamma}=(\gamma_1,\ldots,\gamma_r)$ is a
 quiver Schur sequence of length $r$ such that 
$\alpha=\sum_{i=1}^r a_i \gamma_i$ 
 with 
\begin{enumerate}
\renewcommand{\theenumi}{\arabic{enumi}}
\item $a_i$ a positive integer for all $i$;
\item if $\gamma_i$ is imaginary and not isotropic, then $a_i=1$.
\end{enumerate}
Let ${\mathcal F}_r (Q, \alpha )$ be the set of
faces of dimension $n-r$ of $\RR_+\Sigma (Q,\alpha)$.
\begin{theorem}\label{bijection}
Let $Q$ be a quiver without oriented cycles and let $\alpha$ be a 
dimension vector.
For each $r$, $1\le r\le n-1$ there is a natural bijection
$$\psi(r): {\mathcal W}_r(Q,\alpha)\rightarrow {\mathcal F}_r (Q,\alpha)$$ 
which sends the quiver Schur sequence 
${\underline \gamma}=(\gamma_1 ,\ldots ,\gamma_r )$ to the face
$$\RR_+\Sigma(Q,\gamma_1)\cap\ldots\cap\RR_+ \Sigma(Q,\gamma_r)=
\RR_+\Sigma(Q,\alpha)\cap 
\{
\sigma\in (\RR^n)^\star\mid  \sigma(\gamma_1)=\cdots=\sigma(\gamma_r)=0\}
.$$
The inverse bijection is obtained as follows. 
For a given face $F$ we take a weight $\sigma\in \Gamma^\star$ in the
relative interior of $F$ and associate to it the quiver Schur sequence coming 
from the $\sigma$-stable decomposition of $\alpha$.
\end{theorem}
\begin{proof} Suppose that 
${\underline \gamma}= (\gamma_1 ,\ldots ,\gamma_r )\in
{\mathcal W}_r(Q,\alpha)$.
Let us first prove  that 
$$\RR_+\Sigma(Q,\gamma_1)\cap\ldots\cap \RR_+\Sigma(Q,\gamma_r)$$
is a face of codimension $r$ in the space of dimension vectors. 
Clearly
$$
\Sigma(Q,\gamma_1)\cap \ldots \cap \Sigma(Q,\gamma_r)\subseteq
\Sigma(Q,\alpha)\cap\{\sigma\in \Gamma^\star\mid \sigma(\gamma_1)=\cdots=\sigma(\gamma_r)=0\}.
$$
The converse inequality
$$
\Sigma(Q,\gamma_1)\cap \ldots \cap \Sigma(Q,\gamma_r)\supseteq
\Sigma(Q,\alpha)\cap \{\sigma\in \Gamma^\star\mid \sigma(\gamma_1)=\cdots=\sigma(\gamma_r)=0\}.
$$
also holds: If $\alpha$ is $\sigma$-semi-stable, and $\sigma(\gamma_1)=\cdots=\sigma(\gamma_r)=0$
then $\gamma_1,\dots,\gamma_r$ are $\sigma$-semi-stable.
We conclude that
$$
\Sigma(Q,\gamma_1)\cap \ldots \cap \Sigma(Q,\gamma_r)=
\Sigma(Q,\alpha)\cap\{\sigma\in \Gamma^\star\mid \sigma(\gamma_1)=\cdots=\sigma(\gamma_r)=0\}.
$$
Using the saturation property, we get
$$
\RR_+\Sigma(Q,\gamma_1)\cap\ldots\cap\RR_+ \Sigma(Q,\gamma_r)=
\RR_+\Sigma(Q,\alpha)\cap 
\{
\sigma\in (\RR^n)^\star\mid  \sigma(\gamma_1)=\cdots=\sigma(\gamma_r)=0\}
$$
By the refinement Theorem (Theorem~\ref{theorefinement})
there exists an exceptional sequence $\underline \varepsilon=
(\varepsilon_1,\dots,\varepsilon_s)$ which is a refinement of 
$\underline \gamma$,
i.e., there exists a sequence $0=b_0 < b_1<\ldots <b_{r-1}<n=b_r$ such that 
for each $j=1,\ldots ,r$ the dimension
vector $\gamma_j$ is a positive linear combination of 
$\varepsilon_{b_{j-1}+1},\ldots ,\varepsilon_{b_j}$.

We proceed by induction on $r$. 
Suppose that $r=1$. Then $\alpha$ is a  integer multiple of $\gamma_1$.
We have $\Sigma(Q,\alpha)=\Sigma(Q,\gamma_1)$.
The cone $\RR_+\Sigma(Q,\alpha)=\RR_+\Sigma(Q,\gamma_1)$
is given by one equality, $\sigma(\gamma_1)=0$ and
many inequalities of the form $\sigma(\beta)\leq 0$
for all $\beta\hookrightarrow \gamma_1$.
There exists a weight $\tau$ such that $\gamma_1$ is $\tau$-stable.
This means that $\tau(\beta)<0$ for all $\beta\hookrightarrow \gamma_1$
and $\beta\neq 0,\gamma_1$. If we view $\sigma$, $\RR_+\Sigma(Q,\alpha)$
inside the space
$$
\{\sigma\in (\RR^n)^\star\mid \sigma(\gamma_1)=0\}\cong \RR^{n-1}
$$
then $\RR_+\Sigma(Q,\alpha)$ contains an open neighborhood of $\tau$.
This shows that $\RR_+\Sigma(Q,\alpha)=\RR_+\Sigma(Q,\gamma_1)$
is a cone of dimension $n-1$. In particular $\RR_+\Sigma(Q,\alpha)$
is a face of dimension $n-1$.

Suppose that $r>1$. Let $V_i$ be the unique representation
corresponding to the dense orbit in $\Rep(Q,\varepsilon_i)$, for
$i=1,2,\dots,b_1$. The representations $V_1,\dots,V_{b_1}$
generate the full  subcategory ${\mathcal C}(V_1,\dots,V_{b_1})$
 of $\Rep_K(Q)$ which is closed under
extensions, direct sums, and taking kernels and cokernels.
This category is equivalent to the category of representations
of some quiver $Q'$ with $b_1$ vertices and without oriented cycles.
(See Lemma~\ref{experp} and Theorem~\ref{rightperp}). 
The right orthogonal category
$$(V_1,\dots,V_{b_1})^\perp=V_1^\perp\cap \cdots \cap V_{b_1}^{\perp}$$ 
is the category
of representations of a quiver $Q''$ with $n-b_1$ vertices and
without oriented cycles (Theorem~\ref{rightperp}).
Define $I':\NN^{Q_0'}\to \NN^{Q_0}$ and $I'':\NN^{Q_0''}\to \NN^{Q_0}$ as in
Theorem~\ref{theoembedding}.


By the induction hypothesis, we can find linearly independent
$$\tau'_1,\dots,\tau'_{b_1-1}\in \Sigma(Q',(I')^{-1}(\gamma_1)).$$
and linearly independent
$$
\tau''_1,\dots,\tau''_{n-b_1-r+1}\in \Sigma(Q'',(I'')^{-1}(\gamma_2))\cap \cdots\cap
\Sigma(Q'',(I'')^{-1}(\gamma_r)).
$$
Define $\sigma'_1,\dots,\sigma'_{b_1-1}\in \Gamma^\star$
such that $\sigma'_i\circ I'=\tau_i'$ 
and $\sigma'_i\circ I''=0$ for $i=1,2,\dots,b_1-1$.
Define $\sigma''_1,\dots,\sigma''_{n-b_1-r-1}\in \Gamma^\star$
by $\sigma''_i\circ I''=\tau''_i$ and $\sigma''_i\circ I=0$ for $i=1,2,\dots,n-b_1-r-1$.
Then
$$\sigma'_1,\dots,\sigma'_{b_1-1},\sigma''_1,\dots,\sigma''_{n-b_1-r+1}\in
\Sigma(Q,\gamma_1)\cap \cdots\cap \RR_+\Sigma(Q,\gamma_r)$$
are $n-r$ linearly independent weights. This shows that the face
$$\RR_+\Sigma(Q,\gamma_1)\cap \cdots\cap \RR_+\Sigma(Q,\gamma_r)$$ 
has
dimension $n-r$.

We have proved that the map $\psi(r)$ defined above is well defined. 
Let us show that the inverse map is well defined
and that it is indeed an inverse. Let ${F}$ be a face of dimension $n-r$ of 
$\RR_+\Sigma (Q,\alpha)$. Take a weight
$\sigma\in \Gamma^\star$ in the relative interior of ${F}$. 
Let 
$$\alpha =c_1\cdot \delta_1\pp c_2\cdot \delta_2\pp\cdot\pp c_l\cdot \delta_l
$$
be the $\sigma$-stable decomposition of 
$\alpha$. Define
${F}^\prime =\RR_+\Sigma(Q,\delta_1)\cap\cdots\cap \Sigma(Q,\delta_l)$. 
Since $\sigma\in {F}^\prime$ and $\sigma$ is
in the relative interior of $F$, we have ${F}\subseteq {F}^\prime$. 

Suppose that $\gamma\hookrightarrow\alpha$ and $\sigma(\gamma)=0$.
Then $\gamma$ is a linear combination
of the $\delta_i$'s by the definition of $\sigma$-stable decomposition. 
But the description of 
$\Sigma (Q,\alpha )$ given in Theorem~\ref{theosubrep} implies that
$$
F'=\RR_+\Sigma(Q,\alpha)\cap \{\sigma\in (\RR^n)^\star\mid \sigma(\delta_1)=\cdots=\sigma(\delta_l)=0\}\subseteq$$
$$
\subseteq
{F}= \RR_+\Sigma(Q,\alpha)\cap
\bigcap_{\gamma,\gamma\hookrightarrow\alpha,\langle\gamma,\beta\rangle=0}
\lbrace \sigma\in (\RR^n)^\star\mid \sigma(\gamma)=0  \rbrace
$$

so we have ${F}^\prime\subseteq {F}$. This concludes the proof of the Theorem.
\end{proof} 

Let us state the meaning of Theorem~\ref{bijection} in two extreme cases: for the walls of maximal dimension of $\Sigma
(Q,\alpha )$ and for extremal rays.

\begin{corollary}
Let $Q$ be a quiver without oriented cycles and let $\alpha$ be a Schur root. Then the
walls of $\Sigma(Q,\alpha)$ (i.e., faces of dimension $n-2$) 
are in one to one correspondence 
with the ways of writing
$$\alpha = c_1\gamma_1 +c_2 \gamma_2$$
where $(\gamma_1 ,\gamma_2 )$ is a quiver Schur sequence, 
$c_1 ,c_2$ positive integers with $c_i=1$ whenever $\gamma_i$ is imaginary
and non-isotropic.
\end{corollary}
 Let us consider an extremal ray $\sigma$ in $\Sigma(Q,\alpha)$. It corresponds to
the linear combination
$$\alpha = c_1\gamma_1 +\ldots +c_{n-1}\gamma_{n-1}$$ 
where $(\gamma_1 ,\ldots ,\gamma_{n-1})$ is a quiver Schur sequence. 
Theorem~\ref{theorefinement} implies that at least $n-2$ of the roots $\gamma_1 ,\ldots
,\gamma_{n-1}$ are real Schur roots.
Consider the subring 
$$\SI(Q,\alpha ,\sigma)=\bigoplus_{m\ge 0} \SI(Q,\alpha)_{m\sigma}$$
By peeling off real roots from the left and from the right we can reduce 
the calculation of this ring to the
ring of semi-invariants for a quiver $\theta(\ell)$
where $\theta(\ell)$ is the
Kronecker quiver $\theta(\ell)$ with
two vertices and $\ell$ equioriented arrows.
\begin{corollary} 
Let $Q$ be a quiver with no oriented cycles and let $\alpha$ be a Schur root.
Suppose that $\sigma\in \Sigma(Q,\alpha)$ is indivisible and
spans an extremal ray. Then there exists a positive integer $n$ and a 
Schur root
$\beta$ for the quiver $\theta(\ell)$ such that
$$
\SI(Q,\alpha,\sigma)\cong\SI(\theta(\ell),\beta).
$$
\end{corollary}

\begin{example}
Let $Q$ be the quiver
$$
\xymatrix{
& \circ\ar[d] & \\
\circ\ar[r] & \circ & \circ\ar[l]\\
& \circ\ar[u] &}
$$

and let $\alpha$ be the dimension vector
$$
\begin{matrix}
& 1 & \\
1 & 2 & 1 \\
& 1 & \\
\end{matrix}
$$
There are 8 walls of $\Sigma(Q,\alpha)$. They are given by
the Schur sequences
$$
\begin{matrix}
& 0 & \\
1 & 2 & 1 \\
& 1 & \\
\end{matrix}\ ,\ \ 
\begin{matrix}
& 1 & \\
0 & 0 & 0 \\
& 0 & \\
\end{matrix}
\qquad \mbox{(4 by symmetry)}
$$
$$
\begin{matrix}
& 1 & \\
0 & 1 & 0 \\
& 0 & \\
\end{matrix}\ ,\ \ 
\begin{matrix}
& 0 & \\
1 & 1 & 1 \\
& 1 & \\
\end{matrix}
\qquad \mbox{(4 by symmetry).}
$$
There are 12 two-dimensional faces of the cone given
by the sequences
$$
\begin{matrix}
& 1 & \\
0 & 1 & 0 \\
& 0 & \\
\end{matrix}\ ,\ \ 
\begin{matrix}
& 0 & \\
1 & 1 & 0 \\
& 1 & \\
\end{matrix}\ ,\ \ 
\begin{matrix}
& 0 & \\
0 & 0 & 1 \\
& 0 & \\
\end{matrix}
\qquad \mbox{(12 by symmetry).}
$$
There are 6 extremal rays, given by the Schur sequences
$$
\begin{matrix}
& 1 & \\
0 & 1 & 0 \\
& 0 & \\
\end{matrix}\ ,\ \
\begin{matrix}
& 0 & \\
1 & 1 & 0 \\
& 0 & \\
\end{matrix}\ ,\ \ 
\begin{matrix}
& 0 & \\
0 & 0 & 0 \\
& 1 & \\
\end{matrix}\ ,\ \ 
\begin{matrix}
& 0 & \\
0 & 0 & 1 \\
& 0 & \\
\end{matrix}
\qquad \mbox{(6 by symmetry).}
$$
The set $\Sigma(Q,\alpha)$ is a cone over a regular octahedron.\\
\centerline{\includegraphics[width=3in]{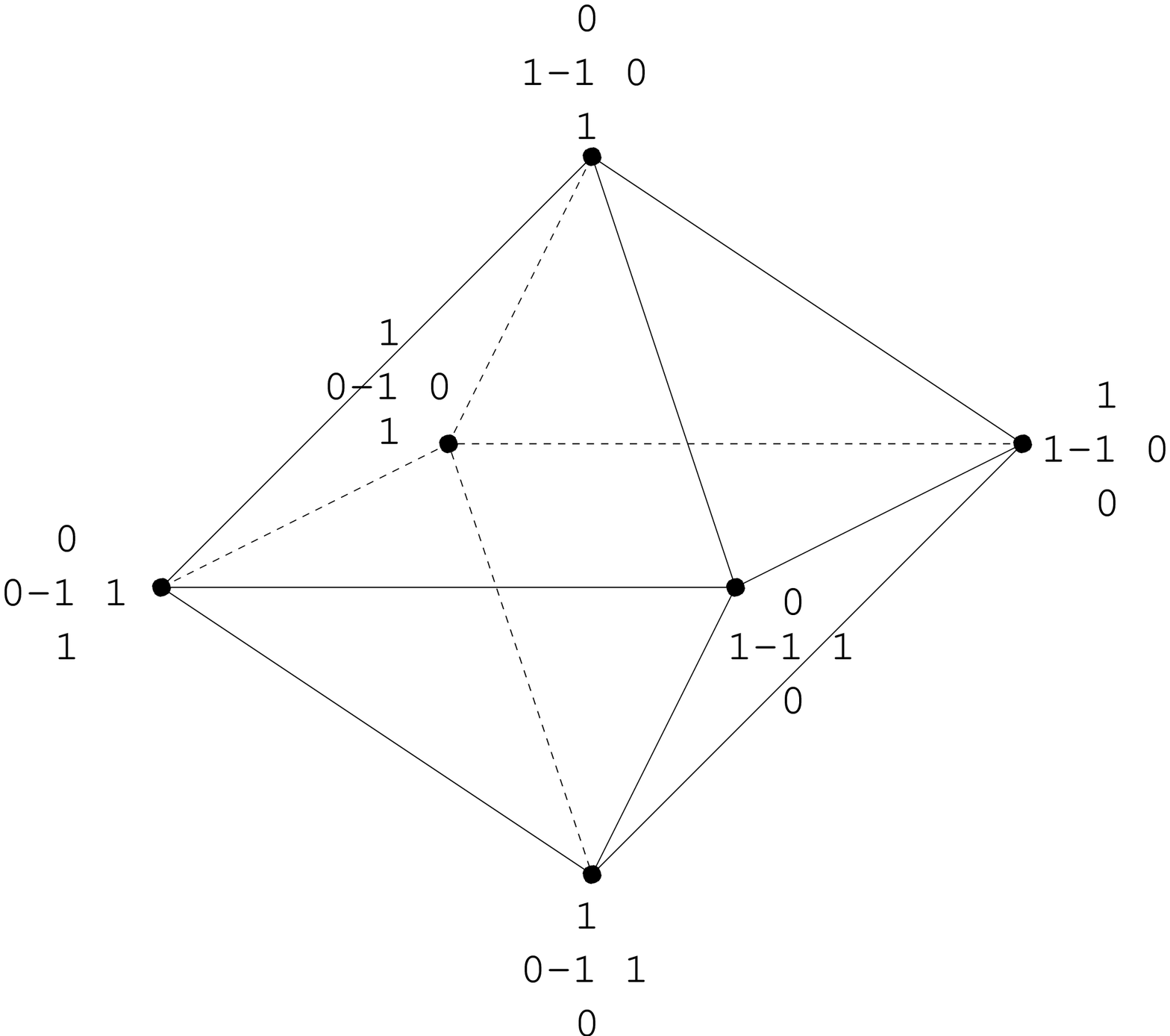}}\\[10pt]
\end{example}
\begin{example}
Let $Q$ be the quiver
$$
\xymatrix{
& \circ\ar@<.5ex>[d]\ar@<-.5ex>[d] & \\
\circ\ar@<.5ex>[r]\ar@<-.5ex>[r] & \circ & \circ\ar@<.5ex>[l]\ar@<-.5ex>[l]}
$$

and let $\alpha$ be the dimension vector
$$
\begin{matrix}
& 1 & \\
1 & 3 & 1\end{matrix}.
$$
Now $\Sigma(Q,\alpha)$ has 6 walls corresponding to the Schur sequences
$$
\begin{matrix}
& 0 & \\
1 & 3 & 1
\end{matrix},\ \ 
\begin{matrix}
& 1 & \\
0 & 0 & 0
\end{matrix}
\quad\mbox{(3 by symmetry)}
$$
$$
\begin{matrix}
& 1 & \\
0 & 2 & 0
\end{matrix},\ \ 
\begin{matrix}
& 0 & \\
1 & 1 & 1
\end{matrix}
\quad\mbox{(3 by symmetry)}
$$
The are also 6 extremal rays which correspond to the Schur sequences
$$
\begin{matrix}
& 1 & \\
0 & 2 & 0
\end{matrix},\ \ 
\begin{matrix}
& 0 & \\
1 & 1 & 0
\end{matrix},\ \ 
\begin{matrix}
& 0 & \\
0 & 0 & 1
\end{matrix}\quad\mbox{(6 by symmetry)}.
$$
The cone $\Sigma(Q,\alpha)$ is a cone over a hexagon.\\
\centerline{\includegraphics[width=2in]{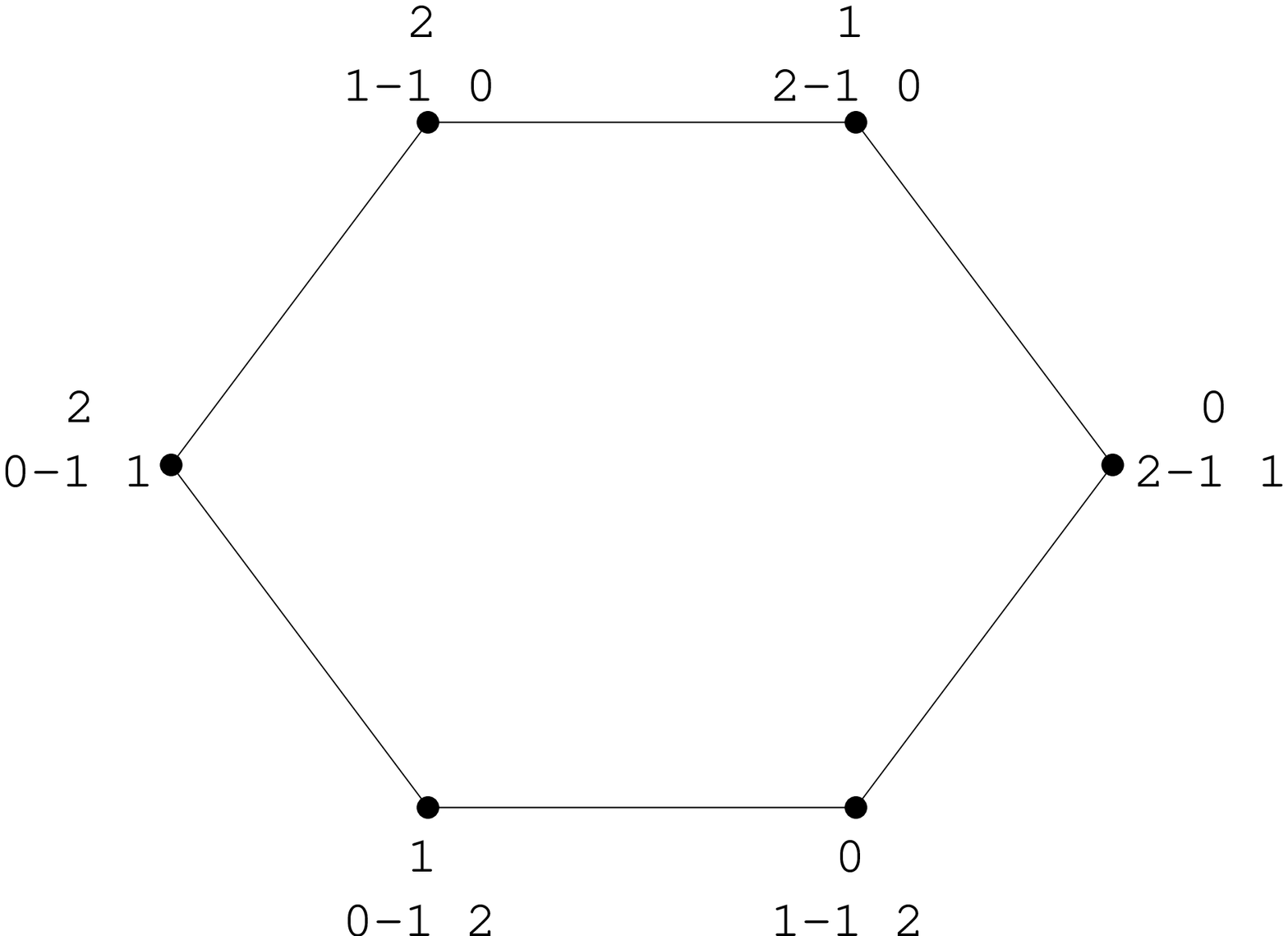}}\\[10pt]
\end{example}
\section{More on the $\sigma$-stable decompositions}\label{sec6}
\subsection{The set of $\sigma$-stable dimension vectors}
In the previous section we studied, how the $\sigma$-stable decomposition
of $\alpha$ varies, when $\sigma$ varies and $\alpha$ is fixed. This
led to the description of the faces of $\Sigma(Q,\alpha)$. In this
section, we will study how the $\sigma$-stable decomposition looks like
for a fixed weight $\sigma$. Let us define $\overline{\Sigma}(Q,\sigma)$
as the set of all $\sigma$-semi-stable dimension vectors.
Notice that
$$
\alpha\in \overline{\Sigma}(Q,\sigma)\quad\Leftrightarrow\quad
\sigma\in \Sigma(Q,\alpha).
$$

\begin{lemma}\label{extreme}
Suppose that $\alpha\in \overline{\Sigma}(Q,\sigma)$.
There exists a sequence $\underline{\delta}=\delta_1,\dots,\delta_s$ of dimension vectors
such that
\begin{enumerate}
\renewcommand{\theenumi}{\arabic{enumi}}
\item $\alpha=\sum_{i=1}^s a_i\delta_i$ for some positive rational numbers 
$a_1,\dots,a_s$;
\item $\delta_1,\dots,\delta_s$ are linearly independent dimension vectors;
\item each $\delta_i$ generates an extremal ray, i.e., 
there exists an extremal ray ${\mathcal R}_i$ of $\RR_+\overline{\Sigma}(Q,\sigma)$ 
such that ${\mathcal R}_i\cap \overline{\Sigma}(Q,\sigma)=\NN\delta_i$.
\end{enumerate}
\end{lemma}
\begin{proof} This is trivial.\end{proof}
\begin{lemma}\label{ratcomb}
Suppose that $\alpha,\beta,\delta_1,\dots,\delta_s$ are $\sigma$-stable, 
and $\beta\hookrightarrow \alpha$.
\begin{enumerate}
\renewcommand{\theenumi}{\alph{enumi}}
\item If $\alpha$ is a nonnegative integral combination of $\delta_1,\dots,\delta_s$
then so is $\beta$.
\item If $\alpha$ is a nonnegative rational  combination of $\delta_1,\dots,\delta_s$
then so is $\beta$.
\end{enumerate}
\end{lemma}
\begin{proof}
Suppose that $\alpha=\sum_{i=1}^s a_i\delta_i$ for some integers $a_1,\dots,a_s$.
Note that $\delta_i$ is $\sigma$-stable for all $i$ because $\delta_i$ is not
the sum of two smaller $\sigma$-semi-stable dimension vectors.
Let $V_i$ be $\sigma$-stable of dimension $\delta_i$ for all $i$.
Consider the representation 
$$V=V_1^{a_1}\oplus V_2^{a_2}\oplus \cdots \oplus V_s^{a_s}.$$
Now $V$ has a semi-stable subrepresentation $W$ of dimension $\beta$
by Lemma~\ref{semisub}.
In a Jordan-H\"older filtration of $W$, only $V_1,\dots,V_s$ will
appear,  so $\beta$ must
be a nonnegative integral combination of $\delta_1,\dots,\delta_s$.

The second statement follows from the fact that for each positive integer $m$ we
have (see Theorem~\ref{lemScho})
$$
\beta\hookrightarrow\alpha\quad\Rightarrow \quad m\beta\hookrightarrow m\alpha.
$$
\end{proof}
\begin{definition}
For a sequence of roots $\underline{\alpha}=(\alpha_1,\dots,\alpha_s)$
(all $\alpha_i$ distinct)
with $\langle \alpha_i,\alpha_j\rangle \leq 0$ for $i\neq j$,
we define a quiver $Q(\underline{\alpha})$ as follows.
The set of vertices of $Q(\underline{\alpha})_0$ is equal to $\{1,2,\dots,s\}$.
For $i\neq j$ there are $-\langle \alpha_i,\alpha_j\rangle$
arrows from $i$ to $j$. 
There are $1-\langle \alpha_i,\alpha_i\rangle$
arrows (loops) from $i$ to $i$.
\end{definition}
\begin{theorem}\label{sigmastable}
Under the assumptions of Lemma~\ref{extreme} above,
$\alpha$ is $\sigma$-stable if and only if
\begin{enumerate}
\renewcommand{\theenumi}{\arabic{enumi}}
\item either $\alpha=\delta_i$ and $\delta_i$ is a real Schur root for some $i$,
\item or $\langle \delta_i,\alpha\rangle\leq 0$ and
$\langle \alpha,\delta_i\rangle \leq 0$ for all $i$,
$Q(\underline{\delta})$ is path connected
and $\alpha$ is indivisible if $\alpha$ is isotropic.
\end{enumerate}
\end{theorem}
\begin{proof}
First we prove that the conditions are necessary.
Assume that $\alpha$ is $\sigma$-stable. For $i\neq j$ we have $\hom_Q(\delta_j,\delta_i)=0$
because $\delta_i,\delta_j$ are $\sigma$-stable (see Lemma~\ref{lemmorphism}).
Suppose that $\langle \alpha,\delta_i\rangle>0$. This is only possible
when $\delta_i$ is a real Schur root because $\langle\delta_j,\delta_i\rangle\leq 0$
for all $j\neq i$. In that case we have
$\hom_Q(\alpha,\delta_i)\neq 0$, so $\alpha=\delta_i$ because $\alpha$
and $\delta_i$ are $\sigma$-stable (see Lemma~\ref{lemmorphism}).
Similarly, if  $\langle \delta_i,\alpha\rangle>0$ then we have $\alpha=\delta_i$.

Consider the quiver $Q(\underline{\delta})$. 
Let $S_1$ be the set of all $k$ such that there is a path from
$1$ to $k$ and let $S_2=S\setminus S_1$.
Define 
$$
\alpha_1=\sum_{i\in S_1}a_i\delta_i,
\qquad
 \alpha_2=\sum_{i\in S_2}a_i\delta_i.
$$
There are no arrows from $S_1$ to $S_2$.
This shows that $\langle \alpha_1,\alpha_2\rangle=0$.
Choose an integer $m$ such that $ma_i$ is a positive integer for all $i$.
We have $m\alpha_1\hookrightarrow m\alpha$.
Since $\alpha$ is $\sigma$-stable, the $\sigma$-stable decomposition
of $m\alpha$ is either $m\alpha$ or $\alpha\pp\dots\pp\alpha$ by Proposition~\ref{multiplestable}.
Because  $m\alpha_1$ is $\sigma$-semi-stable, 
$m\alpha_1$ must be proportional
to $\alpha$. This can only happen if $S_2=\emptyset$.
We have shown that the quiver  $Q(\underline{\delta})$
is path connected.

If $\alpha$ is isotropic, then $\alpha$ must be indivisible because of Proposition~\ref{multiples}.

Clearly, if condition (1) is satisfied, then $\alpha$ is $\sigma$-stable.
Suppose that (2) is satisfied.
Suppose that $\beta\hookrightarrow\alpha$,  $\beta$ is $\sigma$-stable
and $\beta\neq 0,\alpha$. We will derive a contradiction.
By Lemma~\ref{ratcomb}, $\beta=\sum_{i=1}^s b_i \delta_i$ such that
the $b_i$'s are nonnegative rational numbers.
Define
$$
T_1=\supp(\beta):=\{i\mid 1\leq i\leq s,\ b_i\neq 0\}
$$
Let 
$$T_2=\{1,2,\dots,s\}\setminus T_1=\{i\mid 1\leq i\leq s,\ b_i=0\}.
$$
Assume that $T_2\neq\emptyset$.
Define
$$
\alpha_1=\sum_{i\in T_1}a_i\delta_i,\quad
\alpha_2=\sum_{i\in T_2}a_i\delta_i.
$$
Now $\alpha=\alpha_1+\alpha_2$. Because $\langle\delta_i,\delta_j\rangle\leq 0$
for all $i\neq j$, and because there must be at least one arrow from $T_1$
to $T_2$, we get $\langle \beta,\alpha_2\rangle<0$.
Since $\ext_Q(\beta,\alpha-\beta)=0$, we get $\langle \beta,\alpha-\beta\rangle\geq 0$,
so $\langle\beta,\gamma\rangle>0$ with 
$$
\gamma=(\alpha-\beta)-\alpha_2=\alpha_1-\beta=\sum_{\alpha_i\in T_1}(a_i-b_i)\delta_i.
$$
If
$$
\gamma=\gamma_1\pp\gamma_2\pp\cdots\pp\gamma_r
$$
is the $\sigma$-stable decomposition of $\gamma$, then
$\langle\beta,\gamma_j\rangle>0$ for some $j$. We get
that $\beta=\gamma_j$ by Lemma~\ref{lemmorphism}. But then
$$\langle\beta,\gamma_j\rangle=\langle\beta,\beta\rangle=\langle\beta,\alpha\rangle-
\langle\beta,\alpha-\beta\rangle\leq 0.$$
Contradiction.

This shows that $T_2=\emptyset$ and $T_1=\supp(\beta)=\{1,2,\dots,s\}$. Let $\gamma=\alpha-\beta$.
We have $\gamma\neq 0,\alpha$.
We can find a $\sigma$-stable dimension vector $\gamma'$ such
that $\gamma\twoheadrightarrow \gamma'$ and therefore $\alpha\twoheadrightarrow
\gamma'$. By a similar argument as before we obtain
$\supp(\gamma')=\supp(\gamma)=\{1,2,\dots,s\}$.
Write $\gamma=\sum_{i=1}^s c_i\delta_i$ with $c_i$ a positive rational
number for all $i$.
We have
\begin{equation}\label{ci}
0=\langle\beta,\gamma\rangle=\sum_{i=1}^s c_i\langle \beta,\delta_i\rangle.
\end{equation}
Since $\beta$ and $\delta_i$ are $\sigma$-stable, and $\beta\neq \delta_i$
we get that $\hom_Q(\beta,\delta_i)=0$ by Lemma~\ref{lemmorphism}.
In particular, $\langle\beta,\delta_i\rangle\leq 0$ for all $i$.
Combined with \ref{ci} 
we conclude that $\langle\beta,\delta_i\rangle=0$ for all $i$.

Let $B=\max\{b_1,\dots,b_s\}$. Let $S_1=\{i\mid b_i=B\}$ and
let $S_2=\{1,2,\dots,s\}\setminus S_1$.
Suppose $S_2\neq \emptyset$. There must be an arrow from $S_1$ to $S_2$,
say $j\to k$ with $j\in S_1$ and $k\in S_2$.
Then
$$
0=\langle \beta,\delta_k\rangle\leq
b_j\langle\delta_j,\delta_k\rangle+b_k\langle\delta_k,\delta_k\rangle\
$$
We know that $b_j>b_k$, $\langle\delta_j,\delta_k\rangle\leq -1$, and
$\langle\delta_k,\delta_k\rangle\leq 1$. This leads to a contradiction.

So $S_2=\emptyset$ and $b_i=B$ for all $i$.
From $\langle \beta,\delta_i\rangle=0$ follows that for every
$i$ there is exactly one arrow with tail $i$.
Since $Q(\underline{\delta})$ is path connected, the quiver 
$Q(\underline{\delta})$ has to be a cycle.
Now it easily follows from $\langle\alpha,\delta_i\rangle\leq 0$
that $\alpha$ must be proportional to $\beta$ and $\langle \alpha,\alpha\rangle=0$.
Since $\alpha$ is indivisible in that case, $\alpha=\beta$.
\end{proof}
\begin{remark}
The previous theorem also provides us an inductive way of finding the $\sigma$-stable
decomposition, if $\alpha$ is $\sigma$-semi-stable, but not $\sigma$-stable.
There are two cases.

In the first case 
$\langle \delta_i,\alpha\rangle>0$ or $\langle \alpha,\delta_i\rangle>0$
for some extremal dimension vector $\delta_i\in\overline{\Sigma}(Q,\sigma)$.
If we know the $\sigma$-stable decomposition of
the smaller dimension vector $\alpha-\delta_i$,
we know the $\sigma$-stable decomposition of $\alpha$.

In the second case
$Q(\underline{\delta})$ is not path-connected, say there is no path
from $i$ to $j$. Let $S_1$ be the set of all $k$ such that
there is a path from $i$ to $k$ and let $S_2$ be the complement.
There are no arrows from $S_1$ to $S_2$.
If we define 
$$\alpha_1=\sum_{i\in S_1} a_i\delta_i,\qquad
\alpha_2=\sum_{i\in S_2} a_i\delta_i.
$$
For some $m$, $m\alpha_1,m\alpha_2$ are dimension vectors, and
$m\alpha_1\hookrightarrow m\alpha$.
If we know the $\sigma$-stable decomposition of $m\alpha_1$ and $m\alpha_2$
then we know the $\sigma$-stable decomposition of $m\alpha$ and $\alpha$. 
Now $m\alpha_1$ and $m\alpha_2$ have smaller support than $\alpha$.
\end{remark}
\begin{example}
Let $Q$ be the quiver
$$
\xymatrix{
\circ\ar[r]\ar[rd]\ar[d] & \circ\ar[d]\\
\circ \ar[r]\ar[ru] & \circ}
$$

and let $\sigma$ be the weight
$$
\begin{matrix}
1 & -1\\
1 & -1
\end{matrix}
$$
The extremal rays of the cone $\overline{\Sigma}(Q,\sigma)$ are given
by the dimension vectors
$$
\delta_1=\begin{matrix}
1 & 1\\
0 & 0
\end{matrix},\ \ 
\delta_2=\begin{matrix}
1 & 0\\
0 & 1
\end{matrix},\ \ 
\delta_3=\begin{matrix}
0 & 0\\
1 & 1
\end{matrix},\ \ 
\delta_4=\begin{matrix}
0 & 1\\
1 & 0
\end{matrix}. 
$$
The quiver $Q(\underline{\delta})$ is
$$
\xymatrix{
1\ar[rr]\ar[dd]\ar@<1.5ex>[rrdd]\ar@<.5ex>[rrdd]\ar@<-.5ex>[rrdd]
& & 2\ar[dd]\ar@<1.5ex>[lldd]\ar@<.5ex>[lldd]\\
& & \\
4\ar@<.5ex>[rruu]\ar@<1.5ex>[rruu]\ar[rr] & & 3\ar@<1.5ex>[lluu]}
$$

Any $\sigma$-stable dimension vector is a nonnegative rational
combination of $\delta_1,\delta_2,\delta_4$ or a
nonnegative rational combination of $\delta_2,\delta_3,\delta_4$.
Suppose $\alpha$ is $\sigma$-stable and not
equal to $\delta_1,\delta_2,\delta_3,\delta_4$.
If $\alpha$ is a nonnegative rational combination
$\delta_1,\delta_2,\delta_4$, then because
the support has to be path connected, it
is actually a nonnegative rational combination
of $\delta_2,\delta_4$.
Similarly, if $\alpha$ is a nonnegative rational
combination of $\delta_2,\delta_3,\delta_4$, then
it must be in fact a nonnegative rational combination
of $\delta_2$ and $\delta_4$.

Now it easily follows that the set of $\sigma$-stable
dimension vectors is
$$
\delta_1,\ \ \delta_2,\ \ \delta_3,\ \ \delta_4,\ \ a\delta_2+b\delta_4\  
(a,b>0,a,b\in \ZZ,a\leq 2b,b\leq 2a).
$$
The cone $\overline{\Sigma}(Q,\alpha)$ is a cone over a square.
In the diagram below, the coordinates should be interpreted as
projective coordinates.\\
\centerline{\includegraphics[width=1.5in]{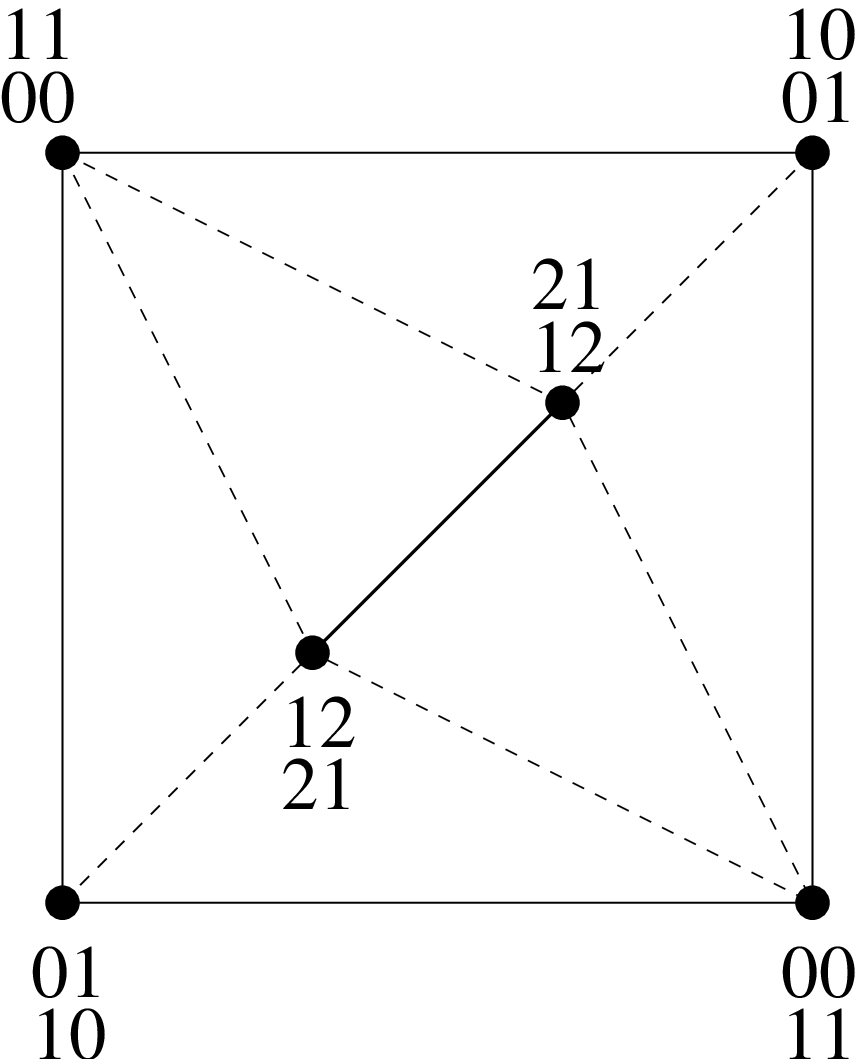}}\\[10pt]
The fat line in the middle of the square corresponds to the
imaginary $\sigma$-stable dimension vectors. The dashed lines distinguish
the regions where the $\sigma$-stable decomposition looks different.

Some examples of $\sigma$ stable decomposition are:
$$
\begin{matrix}
5 & 4\\
3 & 4
\end{matrix}
=\begin{matrix}
1 & 1\\
0 & 0
\end{matrix}\ \pp\ 
\begin{matrix}
4 & 3\\
3 & 4
\end{matrix}
$$
$$
\begin{matrix}
1 & 4\\
5 & 2
\end{matrix}
=\begin{matrix}
0 & 0\\
1 & 1
\end{matrix}\ \pp\ 
\begin{matrix}
1 & 2\\
2 & 1
\end{matrix}\ \pp\ 
2\cdot
\begin{matrix}
0 & 1\\
1 & 0
\end{matrix}.
$$

\end{example}

\subsection{$\sigma$-stable decomposition for quivers with oriented cycles}
Doubling of the quiver, reduces the $\sigma$-stable decomposition
for quivers with oriented cycles, to the case of quivers without
oriented cycles. Suppose that $Q$ is a quiver with oriented cycles.
We define a new quiver $\widehat{Q}$ by $\widehat{Q}_0=Q_0\times \{0,1\}$.
For every $a\in Q_1$ we define an arrow $\widehat{a}\in \widehat{Q}_1$ with
$t\widehat{a}=(ta,0)$ and $h\widehat{a}=(ha,1)$ and
for every $x\in Q_0$ we define an arrow $\widehat{x}\in \widehat{Q}_1$ with
$t\widehat{x}=(x,0)$ and $h\widehat{x}=(hx,1)$.
For example, if $Q$ is the quiver 
$$
\xymatrix{
\circ \ar@<.5ex>[r] & \circ \ar@<.5ex>[r]\ar@<-.5ex>[r]\ar@<.5ex>[l] & \circ}
$$
then $\widehat{Q}$ is the quiver
$$
\xymatrix{
\circ\ar[rd]\ar[d] & \circ \ar@<.5ex>[rd]\ar@<-.5ex>[rd]\ar[ld]
\ar[d] & \circ\ar[d]\\
\circ & \circ & \circ}
$$

For a $Q$-dimension vector $\alpha$, we define a dimension vector $\widehat{\alpha}$
of $\widehat{Q}$ by $\widehat{\alpha}(x,0)=\widehat{\alpha}(x,1)=\alpha(x)$
for all $x\in Q_0$.
Similarly, if $\sigma$ is a weight of $Q$, we define a weight $\widehat{\sigma}$ 
of $\widehat{Q}$ by
$\widehat{\sigma}(x,0)=\widehat{\sigma}(x,1)=\sigma(x)$.
We define the weight $\tau$ of $\widehat{Q}$ by $\tau(x,0)=1$ and
$\tau(x,1)=-1$ for all $x\in Q_0$.
Note that for any $\alpha\in Q_0$, $\widehat{\alpha}$ is $\tau$-stable.
\begin{proposition}
Suppose that $\alpha$ is a dimension vector and $\sigma$ is a weight for $Q$. Then
$\alpha$ is $\sigma$-semi-stable (stable) if and only if for some large positive integer
$m$, $\widehat{\alpha}$ is $\widehat{\sigma}+m\tau$-semi-stable (stable).
\end{proposition}
\begin{proof}
Suppose that $\alpha$ is $\sigma$-semi-stable. 
If $\gamma\hookrightarrow \widehat{\alpha}$
for some $\widehat{Q}$-dimension vector $\gamma$.
Note that $\gamma(x,0)\leq \gamma(x,1)$ for all $x\in Q_0$ because
$\widehat{\alpha}(x,0)=\widehat{\alpha}(x,1)$, and for a general representation
$V$ of dimension $\widehat{\alpha}$ the map
$V(\widehat{x}):V(x,0)\to V(x,1)$ is injective.
If $\gamma(x,0)=\gamma(x,1)$ for all $x\in Q_0$ then $\gamma$ is of the form
$\widehat{\beta}$ and $\beta\hookrightarrow \alpha$. 
Then we have $\widehat{\sigma}(\gamma)=\sigma(\beta)\leq 0$.
Also we have $(\widehat{\sigma}+m\tau)(\gamma)=\widehat{\sigma}(\gamma)\leq 0$.

Suppose that $\gamma(x,0)< \gamma(x,1)$ for some $x\in Q_0$.
Then $\tau(\gamma)<0$ so in particular for $m$ large enough we will
have $(\widehat{\sigma}+m\tau)(\gamma)<0$.

Since there are only finitely many subdimension vectors $\gamma$, we can choose
$m$ large enough such that $(\widehat{\sigma}+m\tau)(\gamma)\leq 0$ 
for all $\gamma\hookrightarrow \widehat{\alpha}$.
This shows that $\widehat{\alpha}$ is $(\widehat{\sigma}+m\tau)$-semi-stable.

Conversely, assume that $\widehat{\alpha}$ is $(\widehat{\sigma}+m\tau)$-semi-stable
for some $m$
and $\beta\hookrightarrow\alpha$. Then
$\widehat{\beta}\hookrightarrow\widehat{\alpha}$, so $0\geq (\widehat{\sigma}+m\tau)
(\widehat{\beta})=\widehat{\sigma}(\widehat{\beta})=\sigma(\beta)$.
This shows that $\alpha$ is $\sigma$-semi-stable.

A similar statement with stable instead of semi-stable is easy to prove.
\end{proof}
Suppose now that $Q$ is quiver, possibly with oriented cycles. Let us consider
the $0$-stable decomposition. Clearly, every representation of $Q$ is $0$-semi-stable
in the sense of Theorem~\ref{KingCrit}.
A representation $V$ is $0$-stable if the only subrepresentations
are $0$ and $V$ itself. In other words, $0$-stable representations are exactly
simple representations. Notice that if there exists an $\alpha$-dimensional
simple representation, then the general representation of dimension $\alpha$
is simple. We will call such dimension vectors simple.
\begin{corollary}\label{simpledec}
Suppose that $Q$ is an arbitrary quiver. For each $x\in Q_0$ we define
a dimension vector $\delta_x$ by $\delta_x(y)=0$ for $y\neq x$ and
$\delta_x(x)=1$.
 A dimension vector
$\alpha$ is simple if
\begin{enumerate}
\renewcommand{\theenumi}{\arabic{enumi}}
\item either $\alpha=\delta_x$ and $\delta_x$ is real (i.e., $\langle
\delta_x,\delta_x\rangle=1$);
\item or $\langle\delta_x,\alpha\rangle\leq 0$ and $\langle\alpha,\delta_x\rangle\leq
0$ for all $x\in Q_0$, the full subquiver of $Q$ with vertices
$$
\supp(\alpha):=\{x\in Q_0\mid \alpha(x)\neq 0\}
$$
is path connected, and if $\alpha$ is isotropic, then $\alpha$ is indivisible.
\end{enumerate}
\end{corollary}
\begin{remark}
A statement similar to Corollary~\ref{simpledec}
was proved in \cite{LP}. Theorem~\ref{sigmastable} is
a generalization of this result.
\end{remark}
\begin{example}
Consider the quiver
$$
\xymatrix{
& 1\ar@<.5ex>[ld]\ar[rd]\ar@<1ex>[rd] & \\
\circ\ar@<.5ex>[ru]\ar[rr] & & 3\ar@<1ex>[lu]}
$$
Suppose $\alpha$ that is the dimension vector $(a_1,a_2,a_3)$.
We will find necessary and sufficient conditions for $\alpha$ to 
be a simple dimension vector.
Of course $\alpha$ can be equal to $\delta_1$, $\delta_2$, $\delta_3$.
The conditions $\langle\delta_i,\alpha\rangle\leq 0$ and
$\langle\alpha,\delta_i\rangle\leq 0$ give the inequalities
$$
a_1\leq a_2+a_3,\ a_2\leq a_1,\ a_3\leq a_1+a_2
$$
(other inequalities turn out to be redundant).
If $a_3=0$ and $a_1=a_2$, then $\alpha$ is isotropic, so
we must have that $a_1=a_2=1$ in that case.
The only support of $\alpha$ which is not possible 
(because it is not path connected) is $\{2,3\}$,
but this is already excluded by the inequalities.

From the inequalities and the remarks above it
is easy to deduce that set of simple dimension vectors is given by
$$
(1,0,0),\ (0,1,0),\ (0,0,1),\ (1,0,0),\ (1,1,0),\ \mbox{and all}
$$
$$
\{(a_1,a_2,a_3)\in \ZZ^3\mid a_1\leq a_2+a_3,\ a_2\leq a_1,\ a_3\leq a_1+a_2,a_3>0\}
$$

In the picture shows how the simple decomposition looks like.
We use projective coordinates.\\
\centerline{\includegraphics[width=3in]{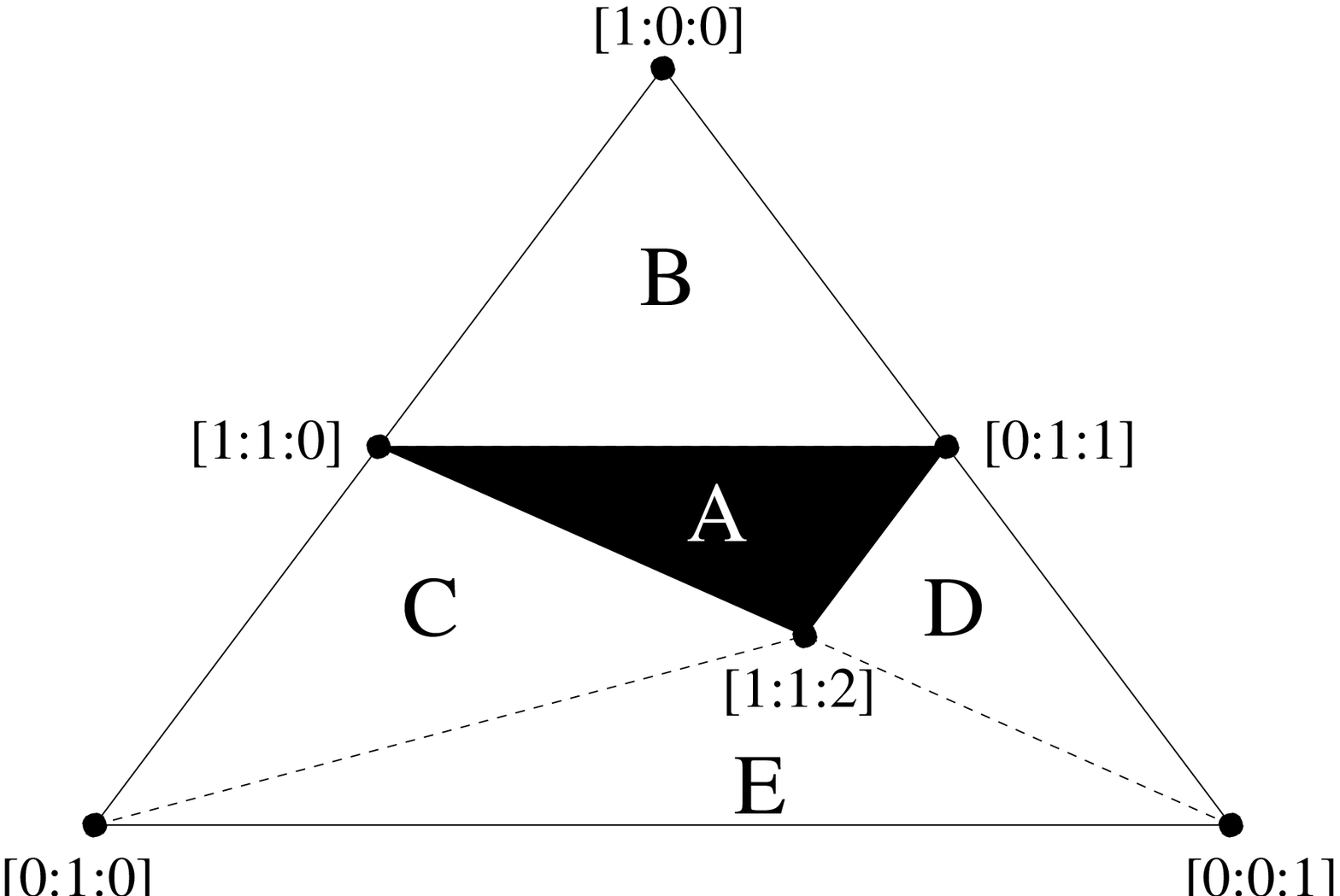}}\\[10pt]

Region A is defined by $a_1\leq a_2+a_3,\ a_2\leq a_1,\ a_3\leq a_1+a_2$. 
This will always define a simple dimension vector except when $a_3=0$
(and $a_1=a_2$). In that case, the simple decomposition is
$$
(a,a,0)=a\cdot (1,1,0).
$$

Region B is defined by $a_2+a_3\leq a_1$. The simple decomposition in this
region is
$$
(a_1,a_2,a_3)=(a_1-a_2-a_3)\cdot (1,0,0)\pp (a_2+a_3,a_2,a_3) \qquad\mbox{if $c_3>0$
and}
$$
$$
(a_1,a_2,0)=(a_1-a_2)\cdot (1,0,0)\pp a_2\cdot (1,1,0).
$$

Region C is defined by $a_2\geq a_1$, $2a_1\geq a_3$. The simple decomposition
is
$$
(a_1,a_2,a_3)=(a_2-a_1)\cdot (0,1,0)\pp (a_1,a_1,a_3)\qquad\mbox{if $c_3>0$ and}
$$
$$
(a_1,a_2,0)=(a_2-a_1)\cdot (0,1,0)\pp a_1\cdot (1,1,0).
$$

Region D is defined by $a_1\geq a_2$ and $a_3\geq a_1+a_2$.
The simple decomposition here is
$$
(a_1,a_2,a_3)=(a_3-a_1-a_2)\cdot (0,0,1)\pp (a_1,a_2,a_1+a_2).
$$

Region E is defined by $a_2\geq a_1$, $a_3\geq 2a_1$.
The simple decomposition in this region is
$$
(a_1,a_2,a_3)=(a_2-a_1)\cdot (1,0,0)\pp (a_3-2a_1)\cdot (0,1,0)\pp (a_1,a_1,2a_1).
$$ 
\end{example}
\begin{example}
Let $Q$ be the quiver with 3 vertices (labeled 1, 2 and 3),
with a loop at each vertex and with arrows $1\to 2$, $2\to 3$ and
$3\to 1$. The set of simple dimension vectors is
$$
(1,0,0),\ (0,1,0),\ (0,0,1),\ \mbox{ and all }
(a,b,c)\mbox{ with $a,b,c>0$}.
$$
Notice that for example a dimension vector of the form $(a,b,0)$ ($a,b>0$) is not
simple because its support is not path connected.
\end{example}
\section{Littlewood-Richardson coefficients}\label{sec7}
\subsection{The Klyachko cone}
Polynomial representations of $\GL_n$ are parametrized by
 non-increasing integer sequences of length $n$.
 If $\lambda=(\lambda_1,\dots,\lambda_n)$ 
is such a sequence, then we denote the corresponding representation by 
$V_\lambda$.
We define
$$
|\lambda|:=\lambda_1+\cdots+\lambda_n.
$$
The {\it Littlewood-Richardson coefficient \/} $c_{\lambda,\mu}^\nu$
is defined by
$$
\dim(V_\lambda\otimes V_\mu\otimes V_\nu^\star)^{\GL_n}.
$$
We would like to study the set
$$
{\mathcal K}_n=\{(\lambda,\mu,\nu)\in (\ZZ^n)^3\mid \mbox{$\lambda,\mu,\nu$ 
are non-increasing and }c_{\lambda,\mu}^\nu\neq 0\}. 
$$
The cone $\RR_+{\mathcal K}_n$ is the {\it Klyachko cone}. 
The Klyachko cone has dimension $3n-1$.
Let $T_{p,q,r}$ be the quiver with $p+q+r-2$ vertices:
$$
\xymatrix{
x_1\ar[r] &  x_2\ar[r] & \cdots\ar[r] & x_{p-2}\ar[r]  
& x_{p-1}\ar[rd]  &    \\
  y_1\ar[r]& y_2\ar[r] & \cdots\ar[r] & y_{q-2}\ar[r]  & 
y_{q-1}\ar[r] & x_{p} \\
z_1\ar[r] & z_2\ar[r] & \cdots\ar[r] & z_{r-2}\ar[r]  
& z_{r-1}\ar[ru] &      } 
$$

We use the convention $y_{q}=z_{r}=x_p$.
In \cite{DW} we have seen that if we take
the dimension vector 
$$
\beta=
\begin{matrix}
1 & 2 & \cdots & n-1 & \\
1 & 2 & \cdots & n-1 & n \\
1 & 2 & \cdots & n-1 & \end{matrix},
$$
for $T_{n,n,n}$, 
then we can view $\dim \SI(Q,\beta)_{\sigma}$ as a Littlewood-Richardson
coefficient as follows.
If $\sigma$ is given by 
$$
\sigma=
\begin{matrix}
a_1 & a_2 & \cdots & a_{n-1} & \\
b_1 & b_2 & \cdots & b_{n-1} & c_n \\
c_1 & c_2 & \cdots & c_{n-1} & \end{matrix},
$$
then
$$
\dim \SI(Q,\beta)_{\sigma}=c_{\lambda,\mu}^{\nu}
$$
where
\begin{eqnarray*}
\lambda &= & \lambda(\sigma) = (a_1+\cdots+a_{n-1},a_2+\cdots+a_{n-1},\cdots,a_{n-1},0),\\
\mu  &= & \mu(\sigma)=(b_1+\cdots+b_{n-1},b_2+\cdots+b_{n-1},\cdots,b_{n-1},0),\\
\nu &=& \nu(\sigma)=(-c_n,-(c_n+c_{n-1}),\dots,-(c_{n}+c_{n-1}+\dots+c_1)).
\end{eqnarray*}
Conversely, if $\lambda,\mu,\nu\in \ZZ^n$, then 
$$
c_{\lambda,\mu}^\nu=\dim \SI(Q,\beta)_{\sigma}
$$
where
$$
\sigma=\sigma(\lambda,\mu,\nu)=
\begin{matrix}
\lambda_1-\lambda_2 & \lambda_2-\lambda_3, & \cdots & \lambda_{n-1}-\lambda_n & \\
\mu_1-\mu_2 & \mu_2-\mu_3 & \cdots & \mu_{n-1}-\mu_n & \lambda_n+\mu_n-\nu_1\\
\nu_{n-1}-\nu_n & \nu_{n-2}-\nu_{n-1} & \cdots & \nu_{1}-\nu_2 & 
\end{matrix}
$$

The set $\Sigma(Q,\beta)$ is almost equal to ${\mathcal K}_n$. In fact, we define
a bijection 
$$
\psi:\Sigma(Q,\beta)\times \ZZ^2\to {\mathcal K}_n
$$
by
$$
\psi(\sigma,a,b)=(\lambda(\sigma)+a\cdot {\bf 1},\mu(\sigma)+b\cdot {\bf 1},\nu(\sigma)+(a+b)\cdot {\bf 1}).
$$
where ${\bf 1}=(1,1,\dots,1)\in \NN^n$. This bijection extends to an isomorphism
of the cones $\RR_+\Sigma(Q,\beta)\times \RR^2$ and $\RR_+{\mathcal K}_n$.
The inverse of $\psi$ is given by
$$
(\lambda,\mu,\nu)\mapsto (\sigma(\lambda,\mu,\nu),\lambda_n,\mu_n).
$$

Recall that if $\sigma=\langle \alpha,\cdot\rangle$, then $\dim \SI(Q,\beta)_{\sigma}=\alpha\circ \beta$.
The numbers $\alpha\circ\beta$ can be interpreted
as Littlewood-Richardson coefficients  even if
$\alpha,\beta$ are arbitrary dimension vectors  and$Q=T_{p,q,r}$ where
$p,q,r$ are arbitrary.
This can be done as follows.
Let us define
$$
\widetilde{c}_{\lambda,\mu,\nu}=\widetilde{c}_{\lambda,\mu,\nu}^{(n)}=\dim (V_{\lambda}\otimes V_\mu\otimes V_\nu)^{\SL_n}.
$$
If
$\widetilde{c}_{\lambda,\mu,\nu}\neq 0$, then $|\lambda|+|\mu|+|\nu|$ must be
a multiple of $n$, say $mn$. In that case
$$\widetilde{c}_{\lambda,\mu,\nu}=c_{\lambda,\mu}^{\nu^\star}$$ 
where
$$
\nu^\star=(m-\nu_n,m-\nu_{n-1},\dots,m-\nu_1).
$$
\begin{definition}
Let $\underline{x}=(x_1,x_2,\dots,x_n)$, $\underline{y}=(y_1,y_2,\dots,y_n)$ be two 
nondecreasing
sequences of nonnegative integers. We define a partition
$P(\underline{x},\underline{y})$ by 
$$
P(\underline{x},\underline{y})=
(x_{n-1}^{y_n-y_{n-1}},x_{n-2}^{y_{n-1}-y_{n-2}},\dots,x_1^{y_2-y_1},0^{y_1}).
$$
Graphically, this partition can be found as follows.
In the Euclidean plane we draw
a square with vertices $(0,0)$, $(x_n,0)$, $(0,y_n)$, $(x_n,y_n)$.
and  we plot the points
$$(x_1,y_1),(x_2,y_2),\dots,(x_n,y_n).$$ 
We take the region in the
square above and left to those points.
Viewed in the unit grid, the corresponding partition $P(\underline{x},\underline{y})$
can be read off.
\end{definition}
\begin{example}
If $\underline{x}=(2,4,5,6)$ and $\underline{y}=(1,3,3,4)$ then
$P(\underline{x},\underline{y})=(5,2,2,0)$ by the diagram below:\\
\centerline{\includegraphics[width=2in]{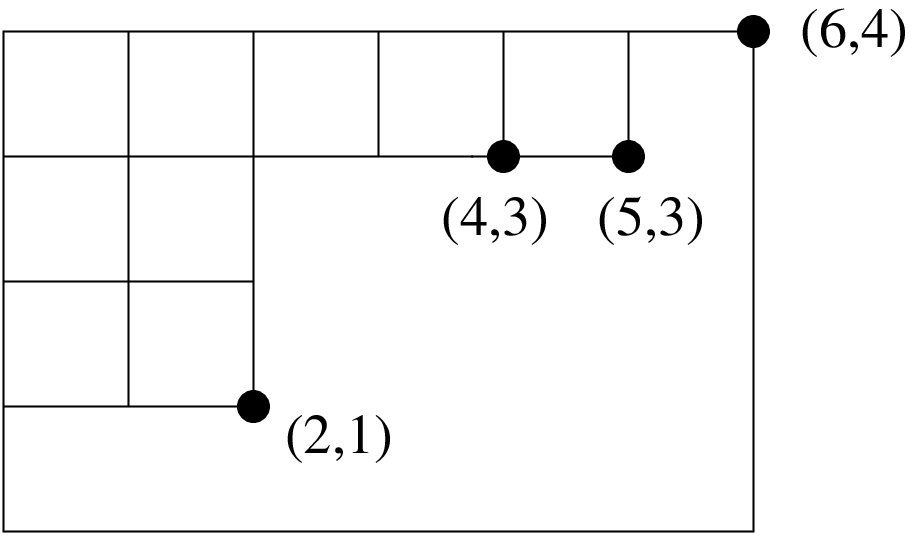}}\\[10pt]
\end{example}
We consider the quiver $T_{p,q,r}$.
Let $\alpha$ and $\beta$ be dimension vectors.
We write $\alpha(x)=(\alpha(x_1),\dots,\alpha(x_p))$,
$\alpha(y)=(\alpha(y_1),\dots,\alpha(y_q))$ and in a similar way
we define $\alpha(z),\beta(x),\beta(y),\beta(z)$.
\begin{lemma}
$$
\alpha\circ \beta=\widetilde{c}_{\lambda,\mu,\nu}
$$
where
$$
\lambda=P(\alpha(x),\beta(x)),\ \mu=P(\alpha(y),\beta(y)),\ 
\nu=P(\alpha(z),\beta(z)).
$$
\end{lemma}
\begin{proof}
This is an easy computation, see \cite{DW}.
\end{proof}
\begin{remark}\label{rempqr}
Note that if $1/p+1/q+1/r>1$, then $T_{p,q,r}$ is a quiver of finite type and
for all dimension vectors $\alpha,\beta$ we have $\alpha\circ \beta=0$
or $\alpha\circ\beta=1$.
For a partition $\lambda=(\lambda_1,\dots,\lambda_n)$ we define
$$
j(\lambda)=\#\{i\mid \lambda_{i+1}\neq \lambda_i,1\leq i\leq n-1\}
$$
If $\lambda$ is trivial then $j(\lambda)=0$, if $\lambda$ is a box
then $j(\lambda)=1$ and if $\lambda$ is L-shaped (fat hook) then $j(\lambda)=2$.
A coefficient $\widetilde{c}_{\lambda,\mu,\nu}$ can be obtained
from the quiver $T_{p,q,r}$ with $p=j(\lambda)+1,q=j(\mu)+1$ and $r=j(\nu)+1$.
In particular, we get the corollary below:
\end{remark}
\begin{corollary}
If
 $$
\frac{1}{j(\lambda)+1}+\frac{1}{j(\mu)+1}+\frac{1}{j(\nu)+1}>1
$$
then $\widetilde{c}_{\lambda,\mu,\nu}$ equals 0 or 1.
\end{corollary}
\begin{example}
Let us consider the quiver $T_{3,3,2}$:
$$
\xymatrix{
\circ\ar[r] & \circ\ar[rd]\\
\circ\ar[r] & \circ\ar[r] & \circ\\
& \circ\ar[ru] &}
$$
Let $\alpha$, $\beta$ be the dimension vectors
$$
\alpha=
\begin{matrix}
1 & 3 &   \\
1 & 2 & 4  \\
  & 2 
\end{matrix},\qquad
\beta=
\begin{matrix}
1 & 2 &   \\
0 & 2 & 3  \\
  & 1 
\end{matrix}.
$$
Now $\alpha\circ\beta$ is equal to the LR-coefficient 
$\widetilde{c}_{\lambda,\mu,\nu}=1$
were $\lambda=(3,1)$, $\mu=(2,1,1)$ and $\nu=(2,2)$.\\
\centerline{\includegraphics[width=4in]{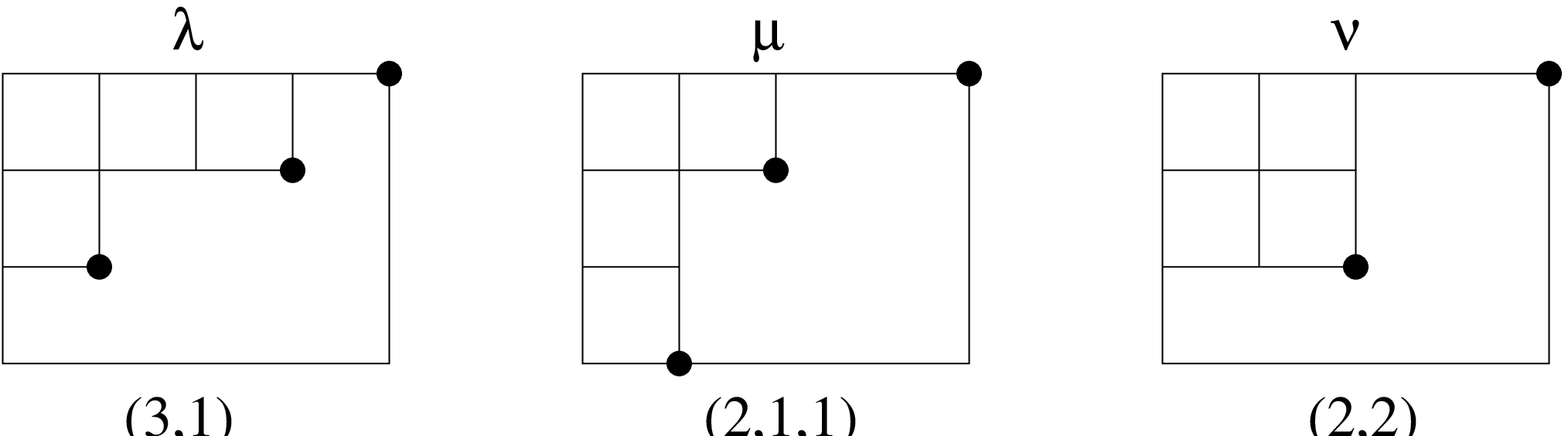}}\\[10pt]
\end{example}
\subsection{Walls of the Klyachko cone}
Let us consider the quiver $Q=T_{n,n,n}$ and the dimension vector
$$
\beta=
\begin{matrix}
1 & 2 & \cdots & n-1 & \\
1 & 2 & \cdots & n-1 & n \\
1 & 2 & \cdots & n-1 & \end{matrix}.
$$
\begin{lemma}\label{Schurroot}
The dimension vector $\beta$ above is a Schur root.
\end{lemma}
\begin{proof}
For $n\leq 2$ this is an easy check that $\beta$ is a real Schur root.
For $n\geq 3$
the dimension vector $\beta$ is indivisible and lies in the {\it fundamental set\/}
(see \cite[\S 1]{K2}). This implies that $\beta$ is a Schur root
by \cite[Theorem B (d)]{K2}).
\end{proof}
We will study the cone $\RR_+\Sigma(Q,\beta)$ which is essentially the Klyachko cone.
The following Theorem is equivalent to a result of Knutson, Tao and Woodward (see~\cite{KTW}).
It gives a precise description of the walls of the Klyachko cone.
\begin{theorem}
For every pair $(\beta_1,\beta_2)$ with $\beta=\beta_1+\beta_2$,
$\beta_1,\beta_2$ nondecreasing along arms,
$\beta_1\circ \beta_2=1$
the inequality $\sigma(\beta_1)\leq 0$ defines
a wall of $\RR_+\Sigma(Q,\beta)$. All nontrivial walls 
can be uniquely obtained this way.
\end{theorem}
\begin{proof}
Clearly $\beta_1$ and $\beta_2$ have at most jumps 1 along the arms
of the quiver $Q$. It follows from
Lemma~\ref{Schurroot} that $\beta_1,\beta_2$ are Schur roots. 
Now either $\ext_Q(\beta_2,\beta_1)=0$ or $\hom_Q(\beta_2,\beta_1)=0$ by Theorem~\ref{theohomorext}
The first would give a nontrivial decomposition of $\beta$, therefore
$\hom_Q(\beta_2,\beta_1)=0$ and $\beta_1,\beta_2$ is a quiver Schur sequence.
This shows that $\sigma(\beta_1)\leq 0$ defines a wall by Theorem~\ref{bijection}.

For every wall, there exists a Schur sequence $(\beta_1,\beta_2)$ such that
that $\beta=c_1\beta_1+c_2\beta_2$ with $c_i$ positive integers
and the wall is defined by $\sigma(\beta_1)\leq 0$ by Theorem~\ref{bijection}. Note
that for $T_{n,n,n}$ a Schur root either has support on one arm
(in which case it corresponds to a positive root of $A_{n-1}$),
or it is nondecreasing along each arm.
Because $\beta_1\hookrightarrow \beta$, it is easy
to see that $\beta_1$ must also be nondecreasing. Indeed, for a general representation
of dimension $\beta$, the linear maps along the arms are injective, so these
maps are also injective for every subrepresentation.
But $\beta_2$ could
have support on one arm. In that case it follows from $\langle
\beta_1,\beta_2\rangle=0$ that $\beta_2$ is simple. 
The inequalities following from such a quiver sequence are trivial,
they say that the partitions $\lambda$, $\mu$ and $\nu$ must be weakly decreasing.
If $\beta_1$ and $\beta_2$ are both nondecreasing along the arms, then
$c_1=c_2=1$ because $\beta=\beta_1+\beta_2$ and $\beta$ jumps
only by steps of 1 along the arms.
\end{proof}
We need to find all $\beta_1,\beta_2$ such that $\beta=\beta_1+\beta_2$
and $\beta_1\circ \beta_2=1$. In that case $\sigma(\beta_1)\leq 0$ defines
a wall of $\Sigma(Q,\beta)$. If we are dealing with a nontrivial
wall, i.e., {\it not\/} $\lambda_i\geq\lambda_{i+1}$, $\mu_i\geq \mu_{i+1}$
or $\nu_{i}\geq \nu_{i+1}$ for some $i$, then both $\beta_1$ and $\beta_2$
are increasing along the arms, with jumps at most 1.
The fact that $\beta_1$ has jumps of at most 1, gives these
inequalities a special form, namely, if we take
$$
I=\{i\mid \beta_1(x_{i-1})=\beta_1(x_{i}),1\leq i\leq n\},
$$
$$ 
J=\{i\mid \beta_1(y_{i-1})=\beta_1(y_{i}),1\leq i\leq n\},
$$
$$ 
K=\{i\mid \beta_1(z_{i-1})=\beta_1(z_{i}),1\leq i\leq n\}, 
$$
(by convention $\beta_1(x_0)=\beta_1(y_0)=\beta_1(z_0)=0$),
then the inequality corresponding to $\beta_1$ is
\begin{equation}\label{IJK}
\sum_{i\in I}\lambda_i+\sum_{i\in J}\mu_i\leq \sum_{i\in K}\nu_{n-i}
\end{equation}
Note that $\#I=\#J=\#K=\beta_1(x_n)$.
Now (\ref{IJK}) is a necessary inequality for the
Klyachko cone if $\beta_1\circ \beta_2=1$. If $\beta_1\circ\beta_2>0$
then (\ref{IJK}) still defines a true inequality.
Now $\beta_1\circ\beta_2$ is the value
of a Littlewood-Richardson coefficient for $\SL_{\beta_2(x_n)}$.
Since $\beta_2(x_n)<n$ we know necessary and sufficient inequalities
for the corresponding LR-coefficient to be nonzero.
{\it This explains the inductive nature of the Klyachko inequalities}.
\begin{example}
Consider the quiver $T_{3,3,3}$ and 
$$
\beta=\begin{matrix}
1 & 2 &\\
1 & 2 & 3 \\
1 & 2 &
\end{matrix}.
$$
The LR-coefficient $c_{\lambda,\mu}^\nu$ corresponds
to $\dim \SI(Q,\beta)_\sigma$ where $\sigma$ is given by
$$
\sigma=
\begin{matrix}
\lambda_1-\lambda_2 & \lambda_2-\lambda_3 &\\
\mu_1-\mu_2 & \mu_2-\mu_3 & \lambda_3+\mu_3-\nu_1 \\
\nu_2-\nu_3 & \nu_1-\nu_2 &
\end{matrix}.
$$
For example the sequence
$$
\begin{matrix}
1 & 2 &\\
0 & 1 & 2 \\
0 & 1 &
\end{matrix},\ 
\begin{matrix}
0 & 0 &\\
1 & 1 & 1 \\
1 & 1 &
\end{matrix}
$$
is a Schur sequence, because $\widetilde{c}^1_{2,0,0}=c_{2,0}^2=1$.
This Schur sequence corresponds to the wall
$$
\lambda_1+\lambda_2+\mu_2+\mu_3\leq \nu_1+\nu_2
$$
By permuting the arms, we get the inequalities
\begin{eqnarray*}
\lambda_2+\lambda_3+\mu_1+\mu_2 &\leq& \nu_1+\nu_2,\\
\lambda_2+\lambda_3+\mu_2+\mu_3 &\leq& \nu_2+\nu_3.
\end{eqnarray*}
Other walls are given by the Schur sequences
$$
\begin{matrix}
1 & 1 &\\
1 & 1 & 2 \\
0 & 1 &
\end{matrix},\ 
\begin{matrix}
0 & 1 &\\
0 & 1 & 1 \\
1 & 1 &
\end{matrix}\quad (\widetilde{c}^{(1)}_{1,1,0}=c_{1,1}^2=1).
$$
By permuting arms we get the inequalities
\begin{eqnarray*}
\lambda_1+\lambda_3+\mu_1+\mu_3 &\leq& \nu_1+\nu_2,\\
\lambda_1+\lambda_3+\mu_2+\mu_3 &\leq& \nu_1+\nu_3,\\
\lambda_2+\lambda_3+\mu_1+\mu_3 &\leq& \nu_1+\nu_3.
\end{eqnarray*}
$$
\begin{matrix}
1 & 1 &\\
0 & 0 & 1 \\
0 & 0 &
\end{matrix},\ 
\begin{matrix}
0 & 1 &\\
1 & 2 & 2 \\
1 & 2 &
\end{matrix} \quad(\widetilde{c}^{(2)}_{11,0,0}=c_{11,0}^{11}=1)
$$
gives the inequalities
\begin{eqnarray*}
\lambda_1+\mu_3 &\leq& \nu_1\\
\lambda_3+\mu_1 & \leq & \nu_1\\
\lambda_3+\mu_3 &\leq& \nu_3
\end{eqnarray*}
The Schur sequence
$$
\begin{matrix}
0 & 1 &\\
0 & 1 & 1 \\
0 & 0 &
\end{matrix},\ 
\begin{matrix}
1 & 1 &\\
1 & 1 & 2 \\
1 & 2 &
\end{matrix} \quad(\widetilde{c}^{(2)}_{1,1,0}=c_{1,1}^{11}=1)
$$
gives the inequalities
\begin{eqnarray*}
\lambda_2+\mu_2 &\leq & \nu_1\\
\lambda_2+\mu_3 &\leq & \nu_2\\
\lambda_3+\mu_2 &\leq & \nu_2
\end{eqnarray*}
Besides these, there are 6 trivial walls corresponding to the inequalities
$\lambda_1\geq\lambda_2\geq\lambda_3$, $\mu_1\geq\mu_2\geq\mu_3$
and  $\nu_1\geq\nu_2\geq\nu_3$. The Schur sequence
$$
\begin{matrix}
0 & 2 &\\
1 & 2 & 3 \\
1 & 2 &
\end{matrix},\ 
\begin{matrix}
1 & 0 &\\
0 & 0 & 0 \\
0 & 0 &
\end{matrix}
$$
implies the inequalities
\begin{eqnarray*}
\lambda_1&\geq&\lambda_2,\\
\mu_1&\geq& \mu_2,\\
\nu_2&\geq& \nu_3.\\
\end{eqnarray*}
The Schur sequence
$$
\begin{matrix}
1 & 1 &\\
1 & 2 & 3 \\
1 & 2 &
\end{matrix},\ 
\begin{matrix}
0 & 1 &\\
0 & 0 & 0 \\
0 & 0 &
\end{matrix}
$$
leads to the inequalities
\begin{eqnarray*}
\lambda_3&\geq&\lambda_3\\
\mu_2&\geq& \mu_3\\
\nu_1&\geq &\nu_2
\end{eqnarray*}
\end{example}

\subsection{Faces of the Klyachko cone of arbitrary codimension}
For the quiver $Q=T_{n,n,n}$, Theorem~\ref{bijection} translates to:
\begin{corollary}\label{walls}
There is a 1--1 correspondence between the faces of
codimension $l$ of ${\mathcal K}_n$  and Schur sequences
$(\beta_1,\dots,\beta_{l+1})$ such that $\beta=c_1\beta_1+\cdots+c_{l+1}\beta_{l+1}$
for some positive integers $c_1,\dots,c_{l+1}$, $\beta_1,\dots,\beta_{l+1}$
are Schur roots, and $\beta_i\circ \beta_j=1$ for all $i<j$.
\end{corollary}
However, If $l>1$ then the $c_i$ may be larger than 1.
There is no easy criterion for a dimension vector to be a Schur root,
but a fast algorithm for determining whether a dimension
vector is a Schur root was given in \cite{DW2}.
This makes it more difficult to find the faces of higher codimension.
Still, we obtain some interesting features.
\begin{corollary}\label{corjumps}
Suppose that $(\lambda,\mu,\nu)$ lies in a face $F$ of ${\mathcal K}_n$ of 
codimension $l$. Let $j(\lambda),j(\mu),j(\nu)$ be the number
of jumps in $\lambda,\mu,\nu$ respectively.
\begin{enumerate}
\renewcommand{\theenumi}{\alph{enumi}}
\item
$$
j(\lambda)+j(\mu)+j(\nu)\leq 4n-4-l,
$$
\item if $\widetilde{c}_{\lambda,\mu,\nu}>1$ then
$$
j(\lambda)+j(\mu)+j(\nu)\leq 4n-6-l.
$$
\end{enumerate}
\end{corollary}
\begin{proof}
Let $\beta$ be the usual dimension vector for $T_{n,n,n}$.
Now $\alpha\circ\beta=\widetilde{c}_{\lambda,\mu,\nu}$ for some
dimension vector $\alpha$. Let $\sigma=\langle\alpha,\cdot\rangle$.
Because
$\sigma$ is in a face of $\Sigma(Q,\beta)$ of codimension
$l$, exactly
$l+1$ distinct dimension vectors appear in the $\sigma$-stable
decomposition of $\beta$.
Suppose that the $\sigma$-stable decomposition of $\beta$ is
$$
\beta=c_1\cdot\beta_1\pp c_2\cdot\beta_2 \cdots \pp c_{l+1}\cdot \beta_{l+1}.
$$
Whenever $\beta_i(x_n)=0$, then $\beta_i$ has support on one arm
and the equation $\sigma(\beta_i)=0$ corresponds to the equation
$\lambda_j=\lambda_k$, $\mu_j=\mu_k$ or $\nu_j=\nu_k$ for some $j\neq k$.
Notice also that all $\beta_j$'s correspond to linearly independent
equations.
It now follows that
$$
j(\lambda)+j(\mu)+j(\nu)+\#\{i\mid \beta_i(x_n)=0\}\leq 3n-3.
$$
Since $\sum_j c_j\beta_j(x_n)=n$, we have 
$\#\{j\mid \beta_j(x_n)>0\}\leq n$ and
$\#\{j\mid \beta_j(x_n)=0\}\geq l+1-n$.
Now (a) follows.

If $\widetilde{c}_{\lambda,\mu,\nu}>1$, then at least
one of the $\beta_i$'s is an imaginary Schur root, so $\beta_{i}(x_n)\geq 3$.
This shows that $\#\{j\mid \beta_j(x_n)>0\}\leq n-2$ and (b) follows.
\end{proof}
The cone $\RR_+\Sigma(Q,\beta)$ has one 0-dimensional face, namely $\{0\}$.
This corresponds to the the 2-dimensional face of 
$\RR_+{\mathcal K}_n$ consisting of all $(\lambda,\mu,\nu)$,
$\lambda=(a,\dots,a)$, $\mu=(b,\dots,b)$, $\nu=(a+b,\dots,a+b)$.
It would be also interesting  to study the extremal rays of the cone
$\RR_+\Sigma(Q,\beta)$ (or equivalently the 3-dimensional faces of $\RR_+{\mathcal K}_n$).
They span the cone $\Sigma(Q,\beta)$ (or the Klyachko cone).
The codimension of the extremal rays is $3n-4$.
One interesting question is, whether $c_{\lambda,\mu}^{\nu}=1$ 
whenever $(\lambda,\mu,\nu)$ is on an extremal ray of $\RR_+{\mathcal K}_n$ (here
by {\it extremal ray\/} we mean a 3-dimensional face of ${\mathcal K}_n$).
We first give a positive result in this direction:
\begin{corollary}
If $(\lambda,\mu,\nu)$ is in an extremal ray of $\RR_+{\mathcal K}_n$, and
$n\leq 7$, then $c_{\lambda,\mu}^{\nu}=1$.
\end{corollary}
\begin{proof}
Suppose that $c_{\lambda,\mu}^{\nu}>1$.
Then $\widetilde{c}_{\lambda,\mu,\nu^\star}=c_{\lambda,\mu}^\nu>1$ for some partition $\nu^\star$.
Then $j(\lambda)+j(\mu)+j(\nu^\star)\geq 6$ (this
follows from Remark~\ref{rempqr}). So
$$
6\leq 4n-6-(3n-4)
$$
by Corollary~\ref{corjumps} and we deduce that $n\geq 8$. Contradiction.
\end{proof}
\begin{example}
We study $T_{8,8,8}$ with the weight
$$
\sigma=
\begin{matrix}
1 & 0 & 0 & 0 & 1 & 0 & 0 & \\
0 & 0 & 1 & 0 & 0 & 1 & 0 & -3 \\
0 & 0 & 1 & 0 & 0 & 1 & 0 & 
\end{matrix}
$$
The $\sigma$-stable decomposition of 
$$
\beta=
\begin{matrix}
1 & 2 & 3 & 4 & 5 & 6 & 7 &\\
1 & 2 & 3 & 4 & 5 & 6 & 7 & 8\\
1 & 2 & 3 & 4 & 5 & 6 & 7 &
\end{matrix}
$$
is
$$
\begin{matrix}
1 & 1 & 1 & 1 & 2 & 2 & 2 & \\
0 & 0 & 1 & 1 & 1 & 2 & 2 & 3\\
0 & 0 & 1 & 1 & 1 & 2 & 2 &
\end{matrix}
\ \pp\ 
\begin{matrix}
0 & 0 & 0 & 0 & 1 & 1 & 1 & \\
0 & 0 & 1 & 1 & 1 & 1 & 1 & 1\\
0 & 0 & 0 & 0 & 0 & 0 & 0 &
\end{matrix}
\ \pp
$$
$$
\begin{matrix}
0 & 0 & 0 & 0 & 1 & 1 & 1 & \\
0 & 0 & 0 & 0 & 0 & 0 & 0 & 1\\
0 & 0 & 1 & 1 & 1 & 1 & 1 & \\
\end{matrix}
\ \pp\ 
\begin{matrix}
0 & 0 & 0 & 0 & 0 & 0 & 0 & \\
0 & 0 & 0 & 0 & 0 & 1 & 1 & 1\\
0 & 0 & 1 & 1 & 1 & 1 & 1 & \\
\end{matrix}
\ \pp
$$
$$
\begin{matrix}
0 & 0 & 0 & 0 & 0 & 0 & 0 & \\
0 & 0 & 1 & 1 & 1 & 1 & 1 & 1\\
0 & 0 & 0 & 0 & 0 & 1 & 1 & \\
\end{matrix}
\ \pp\ 
\begin{matrix}
0 & 0 & 0 & 0 & 1 & 1 & 1 & \\
0 & 0 & 0 & 0 & 0 & 1 & 1 & 1\\
0 & 0 & 0 & 0 & 0 & 1 & 1 & \\
\end{matrix}\ \pp
$$
$$
\delta_{x_2}\pp 2\cdot \delta_{x_3}\pp 3\cdot \delta_{x_4}+\delta_{x_6}+2\cdot
\delta_{x_7}\pp 
$$
$$
\delta_{y_1}\pp 2\cdot \delta_{y_2}\pp \delta_{y_4}\pp 2\cdot \delta_{y_5}\pp
\delta_{y_7}\pp
$$
$$
\delta_{z_1}\pp 2\cdot \delta_{z_2}\pp \delta_{z_4}\pp 2\cdot \delta_{z_5}\pp
\delta_{z_7}.
$$
Here for any vertex $p$ of the quiver $Q=T_{8,8,8}$ $\delta_p$ denotes
the dimension vector of the simple representation corresponding
to the vertex $p$.
In the $\sigma$-stable decomposition of $\beta$, there are 21 distinct Schur roots.
The quiver $T_{8,8,8}$ has 22 vertices, so this proves that $\sigma$ is
in an extremal ray of $\Sigma(Q,\beta)$.
The LR-Richardson coefficient corresponding to $\beta$ and $\sigma$
is $\widetilde{c}_{\lambda,\mu,\nu}$ where 
$$
\lambda=(2,1,1,1,1,0,0,0),\ \mu=(2,2,2,1,1,1,0,0),\ \nu=(2,2,2,1,1,1,0,0).
$$
The value of $\widetilde{c}_{\lambda,\mu,\nu}$ is 2. In fact,
for any $N$ we have $\widetilde{c}_{N\lambda,N\mu,N\nu}=N+1$.
\end{example}
\subsection{A multiplicative formula for Littlewood-Richardson coefficients}
Let $\beta$ and $T_{n,n,n}$ as before. 
\begin{theorem}\label{propprod}
Suppose $\beta=\beta_1+\beta_2$ and $\beta_1\circ\beta_2=1$.
Let $\alpha$ be another dimension vector with $\alpha\circ
\beta=\widetilde{c}_{\lambda,\mu,\nu}$. Put $\sigma=\langle \alpha,\cdot\rangle$.
The inequality $\sigma(\beta_1)\leq 0$ translates to
\begin{equation}\label{eqIJK}
\sum_{i\in I}\lambda_i+\sum_{i\in J}\mu_i\leq \sum_{i\in K}\nu_i.
\end{equation}
where $I,J,K$ are subsets of $S=\{1,2,\dots,n\}$ of the same cardinality.
Suppose that equality in (\ref{eqIJK}) holds for $(\lambda,\mu,\nu)\in (\ZZ^n)^3$.
Define
$$\lambda^*=(\lambda_{i_1},\cdots,\lambda_{i_r}),\quad I=\{i_1,i_2,\dots,i_r\},$$
$$\lambda^\#=(\lambda_{\widetilde{\i}_1},\cdots,\lambda_{\widetilde{\i}_{n-r}}),
\quad S\setminus I=\{\widetilde{\i}_1,\widetilde{\i}_2,\dots,\widetilde{\i}_{n-r}\},$$
$$\mu^*=(\mu_{j_1},\cdots,\mu_{j_r}),\quad J=\{j_1,j_2,\dots,j_r\},$$
$$\mu^\#=(\mu_{\widetilde{\j}_1},\cdots,\mu_{\widetilde{\j}_{n-r}}),
\quad S\setminus J=\{\widetilde{\j}_1,\widetilde{\j}_2,\dots,\widetilde{\j}_{n-r}\},$$
$$\nu^*=(\nu_{k_1},\cdots,\nu_{j_r}),\quad K=\{k_1,k_2,\dots,k_r\},$$
$$\nu^\#=(\nu_{\widetilde{k}_1},\cdots,\nu_{\widetilde{k}_{n-r}}),
\quad S\setminus K=\{\widetilde{k}_1,\widetilde{k}_2,\dots,\widetilde{k}_{n-r}\}.$$
Then we have
$$
c_{\lambda,\mu}^{\nu}=c_{\lambda^*,\mu^*}^{\nu^*}c_{\lambda^\#,\mu^\#}^{\nu^\#}.
$$
\end{theorem}
\begin{proof}
If $\sigma$ is in the interior of the wall, then the
$\sigma$-stable decomposition of $\beta$ is $\beta_1\pp \beta_2$ and by
Theorem~\ref{invar} we get
$$
c_{\lambda,\mu}^{\nu}=\alpha\circ\beta=
(\alpha\circ\beta_1)(\alpha\circ\beta_2)=c_{\lambda^*,\mu^*}^{\nu^*}
c_{\lambda^{\#},\mu^{\#}}^{\nu^{\#}}.
$$

Assume $\sigma$ is not in the interior of the wall.
Suppose that the $\sigma$-stable decomposition of $\beta_1$ is 
$$
c_1\cdot \gamma_1\pp \cdots \pp c_r\cdot \gamma_r
$$
and that the $\sigma$-stable decomposition of $\beta_2$ is 
$$
d_1\cdot \delta_1\pp \cdots \pp d_s\cdot \delta_s
$$
Then the $\sigma$-stable decomposition of $\beta$ is
$$
c_1\cdot \gamma_1\pp \cdots \pp c_r\cdot \gamma_r\pp
d_1\cdot \delta_1\pp \cdots \pp d_s\cdot \delta_s.
$$
Note that $\{\gamma_1,\dots,\gamma_r\}$ and $\{\delta_1,\dots,\delta_s\}$
are disjoint because $\gamma_i\pperp\delta_j$ for all $i,j$. 
By Theorem~\ref{invar} we get
$$
c_{\lambda,\mu}^\nu=\alpha\circ\beta=\prod (\alpha\circ (c_i\gamma_i))\prod (\alpha\circ
(d_i\delta_i))=(\alpha\circ\beta_1)(\alpha\circ\beta_2)=c_{\lambda^*,\mu^*}^{\nu^*}
c_{\lambda^{\#},\mu^{\#}}^{\nu^{\#}}.
$$
\end{proof}
\begin{example}
For $n=8$, $\beta=\beta_1+\beta_2$ with
$$
\beta_1=\begin{matrix}
1 & 1 & 2 & 2 & 3 & 3 & 4 & \\
1 & 1 & 2 & 2 & 3 & 3 & 4 & 5\\
0 & 0 & 1 & 2 & 3 & 3 & 4 & \\
\end{matrix},\ 
\beta_2=\begin{matrix}
0 & 1 & 1 & 2 & 2 & 3 & 3 & \\
0 & 1 & 1 & 2 & 2 & 3 & 3 & 3\\
1 & 2 & 2 & 2 & 2 & 3 & 3 & \\
\end{matrix}
$$
Now $\beta_1\circ \beta_2=\widetilde{c}^{(3)}_{321,321,3}=c_{321,321}^{552}$, so
$\beta_1,\beta_2$ is a Schur sequence.
The corresponding inequality for the Klyachko cone is
$$
\lambda_1+\lambda_3+\lambda_5+\lambda_7+\lambda_8+
\mu_1+\mu_3+\mu_5+\mu_7+\mu_8\leq 
\nu_1+\nu_2+\nu_4+\nu_5+\nu_6.
$$
or equivalently
$$
\lambda_2+\lambda_4+\lambda_6+\mu_2+\mu_4+\mu_6\geq \nu_3+\nu_7+\nu_8.
$$
If now these inequalities are equalities, then
$$
c_{\lambda,\mu}^{\nu}=c_{\lambda^*,\mu^*}^{\nu^*}c_{\lambda^\#,\mu^\#}^{\nu^\#}
$$
where
$$
\lambda^*=(\lambda_1,\lambda_3,\lambda_5,\lambda_7,\lambda_8),\ 
\mu^*=(\mu_1,\mu_3,\mu_5,\mu_7,\mu_8),\ 
\nu^*=(\nu_1,\nu_2,\nu_4,\nu_5,\nu_6)$$
and
$$
\lambda^\#=(\lambda_2,\lambda_4,\lambda_6),\ 
\mu^\#=(\mu_2,\mu_4,\mu_6),\ 
\nu^\#=(\nu_3,\nu_7,\nu_8).
$$
For example, take $\lambda=\mu=(8,4,4,2,2,0,0,0)$, $\nu=(10,8,7,4,3,3,3,2)$.

Then $c_{\lambda,\mu}^{\nu}=c_{\lambda^*,\mu^*}^{\nu^*}c_{\lambda^\#,\mu^\#}^{\nu^\#}$
where $\lambda^*=\mu^*=(8,4,2,0,0)$,
$\nu^*=(10,8,4,3,3)$,
$\lambda^\#=\mu^\#=(4,2,0)$,
$\nu^*=(-7,3,2)$.
Indeed, $c_{\lambda^*,\mu^*}^{\nu^*}=5$, $c_{\lambda^\#,\mu^\#}^{\nu^\#}=2$
and $c_{\lambda,\mu}^{\nu}=10$.
\end{example}

\section*{Appendix: Belkale's proof of Fulton's conjecture}
At the {\it AMS Summer Institute on Algebraic Geometry Meeting} in Seattle, Belkale
explained his geometric proof of Theorem~\ref{FultonConjecture} 
to the first author and how it generalizes to the more
general quiver setting. What follows below is a reconstruction of that proof.
We are grateful to Prakash Belkale for letting us include his proof in our paper.

For  a nonnegative integer $d$, the Grassmannian of $d$-dimensional subspaces of
the $n$-dimensional vector space $V$ is denoted by
$$
\textstyle\Gr{V\choose d}.
$$
Suppose that $\alpha,\gamma$ are dimension vectors for a quiver $Q$. We define
$$
{\textstyle \Gr{\alpha\choose \gamma}}=\prod_{x\in Q_0}{\textstyle \Gr{K^{\alpha(x)}\choose \gamma(x)}}.
$$
$$
\Hom(K^{\alpha},K^{\beta})=\prod_{x\in Q_0}\Hom(K^{\alpha(x)},K^{\beta(x)}).
$$
Following Schofield, we need to introduce the notion of the {\it general rank} of a
morphism.
For a morphism $\phi:V\to W$ between two representations of $Q$, 
define the rank function $\rk(\phi)\in \NN^{Q_0}$
by
$$
\rk(\phi)(x)=\dim_K \phi(x)(V(x)),\quad x\in Q_0.
$$
Define the variety 
\begin{multline*}
H(\alpha,\beta)=\{(V,W,\phi)\in \Rep(Q,\alpha)\times \Rep(Q,\beta)\times \Hom(K^\alpha,K^\beta)\mid\\
\forall a\in Q_1 \phi(ha)V(a)=W(a)\phi(ta)\}.
\end{multline*}
There is a natural projection $p:H(\alpha,\beta)\to \Rep(Q,\alpha)\times \Rep(Q,\beta)$
such that the fiber of $(V,W)\in \Rep(Q,\alpha)\times \Rep(Q,\beta)$ can be identified
with $\Hom_Q(V,W)$. Suppose that $Z\subseteq \Rep(Q,\alpha)\times \Rep(Q,\beta)$ is an irreducible
closed subset. There exists an open dense subset $U\subseteq Z$ such
$U$ is smooth and the fibers of  $p^{-1}(U)\to U$ have constant dimension. 
Then $p^{-1}(U)$ is irreducible. The rank depends semi-continuously on $(V,W,\phi)\in p^{-1}(U)$.
So there exists an open dense subset $U'\subseteq p^{-1}(U)$ such that
$\rk(\phi)$ is constant on $U'$. The value of this rank is called the {\it general rank} of
a morphism $\phi:V\to W$ with $(V,W)\in Z$.

By symmetry, The Generalized Fulton Conjecture (Theorem~\ref{thmBelkale}) follows from
\begin{theorem}
If $\alpha\circ \beta=1$, then $\alpha\circ (n\beta)=1$ for all $n\geq 0$.
\end{theorem}
\begin{proof}
Assume that $\alpha\circ\beta=1$ and $\alpha\circ (n\beta)>1$ for some
$n\geq 0$. Let us put $\sigma=-\langle \cdot, \beta\rangle$.
Then
$$
\dim \SI(Q,\alpha)_{n\sigma}=\alpha\circ (n\beta).
$$
We construct a quotient as in Section~\ref{secKing}.
Define the projective variety
$$
Y=\operatorname{Proj} \bigoplus_{n\geq 0} \SI(Q,\alpha)_{n\sigma}.
$$
It has dimension $>0$. 
Let 
$$
\pi:\Rep(Q,\alpha)^{\rm ss}_\sigma\to Y.
$$
be the GIT quotient with respect to $\sigma$, where $\Rep(Q,\alpha)^{\rm ss}_\sigma\subseteq \Rep(Q,\alpha)$
is the dense, open subset of $\sigma$-semi-stable representations.
Let $W\in \Rep(Q,\beta)$ in general position. 
The Schofield semi-invariant $c_W\in \SI(Q,\alpha)_{\sigma}$ is nonzero
because $W$ is in general position. By assumption, $\dim_K \SI(Q,\alpha)_\sigma=\alpha\circ \beta=1$.
So
$\SI(Q,\alpha)_\sigma$ is spanned by $c_W$. Also, if $W'\in \Rep(Q,\beta)$ is any
other representation then either $c_{W'}$ is identically 0 or $c_{W'}$ is equal to
$c_W$ up to a scalar.
By assumption, $\dim \SI(Q,\alpha)_{n\sigma}=\alpha\circ (n\beta)>1$ for some $n$.
If $f_1,f_2\in \SI(Q,\alpha)_{n\sigma}$ are linearly independent, then they
must be algebraically independent. This implies that the Krull dimension
of the ring $\bigoplus_{n\geq 0} \SI(Q,\alpha)_{n\sigma}$ is at least $2$,
and the dimension of $Y$ is at least $1$.
The equation $c_W=0$ defines a nonnegative divisor on the projective variety $Y$.
Clearly this divisor is nonzero, because $c_W$ is not a constant.
  Let $y\in Y$ such that $C_W(y)=0$. 
Since $\pi$ is surjective, we can choose $V\in \pi^{-1}(y)\subseteq \Rep(Q,\alpha)^{\rm ss}_\sigma$.
It follows that $c_W(V)=0$.
Let $D\subseteq \Rep(Q,\alpha)$ be an irreducible component
of the divisor
$$
\{Z\in \Rep(Q,\alpha)\mid c_W(Z)=0\}
$$
which contains $V$. In conclusion, $D$ is a divisor of $\Rep(Q,\alpha)$ containing
$\sigma$-semi-stable representations, and $c_W(V)=0$ for all $(V,W)\in D\times \Rep(Q,\beta)$.

Let $\gamma$ be the general rank of a homomorphism
$\phi:V\to W$ where $(V,W)\in D\times \Rep(Q,\beta)$.

Define
\begin{multline*}
\widetilde{M}=\{
(V,W,V_1,W_1,\phi)\in\Rep(Q,\alpha)\times \Rep(Q,\beta)\times
{\textstyle\Gr{\alpha\choose\alpha-\gamma}\times \Gr{\beta\choose\gamma}}\times
\Hom(K^\alpha,K^\beta)\mid\\
\mbox{$\phi:V\to W$ is a morphism},\ \forall\ x\in Q_0\ \phi(V_1(x))=0\mbox{ and }
\phi(K^{\alpha(x)})\subseteq W_1(x)\}.
\end{multline*}
and 
$$
M=\{(V,W,V_1,W_1,\phi)\in \widetilde{M}\mid \phi(K^{\alpha(x)})=W_1(x)\}.
$$
We have that $\widetilde{M}$ is a Zariski closed subset
of 
$$\Rep(Q,\alpha)\times \Rep(Q,\beta)\times
{\textstyle\Gr{\alpha\choose\alpha-\gamma}\times \Gr{\beta\choose\gamma}}\times
\Hom(K^\alpha,K^\beta),$$ 
and $M$ is open in $\widetilde{M}$ and therefore it is locally closed.
So $M$ is a variety. If $(V,W,V_1,W_1,\phi)\in M$ then $\phi:V\to W$ has
kernel $V_1$, image $W_1$ and rank $\gamma=\dd W_1$.

Define
\begin{multline*}
N:=\{
(V_1,W_1,\phi)\in \textstyle \Gr{\alpha\choose \alpha-\gamma}\times \Gr{\beta\choose \gamma}\times
\Hom(K^{\alpha(x)},K^{\beta(x)})\mid\\
\forall x\in Q_0\ \phi(V_1(x))=0,\ \phi(K^{\alpha(x)})=W_1(x)\}.
\end{multline*}
Again $N$ is locally closed, hence a variety.
The projection $p:M\to \Gr{\alpha\choose \alpha-\gamma}\times \Gr{\beta\choose
\gamma}$ factors through a morphism $q:M\to N$ and the projection 
$r:N\to \Gr{\alpha\choose \alpha-\gamma}\times \Gr{\beta\choose \gamma}$.

Let $(V_1,W_1)\in \Gr{\alpha\choose\alpha-\gamma}\times
\Gr{\beta\choose\gamma}$
and consider the fiber $r^{-1}(V_1,W_1)$.
This can be seen as
the set of all $\phi:K^{\alpha}/V_1\to W_1$
which induce an isomorphism at each vertex.
So $r^{-1}(V_1,W_1)\cong \GL(Q,\gamma)$ has dimension
$$
\sum_{x\in Q_0}\gamma(x)^2.
$$

Note that $r$ is actually a fiber bundle.
We get
\begin{multline*}
\dim N=\dim {\textstyle\Gr{\alpha\choose \alpha-\gamma}+\dim \Gr{\beta\choose \gamma}}+
\sum_{x\in Q_0}\gamma(x)^2=\\
=
\sum_{x\in Q_0} \big(\gamma(x)(\alpha-\gamma)(x)+\gamma(x)(\beta-\gamma)(x)+
\gamma(x)^2\big)=\sum_{x\in Q_0}\big(\gamma(x)(\alpha-\gamma)(x)+\gamma(x)\beta(x)\big).
\end{multline*}
The map $q:M\to N$ is a vector bundle. A fiber
$$
q^{-1}(V_1,W_1,\phi)
$$
is the set of all
$$(V,W)\in \Rep(Q,\alpha)\times \Rep(Q,\beta)
$$
such that $V(a)(V_1(ta))\subseteq V_2(ha)$
($(\alpha-\gamma)(ta)\gamma(ha)$ linear constraints)
 for all $a$ and
the restriction of $W(a)$ to $W_1(ta)$ is 
$\phi(ha)V(a)\phi(ta)^{-1}$ ($\gamma(ta)\beta(ha)$ linear constraints)
for all $a$.
The dimension of $q^{-1}(V_1,W_1,\phi)$ is
$$
\dim \Rep(Q,\alpha)+\dim\Rep(Q,\beta)-\sum_{a\in
Q_1}\big((\alpha-\gamma)(ta)\gamma(ha)-\gamma(ta)\beta(ha)\big).
$$
Therefore, 
\begin{multline*}
\dim M=
\dim \mbox{fiber of $q$}+\dim N=\dim \Rep(Q,\alpha)+\dim\Rep(Q,\beta)\\
-\sum_{a\in
Q_1}\big((\alpha-\gamma)(ta)\gamma(ha)-\gamma(ta)\beta(ha)\big)+
\sum_{x\in Q_0} \big(\gamma(x)(\alpha-\gamma)(x)+\gamma(x)\beta(x)\big)=\\
=\dim \Rep(Q,\alpha)+\dim \Rep(Q,\beta)-\langle
\alpha-\gamma,\beta-\gamma\rangle.
\end{multline*}
(Remember that $\langle \alpha,\beta\rangle=0$).
Since $M$ is a fiber bundle over the smooth irreducible
variety  $\Gr{\alpha\choose \alpha-\gamma}\times\Gr{\beta\choose\gamma}$
with smooth irreducible fibers, we have that $M$ is smooth and irreducible.

Consider now the projection $s:M\to \Rep(Q,\alpha)\times \Rep(Q,\beta)$.
Since $\gamma$ is the general rank of a homomorphism $\phi:V\to W$
for  $(V,W)\in D\times \Rep(Q,\beta)$, we see that
$\overline{s(M)}$ contains $D\times \Rep(Q,\beta)$.
Since $\overline{s(M)}$ is irreducible and $D\times \Rep(Q,\beta)$
has codimension 1, we must have
$\overline{s(M)}=D\times \Rep(Q,\beta)$ or $\overline{s(M)}=\Rep(Q,\alpha)\times \Rep(Q,\beta)$.
The latter is impossible, because $\hom(\alpha,\beta)=0$ and $\gamma\neq 0$.
Therefore $\overline{s(M)}=D\times \Rep(Q,\beta)$.
The dimension of a general fiber of $s$ is
\begin{multline*}
\dim M-\dim D\times \Rep(Q,\beta)=\dim M-\dim \Rep(Q,\alpha)-\dim \Rep(Q,\beta)+1=\\
=
-\langle \alpha-\gamma,\beta-\gamma\rangle+1.
\end{multline*}

Choose $(V,W,V_1,W_1,\phi)\in M$ in general position.
Then $(V,W)\in D\times \Rep(Q,\beta)$ is in general position.
Then a general element in $\Hom_Q(V,W)$ has rank $\gamma$. 
The fiber $s^{-1}(V,W)$ is the set of all $\psi\in \Hom_Q(V,W)$ of
rank $\gamma$. Therefore,
$$
\dim \Ext_Q(V,W)=\dim \Hom_Q(V,W)=\dim s^{-1}(V,W)=-\langle \alpha-\gamma,\beta-\gamma\rangle+1.
$$

Suppose that $V_1'\in \Rep(Q,\alpha-\gamma)$ and $W'\in\Rep(Q,\beta)$
are both in general position.
We  know that a general representation of dimension $\beta$
has a $\gamma$-dimensional subrepresentation, because
$\overline{s(M)}=D\times \Rep(Q,\beta)$. So let $W_1'$ be a $\gamma$-dimensional
subrepresentation of $W'$. Define $V'=V_1'\oplus W_1'$ and let
$\phi':V'\to W'$ be the projection $V'\twoheadrightarrow W_1'\subseteq W'$.
Now $(V',W',V_1',W_1',\phi')\in M$.
Moreover $\dim \Ext(V_1',W')$ has the minimum possible value $\ext(\alpha-\gamma,\beta)$.
The set of all $(V',W',V_1',W_1',\phi')\in M$ for which
$\dim \Ext_Q(V_1',W')$ is minimal (or equivalently $\dim \Hom_Q(V_1',W')$
is minimal) is an open dense subset of $M$. Since $(V,W,V_1,W_1,\phi)$ was
assumed to be in general position, we may assume
 $\dim\Ext_Q(V_1,W)=\ext(\alpha-\gamma,\beta)$.
 
From Theorem~\ref{Schofield5.4}
follows that there exists a $\delta\hookrightarrow \alpha-\gamma$ such that
$$
\ext(\alpha-\gamma,\beta)=-\langle\delta,\beta\rangle
$$
Since $V\in D$ is in general position, we know that $V$ is $\sigma$-semi-stable. 
But $V$ has a $\delta$-dimensional subrepresentation $V_2\subseteq V_1\subseteq V$.
Hence
$$
\ext(\alpha-\gamma,\beta)=-\langle \delta,\beta\rangle=\sigma(\delta)\leq 0
$$ 
by semi-stability.
It follows that $\ext(\alpha-\gamma,\beta)=0$ and $\Ext_Q(V_1,W)=0$.
From the long exact sequence of $\Ext$'s, the map
$\Ext_Q(V_1,W)\to \Ext_Q(V_1,W/W_1)$ 
is surjective. So  we
also get $\Ext_Q(V_1,W/W_1)=0$.
Now
\begin{equation}\label{eqstrange}
0=\dim \Ext_Q(V_1,W/W_1)=-\langle\alpha-\gamma,\beta-\gamma\rangle+
\dim \Hom_Q(V_1,W/W_1)
\end{equation}
and so
$$
-\langle\alpha-\gamma,\beta-\gamma\rangle\leq 0
$$
On the other hand, $\Hom_Q(V,W)\neq 0$, so
$$
1\leq \dim \Hom_Q(V,W)=\dim \Ext_Q(V,W)=-\langle \alpha-\gamma,\beta-\gamma\rangle+1
$$
We conclude that
$$
-\langle\alpha-\gamma,\beta-\gamma\rangle=0.
$$
From this follows that $\Hom_Q(V_1,W/W_1)=0$, by (\ref{eqstrange}).
We have an exact sequence
$$
\cdots\to \Hom_Q(V_1,W/W_1)\to \Ext_Q(V_1,W_1)\to \Ext_Q(V_1,W)\to \cdots
$$
Since the outer two are equal to 0, we get $\Ext(V_1,W_1)=0$.
This implies that $\ext(\alpha-\gamma,\gamma)=0$
and $(\alpha-\gamma)\hookrightarrow\alpha$.
If $V'\in \Rep(Q,\alpha)\setminus D$, then
$V'$ has an $(\alpha-\gamma)$-dimensional subrepresentation $V_1'$.
Define $W_1':=V'/V_1'$ and let $W_2'\in \Rep(Q,\beta-\gamma)$ be
any representation of dimension $\beta-\gamma$. Set $W'=W_1'\oplus W_2'$
and let $\phi':V'\twoheadrightarrow V'/V_1'=W_1'\subseteq W'$ be the projection.
We have $(V',W',V_1',W_1',\phi')\in M$, so $(V',W')\in s(M)$ and $V'\in D$.
Contradiction!

\end{proof}

\end{document}